\documentclass[11pt,a4paper]{article}
\usepackage[margin=24mm]{geometry}
\usepackage[T1]{fontenc}
\usepackage{lmodern}
\usepackage{amsthm}
\usepackage{amsmath,amssymb,mathtools,graphicx,xcolor,enumitem,tikz}
\usepackage{microtype,etoolbox}
\usepackage[hypertexnames=false]{hyperref}
\definecolor{linkink}{RGB}{29,66,101}
\hypersetup{colorlinks=true,allcolors=linkink,bookmarksnumbered=true,
  pdftitle={A study guide for Bourgain's 2003 sum-product theorem and 2010 discretized projection theorem},
  pdfauthor={Yudong Liang}}
\allowdisplaybreaks
\DeclareRobustCommand{\paperref}[2]{\hyperref[#1]{#2}}
\newcommand{\annot}[1]{\par\smallskip{\small\noindent\color{blue!65!black}#1\par}\smallskip}

\theoremstyle{plain}
\newtheorem{theorem}{Theorem}[section]
\newtheorem{lemma}[theorem]{Lemma}
\newtheorem{proposition}[theorem]{Proposition}
\newtheorem{corollary}[theorem]{Corollary}
\newtheorem*{theorem*}{Theorem}
\newtheorem*{conjecture*}{Conjecture}
\newtheorem*{claim}{Claim}
\theoremstyle{definition}
\newtheorem*{notation}{Notation}
\theoremstyle{remark}
\newtheorem*{remark}{Remark}
\setlist[itemize]{topsep=.5\baselineskip,itemsep=.2\baselineskip}
\begin{document}

\title{A study guide for Bourgain's 2003 sum-product theorem and 2010 discretized projection theorem}
\author{Yudong Liang}
\date{}
\maketitle
\begin{abstract}
These notes present Bourgain's proofs of the discretized sum-product theorem and its applications. The discussion emphasizes the construction of trees with multiscale porosity from small additive doubling, using either quantitative \paperref{result:1.4}{Freiman theorrem} or the \paperref{result:1.3}{Pl\"unnecke--Ruzsa inequality}, and the amplification of these structures by multiplication. Additive combinatorial tools connect the resulting growth estimates to discretized sum-product problem. The notes also discuss the passage to the exceptional set estimates for orthogonal projections, the exclusion of Borel subrings of intermediate Hausdorff dimension, and Fourier decay for multiplicative convolutions of measures.
\end{abstract}

\begingroup
\small
\makeatletter
\patchcmd{\l@section}{1.0em}{.50em}{}{}
\makeatother
\tableofcontents
\endgroup
\clearpage

\section{Introduction}\label{sec:0}
These notes emphasize the ideas in Bourgain's arguments in \cite{Bourgain2003}\cite{Bourgain2010}. This paper seeks to balance, in its presentation, the ideas behind the proof and its details. Before starting the proof, I will try to give a rough outline and explain the main ideas. Readers who are not interested in the technical details may read only this part and skip the lengthy proof that follows. Absolute constants may change from line to line; sufficiently small parameters may be reselected to absorb preceding losses. Scale-dependent losses and structural hypotheses are recorded explicitly where needed.

In 1966, Erd\H{o}s \cite{ErdosVolkmann1966} considered the following problem:
\begin{quote}
``Does there exist a sub-ring of $\mathbb R$ which is a Borel set and has Hausdorff dimension strictly between $0$ and $1$?''
\end{quote}
This problem is connected to Falconer's distance problem, Furstenberg set problem, and exceptional set estimate for orthogonal projection. For instance, Falconer \cite{Falconer1985} used results on the distance problem to show

\begin{theorem}\label{result:0.1} Let $R\subset\mathbb R$ be a Borel subring and $\dim_H R>\frac12$, then we must have $\dim_H R=1$.
\end{theorem}

\begin{proof} Notice that
\[
\Delta(R\times R)=\sqrt{(R-R)^2+(R-R)^2}\subset\sqrt{R\cap[0,\infty)},
\]
since $R$ is a subring. The map $t\mapsto\sqrt t$ is Lipschitz on each interval $[1/n,n]$, so countable stability of Hausdorff dimension gives $\dim_H\sqrt{R\cap[0,\infty)}\leq\dim_H R$. By Falconer's result, we get

\begin{align*}
\dim_H R&\geq\dim_H\sqrt{R\cap[0,\infty)}\geq\dim_H\Delta(R\times R)\\
&\geq\min\{1,\dim_H(R\times R)-\tfrac12\}\\
&\geq\min\{1,2\dim_H R-\tfrac12\}.
\end{align*}
Writing $s=\dim_H R$, if $\frac12<s<1$, then $\min\{1,2s-\frac12\}>s$, a contradiction. Therefore $s=1$. 
\end{proof}

\begin{remark} If $R\subset\mathbb R$ is a Borel subring with $\dim_H R=\frac12$, then
\[
\dim_H(R\times R)\geq1\quad\&\quad \dim_H\Delta(R\times R)\leq\dim_H R=\tfrac12.
\]
Therefore, if one want to improve Falconer's bound $\dim_H\Delta(E)\geq\frac12$, one have to show such $\frac12$-dim Borel subring doesn't exist.
\end{remark}

In 2003, Bourgain \cite{Bourgain2003} answer this question negatively.

\begin{theorem}\label{result:0.2} A Borel subring of $\mathbb R$ cannot have Hausdorff dimension strictly intermediate between $0$ and $1$.
\end{theorem}

In fact, Bourgain proved a more quantitative version, as considered by Katz and Tao \cite{KatzTao2001}, refer to as the discretized ring conjecture:

\begin{theorem}\label{result:0.3} If $A$ is a $(\delta,\sigma)_1$-set, $0<\sigma<1$ such that $|A|>\delta^{1-\sigma+\varepsilon}$,

then
\begin{equation}
|A+A|+|A\cdot A|>\delta^{1-\sigma-c}\quad(\approx\delta^{-c}|A|)\tag{0.1}\label{eq:0.1}
\end{equation}
for an absolute constant $c=c(\sigma)>0$.

\end{theorem}
Here, a bounded subset $A\subset\mathbb R$ is called a $(\delta,\sigma)_1$-set, provided $A$ is a union of $\delta$-intervals and satisfies
\[
|A\cap I|<\delta^{-\varepsilon}\cdot\delta\cdot\Bigl(\frac r\delta\Bigr)^\sigma
\]
whenever $I\subset\mathbb R$ is an arbitrary interval of size $\delta\leq r\leq1$.

\begin{remark} \paperref{result:0.3}{Theorem~\ref*{result:0.3}} is also a version of discrete sum-product theorem. Let's briefly recall the sum-product phenomenon in Additive Combinatorics.

Suppose $A\subset\mathbb R$ is a finite set. The typical example of $A$ with additive structure: $|A+A|\sim|A|$ is a generalised arithmetic progression (GAP)---and \paperref{result:1.4}{Freiman's theorem} says that the typical example is roughly the only example. Then, what is the typical example of $A$ with multiplicative structure: $|A\cdot A|\sim|A|$? By analogy, this would be the geometric progression or some $d$-dimensional variant of it.

Intuitively, arithmetic and geometric progressions ``look quite different'', so

\centerline{any set cannot "simultaneously look like both".} 

Since a set $A$ with $|A+A|\sim|A|$ implies ``$A$ looks like an arithmetic progression'' by \paperref{result:1.4}{Freiman's theorem}, this heuristic suggests that, $|A\cdot A|\gg|A|$. Erd\H{o}s and Szemer\'edi \cite{ErdosSzemeredi1983} took this heuristic very seriously and conjectured the following in 1983:
\end{remark}

\begin{conjecture*}[Sum-product conjecture] For any $\varepsilon>0$, there is a constant $c_\varepsilon>0$ such that
\[
\max\{|A+A|,|A\cdot A|\}\geq c_\varepsilon|A|^{2-\varepsilon}.
\]
\end{conjecture*}

\begin{theorem*}[Erd\H{o}s \& Szemer\'edi] There are constants $\delta>0$ and $c>0$ such that
\[
\max\{|A+A|,|A\cdot A|\}\geq c|A|^{1+\delta}.
\]
\end{theorem*}

The case $\sigma=\frac12$ has been given special attention in Katz \& Tao \cite{KatzTao2001}, as a corollary, we obtain progress on the Falconer's distance problem and the Furstenberg conjecture.

\begin{corollary}\label{result:0.4} If $K\subset\mathbb R^2$ is a compact set, $\dim_H K\geq1$, then
\[
\dim_H\Delta(K)>\tfrac12+c_0
\]
for an absolute constant $c_0>0$.
\end{corollary}

\begin{corollary}\label{result:0.5} If $K\subset\mathbb R^2$ is a $\frac12$-Furstenberg set (in every direction there is a line $\ell$ such that $\dim_H(K\cap\ell)\geq\frac12$), then there exists an absolute constant $c>0$ such that $\dim_H K>1+c$.
\end{corollary}

Actually, in 2010, Bourgain \cite{Bourgain2010} proved \paperref{result:0.3}{Theorem~\ref*{result:0.3}} again under weaker non-concentration condition.

\begin{notation} We use $|\cdot|$ to denote both Lebesgue measure and cardinality, the meaning will be clear from context.
\end{notation}

\begin{theorem}\label{result:0.6} Given $0<\sigma<1$ and $\kappa>0$, there exist $\varepsilon_0>0$ and $\varepsilon_1>0$ such that if $\delta>0$ is sufficiently small and $A\subset[1,2]$ is a discrete set consisting of $\delta$-separated points, satisfying $|A|=\delta^{-\sigma}$ and
\begin{equation}
|A\cap I|<\rho^\kappa\cdot|A|\tag{$*$}\label{eq:star}
\end{equation}
whenever $I$ is a $\rho$-interval with $\delta<\rho<\delta^{\varepsilon_0}$, then
\begin{equation}
N(A+A,\delta)+N(A\cdot A,\delta)>\delta^{-\varepsilon_1}|A|.\tag{0.2}\label{eq:0.2}
\end{equation}
\end{theorem}

\begin{remark} (1) Note that if $A$ is a finite $\delta$-separated set, then $\widetilde A:=A+B(0,\delta)$ satisfying
\[
\#A=N(A,\delta)=N(\widetilde A,\delta)\quad\&\quad|\widetilde A|=\delta\cdot(\#A).
\]
Therefore, \paperref{result:0.3}{Theorem~\ref*{result:0.3}} is the case $\kappa=\sigma$ of \paperref{result:0.6}{Theorem~\ref*{result:0.6}}. And when $\kappa=\sigma$, the non-concentration assumption $\paperref{eq:star}{(*)}$ is essentially equivalent to $(\delta,\sigma)_1$-set: for $\rho\geq\delta^{\varepsilon_0}$, taking $\rho_0=\delta^{\varepsilon_0}$, we have trivially
\begin{align*}
|A\cap I|&=\sum_{I_{\rho_0}\subset I}|A\cap I_{\rho_0}|
\leq\rho_0^\sigma|A|\cdot N(A\cap I,\rho_0)\\
&\leq\rho_0^\sigma\cdot\delta^{-\sigma}\cdot\frac\rho{\rho_0}
=\Bigl(\frac\rho{\rho_0}\Bigr)^{1-\sigma}\Bigl(\frac\rho\delta\Bigr)^\sigma
\leq\delta^{-\varepsilon_0(1-\sigma)}\Bigl(\frac\rho\delta\Bigr)^\sigma.
\end{align*}
(2) Hence, $\paperref{eq:star}{(*)}$ play a role when $\kappa\ll\sigma$. In fact, $\paperref{eq:star}{(*)}$ can be viewed as some ``two-ends'' condition introduced by Wolff \cite{Wolff1995}.
\begin{itemize}
\item For intervals of size $1$, there is no requirement.
\item For intervals of size $\delta$, or even $(\delta|A|)^{1/(1-\kappa)}=\delta^{(1-\sigma)/(1-\kappa)}$, the requirement is trivial:
\[|A\cap I|\leq\frac\rho\delta\leq\rho^\kappa|A|.\]
\item The most important scale, where $A$ is required not to concentrate, is $\sim\delta^{\varepsilon_0}$, and the requirement is
\[|A\cap I(a,\delta^{\varepsilon_0})|\leq\delta^{\varepsilon_0\kappa}|A|.\]
In other word, $A$ is $(\varepsilon_0,\varepsilon_0\kappa)$ two-ends.
\item Conversely, if $A$ is $(\varepsilon_1,\varepsilon_2)$ two-ends, $(0<\varepsilon_2<\varepsilon_1<1)$, i.e.,
\[|A\cap I(a,\delta^{\varepsilon_1})|\leq\delta^{\varepsilon_2}|A|,\]

then $A$ satisfies $\paperref{eq:star}{(*)}$ with $(\kappa,\varepsilon_0)$ replaced by $(\varepsilon_2,\varepsilon_1)$: for any $\rho$-interval $I$, $\delta<\rho<\delta^{\varepsilon_1}$,
\[|A\cap I|\leq|A\cap I(c_I,\delta^{\varepsilon_1})|\leq\delta^{\varepsilon_2}|A|<\rho^{\varepsilon_2}|A|.\]
\end{itemize}
\end{remark}

\paperref{result:0.6}{Theorem~\ref*{result:0.6}} will be a consequence of the \paperref{result:0.7}{following theorem}, which is more convenient in applications.

\begin{theorem}\label{result:0.7} Given $0<\sigma<1$ and $\kappa>0$, there exist $\varepsilon_0>0$ and $\varepsilon_1>0$ such that the following holds for $\delta>0$ sufficiently small.

Let $\mu\in\mathcal P([0,1])$ be a probability measure on $[0,1]$ satisfying
\[\mu(I)<\rho^\kappa\]
whenever $I$ is a $\rho$-interval, $\delta<\rho<\delta^{\varepsilon_0}$. Let $A\subset[1,2]$ be a discrete set consisting of $\delta$-separated points with $|A|\leq\delta^{-\sigma}$ such that also
\[|A\cap I|<\rho^\kappa|A|\]
whenever $I$ is a $\rho$-interval, $\delta<\rho<\delta^{\varepsilon_0}$. Then there exists $x\in\operatorname{supp}(\mu)$ such that
\begin{equation}N(A+xA,\delta)>\delta^{-\varepsilon_1}|A|.\tag{0.3}\label{eq:0.3}\end{equation}
\end{theorem}

\begin{remark} Although \paperref{eq:0.3}{(\ref*{eq:0.3})} only holds for a single point in conclusion, we may

deduce that \paperref{eq:0.3}{(\ref*{eq:0.3})} holds for all $x\in D\subset\operatorname{supp}(\mu)$, where
\[\mu(D)>1-\delta^{\frac1{10}\varepsilon_0\kappa}.\]
Indeed, let $D=\{x\in\operatorname{supp}(\mu):\text{\paperref{eq:0.3}{(\ref*{eq:0.3})} holds for }x\}\neq\varnothing$, restricting $\mu$ to $\operatorname{supp}(\mu)\setminus D$, taking
\[\nu=\frac1{\mu(D^c)}\mu|_{D^c}.\]
Then, $\nu\in\mathcal P([0,1])$, and for any $\rho$-interval $I$, $\delta<\rho<\delta^{\varepsilon_0}$,
\[\nu(I)\leq\frac1{\mu(D^c)}\mu(I)\leq\frac1{\mu(D^c)}\rho^\kappa.\]
Assume $\mu(D)\leq1-\delta^{\frac1{10}\varepsilon_0\kappa}$, then
\[\mu(D^c)>\delta^{\frac1{10}\varepsilon_0\kappa}\Rightarrow\nu(I)\leq(\delta^{-\varepsilon_0})^{\frac1{10}\kappa}\rho^\kappa<\rho^{\frac9{10}\kappa}.\]
Apply \paperref{result:0.7}{Theorem~\ref*{result:0.7}} to $\sigma,\frac9{10}\kappa,\nu$ and $A$, we find some $x\in\operatorname{supp}(\nu)\subset\operatorname{supp}(\mu)\setminus D$, such that \paperref{eq:0.3}{(\ref*{eq:0.3})} holds for $x$, which contradicts to the choice of $D$.
\end{remark}

It turns out that there is an interpretation of \paperref{result:0.7}{Theorem~\ref*{result:0.7}} as a Marstrand-type projection theorem: the dimension of projections satisfies a non-trivial lower bound outside a very sparse set of exceptional directions.

\begin{theorem}\label{result:0.8} Given $0<\alpha<2$, $\beta>0$ and $\kappa>0$, there exist $\tau_0>0$ and $\eta>\alpha/2$ such that the following holds.

Let $\mu_1\in\mathcal P(S^1)$ such that
\[\max_\theta\mu_1([\theta-\rho,\theta+\rho])<\rho^\kappa\]
for $\delta<\rho<\delta^{\tau_0}$. Let $\delta>0$ be chosen sufficiently small and $\mathcal A\subset[1,2]\times[1,2]$ be a union of $\delta$-squares satisfying $|\mathcal A|=\delta^{2-\alpha}$ and
\[\max_x|\mathcal A\cap B(x,\rho)|<\rho^\beta|\mathcal A|\]
for $\delta<\rho<\delta^{\tau_0}$. Then, there exists $\theta\in\operatorname{supp}(\mu_1)$ such that
\begin{equation}|\pi_\theta(\mathcal A)|>\delta^{1-\eta},\tag{0.4}\label{eq:0.4}\end{equation}
where $\pi_\theta$ denotes the orthogonal projections on the line $y=(\tan\theta)x$.
\end{theorem}

\begin{theorem}\label{result:0.9} Given $0<\alpha<2$ and $\kappa>0$, there exists $\eta>\alpha/2$ such that, if $\mathcal A\subset\mathbb R\times\mathbb R$ is an analytic set (for example, a Borel set) with $\dim_H\mathcal A>\alpha$, then
\[\dim_H\pi_\theta(\mathcal A)\geq\eta\]
for all $\theta\in S^1$ except in an exceptional set $E$ satisfying
\[\dim_H E\leq\kappa.\]
Equivalently,
\begin{equation}\dim_H\{\theta\in S^1:\dim_H\pi_\theta(\mathcal A)<\eta\}\leq\kappa.\tag{0.5}\label{eq:0.5}\end{equation}
\end{theorem}

\begin{remark} In particular, for analytic $\mathcal A$, we have
\[\dim_H\{\theta\in S^1:\dim_H\pi_\theta(\mathcal A)<\tfrac12\dim_H\mathcal A\}=0.\]
Since, for any $\kappa>0$, suppose $\dim_H\mathcal A=\alpha$, we get
\begin{align*}
&\dim_H\{\theta\in S^1:\dim_H\pi_\theta(\mathcal A)<\tfrac12\dim_H\mathcal A\}\\
&\qquad\leq\sup_n\dim_H\{\theta\in S^1:\dim_H\pi_\theta(\mathcal A)<\tfrac12(\alpha-\tfrac1n)\}.
\end{align*}
(Let $\alpha_n=\alpha-1/n$, applying \paperref{result:0.9}{Theorem~\ref*{result:0.9}} with $\alpha_n,\kappa$ and $\mathcal A$, we may choose $\eta=\eta(\kappa,\alpha,n)>\alpha_n/2=\frac12(\alpha-\frac1n)$, such that the following:)
\[\leq\sup_n\dim_H\{\theta\in S^1:\dim_H\pi_\theta(\mathcal A)<\eta\}\leq\kappa\to0.\]
Moreover, in the proof of \paperref{result:0.9}{Theorem~\ref*{result:0.9}}, we can obtain the uniform boundedness of $\eta$: for any closed interval $[a,b]\subset (0,2)$, $\min_{\alpha\in[a,b]}\eta(\alpha,\kappa)>\alpha/2$. And then, we can replace the strictly less-than sign in the above exceptional set estimate with a less-than-or-equal-to sign:
\[
\dim_H\{\theta\in S^1:\dim_H\pi_\theta(\mathcal A)\le\tfrac12\dim_H\mathcal A\}=0.
\]
\end{remark}

In higher dimension, we also have a similar statements, replacing $\alpha/2$ by $\alpha/d$ and an appropriate assumption on ``directions''. Let $\mathbb P^k=\operatorname{Gr}(k+1,1)$ be the $k$-dimensional (real) projective space and for $\bar y\in\mathbb P^k$, $\rho>0$, let
\[V(y^\perp,\rho)=\left\{\bar x\in\mathbb P^k:\frac{|\langle x,y\rangle|}{|x|\cdot|y|}<\rho\right\}\]
be a $\rho$-neighborhood of $y^\perp$.
\annot{The quotient is $\cos\angle(x,y)$; \paperref{fig:projective}{the diagram} indicates angle $\approx\rho$.}
\begin{center}\phantomsection\label{fig:projective}\includegraphics[width=.45\linewidth]{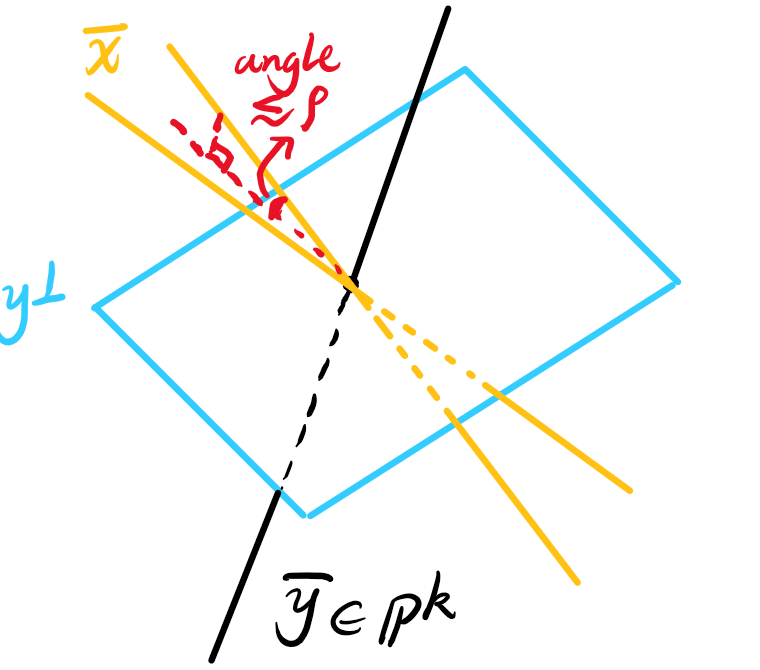}\end{center}

\begin{theorem}\label{result:0.10} Given $0<\alpha<d$, $d\geq2$ and $\kappa>0$, there exist $\tau_0>0$, $\varepsilon>0$ and $\eta>\alpha/d$ such that the following holds for sufficiently small $\delta>0$.

Let $\mu_1\in\mathcal P(\mathbb P^{d-1})$ such that
\[\max_{\bar y}\mu_1(V(y^\perp,\rho))<\rho^\kappa\]
for all $\delta<\rho<\delta^{\tau_0}$. Let $\overline E\subset[0,1]^d$ be a union of $\delta$-balls with $|\overline E|>\delta^{d-\alpha}$, and satisfying a non-concentration condition
\[\max_x|\overline E\cap B(x,\rho)|<\rho^\kappa|\overline E|\]
for all $\delta<\rho<\delta^{\tau_0}$. Then there exist $D\subset\mathbb P^{d-1}$ and $\overline E'\subset\overline E$ with the following properties:
\begin{equation}\begin{gathered}
\mu_1(D)>1-\delta^\varepsilon,\\
|\overline E'|>\delta^\varepsilon|\overline E|,\\
\&\quad |\pi_\theta(\overline E'')|>\delta^{1-\eta}
\end{gathered}\tag{0.6}\label{eq:0.6}\end{equation}
whenever $\theta\in D$ and $\overline E''\subset\overline E'$ satisfies $|\overline E''|>\delta^{2\varepsilon}|\overline E|$.
\end{theorem}

\begin{theorem}\label{result:0.11} Given $0<\alpha<d$, $d\geq2$ and $\kappa>0$, there exists $\eta>\alpha/d$ such that, if $\overline E\subset\mathbb R^d$ is an analytic set (for example, a Borel set) with $\dim_H\overline E>\alpha$, then
\[\dim_H\pi_\theta(\overline E)\geq\eta\]

for all $\theta\in\mathbb P^{d-1}$ except in an exceptional set $F$ satisfying
\[\dim_H F\leq1\cdot(d-1)-1+\kappa=d-2+\kappa.\]
\annot{$1\cdot(d-1)-1$: the dimension of $y^\perp\cap\mathbb P^{d-1}=y^\perp\cap S^{d-1}$.}
Equivalently,
\begin{equation}\dim_H\{\theta\in\mathbb P^{d-1}:\dim_H\pi_\theta(\overline E)<\eta\}\leq d-2+\kappa.\tag{0.7}\label{eq:0.7}\end{equation}
\end{theorem}

As we know, Bourgain, Glibichuk and Konyagin \cite{BGK2006} proved
\[
\begin{gathered}\text{sum-product theorem}\\\text{in }\mathbb F_p\end{gathered}
\quad\Longrightarrow\quad
\begin{gathered}\text{the multilinear exponential}\\\text{sum estimates in}\\\text{finite fields}.\end{gathered}
\]
As a counterpart, we have
\[
\begin{gathered}\text{discretized sum-product}\\\text{theorem in }\mathbb R\end{gathered}
\quad\Longrightarrow\quad
\begin{gathered}\text{estimates on multiplicative}\\\text{convolutions of measures}\\\text{on the real line}.\end{gathered}
\]

\begin{theorem}\label{result:0.12} For all $\kappa>0$ and $C>1$, there exist $k\in\mathbb Z_+$ and $C'>1$ such that the following holds.

Let $\mu\in\mathcal P([\frac12,1])$ satisfying for $\delta<\rho<1$,
\[\max_a\mu(B(a,\rho))<C\rho^\kappa.\]
Denoted by $\nu$ the image measure of $\mu\times\cdots\times\mu$ ($k^2$-fold) under the map
\[
(x_1,\ldots,x_{k^2})\mapsto(x_1\cdots x_k)+(x_{k+1}\cdots x_{2k})+\cdots+(x_{k^2-k+1}\cdots x_{k^2}),
\]
in other word, $\nu$ is the $k$-fold convolutions of measures which are $k$-fold multiplications of $\mu$. Then for all $\delta<\rho<1$,
\begin{equation}\max_a\nu(B(a,\rho))<C'\rho.\tag{0.8}\label{eq:0.8}\end{equation}
\end{theorem}

\begin{remark} \paperref{result:0.12}{Theorem~\ref*{result:0.12}} is equivalent to the statement that if $\mu$ as above, then
\begin{equation}
\left|\bigl(\psi[\underbrace{\mu\otimes\cdots\otimes\mu}_{k}]\bigr)^{\wedge}(\xi)\right|\lesssim|\xi|^{-\tau}\tag{0.9}\label{eq:0.9}
\end{equation}
for all $1\leq|\xi|\leq\delta^{-1}$, some $k=k(\kappa)\in\mathbb Z_+$ and $\tau=\tau(\kappa)>0$, where
\[\psi:(x_1,\ldots,x_k)\mapsto x_1\cdots x_k.\]
This equivalence can be found in the proof.
\end{remark}

The question as to what is the optimal relation between $\kappa$ and $k$ is interesting and remains to be clarified. One should expect that an inequality of the form \paperref{eq:0.9}{(\ref*{eq:0.9})} may hold if $k>1/\kappa$. (The case $k=2$ have been proved in the following.) And in 2009, Bourgain had shown that
the assumption $k>1/\kappa$ is sufficient in the prime field setting. The answer requires more precise version of \paperref{result:0.6}{Theorem~\ref*{result:0.6}}.

\textbf{P.S.} This question has been answered by Saxc\'e, Shmerkin and Orponen \cite{OrponenSaxceShmerkin}. In fact, they proved a more general result under some weaker assumption on $\mu$.

\paperref{eq:0.8}{(\ref*{eq:0.8})} is a generalization of the following classical and elementary inequality for bilinear convolution.

\begin{theorem}\label{result:0.13} Let $\mu,\nu\in\mathcal P([0,1])$ satisfying, for all $\delta\leq\rho\leq1$, the conditions
\[\max_a\mu(B(a,\rho))\lesssim\rho^\alpha\quad\&\quad\max_a\nu(B(a,\rho))\lesssim\rho^\beta.\]
Assume $\alpha+\beta>1$. (Taking $\mu=\nu$, $\alpha=\beta=\kappa$, we have $2>1/\kappa$.)
Let $\psi:[0,1]\times[0,1]\to[0,1]$, $(x,y)\mapsto x\cdot y$ and $\psi[\mu\otimes\nu]$ be the image measure. Then
\begin{equation}
\left|(\psi[\mu\otimes\nu])^{\wedge}(\xi)\right|\lesssim|\xi|^{\frac{1-\alpha-\beta}{2}}\tag{0.10}\label{eq:0.10}
\end{equation}
for all $1\leq|\xi|\leq\delta^{-1}$.
\end{theorem}

I wrote these notes in early 2026 while reading Bourgain’s paper and preparing lecture materials for a seminar. Recently, I utilized chatGPT to assist in performing OCR on the manuscript and converting it into a LaTeX document. Furthermore, much of the explanatory text in \paperref{sec:8}{Sections~8} and \paperref{sec:9}{9} was developed through discussions with GPT; in my original manuscript, my understanding of how to translate the sum-product theorem into Fourier decay estimates for measures was not yet fully formed.

\section{Additive combinatorial tools}\label{sec:1}

All finite sets in this section are understood to be nonempty subsets of an abelian group.

\textbf{(1) Sumset estimates}

\begin{itemize}\item Let $A,B$ be finite sets, $|A|\leq|B|$. Then
\begin{equation}\frac{|A+A|}{|A|}\leq\left(\frac{|A+B|}{|A|}\right)^3.\tag{1.1}\label{eq:1.1}\end{equation}\end{itemize}

\begin{lemma}\label{result:1.1} Let $A,B\subset\mathbb R$ and $\delta>0$ such that $N(A,\delta)\leq N(B,\delta)$. Then
\begin{equation}\frac{N(A+A,\delta)}{N(A,\delta)}\leq\left(\frac{N(A+B,\delta)}{N(A,\delta)}\right)^3.\tag{1.2}\label{eq:1.2}\end{equation}
\end{lemma}

\begin{itemize}\item \paperref{result:1.2}{Ruzsa's triangle inequality}\end{itemize}
\begin{lemma}\label{result:1.2} For any three finite sets $A,B,C$,
\begin{equation}|A-B|\leq\frac{|A-C|\cdot|B-C|}{|C|}.\tag{1.3}\label{eq:1.3}\end{equation}
Furthermore, a ``plus sign'' version of \paperref{eq:1.3}{(\ref*{eq:1.3})} also holds. This argument is due to Todd Cochrane:
\begin{equation}|A+B|\leq\frac{|A-C|\cdot|B-C|}{|C|}.\tag{1.4}\label{eq:1.4}\end{equation}
\end{lemma}

\begin{itemize}\item \paperref{result:1.3}{Pl\"unnecke--Ruzsa inequality}\end{itemize}
\begin{lemma}\label{result:1.3} Suppose $A,B$ are finite sets, $|A+B|\leq K|A|$, then we have

\begin{equation}|nB-mB|\leq K^{n+m}|A|\tag{1.5}\label{eq:1.5}\end{equation}
for all $n,m\geq1$. In particular, if $|A+A|\leq K|A|$, then
\begin{equation}|nA|\leq K^n|A|.\tag{1.6}\label{eq:1.6}\end{equation}
\centerline{``Having small sumset implies having small $n$-sumset.''}
\end{lemma}

\textbf{(2) \paperref{result:1.4}{Freiman's theorem} (quantitative version)}

\begin{theorem}\label{result:1.4} Let $A\subset\mathbb Z$ be a finite set such that
\[|A+A|\leq K|A|.\]
Then $A$ is contained in a translation of a proper $d$-dimensional progression $P$,
\[P=[h]\cdot\xi:=\{b_1\xi_1+\cdots+b_d\xi_d:0<b_i\leq h_i;\ i=1,\ldots,d\},\]
where $\xi\in\mathbb Z^d$ and $h\in\mathbb N^d$ satisfying
\begin{equation}
\begin{gathered}d<K,\\
\&\quad|P|=\prod_{i=1}^d h_i<\exp\bigl(K^2(\log K)^3\bigr)|A|.
\end{gathered}\tag{1.7}\label{eq:1.7}
\end{equation}
\centerline{``A set having small sumset is a dense subset of a low-rank G.A.P.''}
\end{theorem}

\begin{corollary}\label{result:1.5} Let $A\subset\mathbb R$ be bounded and nonempty, and let $0<\delta<1$, such that

\[N(A+A,\delta)<KN(A,\delta).\]
Assume
\begin{equation}\sigma\sim K[\delta N(A,\delta)]^{K^{-3}}.\tag{1.8}\label{eq:1.8}\end{equation}
Then there is $q\in\mathbb Z_+$ such that $\sigma/q>\delta$ and
\begin{equation}A\subset\{\xi_0+a/q:a\in\mathbb Z\}+I(0,\sigma/q).\tag{1.9}\label{eq:1.9}\end{equation}
\end{corollary}

\begin{proof} If $\sigma\geq1$, the conclusion follows with $q=1$, so consider $\sigma<1$. Round $A$ to the lattice $\delta\mathbb Z$: choose a finite $A'\subset\delta\mathbb Z$ such that $A\subset A'+I(0,\delta)$, $|A'|\asymp N(A,\delta)$ and $|A'+A'|\lesssim K|A'|$.

Apply \paperref{result:1.4}{Theorem~\ref*{result:1.4}} to $\delta^{-1}A'\subset\mathbb Z$ and scale the resulting progression back by $\delta$. Thus $A'\subset\xi_0+P=\xi_0+[h]\cdot\xi$, now with real generators $\xi_i$; no bound on their magnitudes is needed. Fix $0<\kappa<1$ and consider on $\mathbb T^d=(\mathbb R/\mathbb Z)^d$ the norm
\[\lvert\!\lvert\!\lvert(\psi_1,\ldots,\psi_d)\rvert\!\rvert\!\rvert=\max_{1\leq i\leq d}(\kappa^{-1}h_i\|\psi_i\|).\]
Thus
\[\operatorname{mes}_{\mathbb T^d}[\lvert\!\lvert\!\lvert\psi\rvert\!\rvert\!\rvert<1]=\left|\prod_{i=1}^d[0,\kappa h_i^{-1})\right|=\kappa^d\prod_{i=1}^d h_i^{-1}.\]
Therefore, consider the sequence $q\xi\pmod1\in\mathbb T^d$,
\[q=0,1,\ldots,\kappa^{-d}\prod_{i=1}^d h_i,\]
by pigeonholing principle, there exists $q\in\mathbb Z_+$, $0<q\leq\kappa^{-d}\prod h_i$ such that
\[\lvert\!\lvert\!\lvert(q\xi_1,\ldots,q\xi_d)\rvert\!\rvert\!\rvert<1.\]
Hence,
\[|q\xi_i-a_i|<\kappa h_i^{-1}\quad\text{for some }a_i\in\mathbb Z,\ i=1,\ldots,d.\]

Thus,
\begin{align*}
\left|(\xi_0+\textstyle\sum b_i\xi_i)-\left(\xi_0+\frac{\sum b_ia_i}q\right)\right|
&\leq\sum|h_i|\left|\xi_i-\frac{a_i}q\right|\\
&\leq\sum_{i=1}^d h_i\frac\kappa{qh_i}=\frac{d\kappa}q.
\end{align*}
\begin{align*}
\Rightarrow\quad A'\subset\xi_0+P&\subset\{\xi_0+a/q:a\in\mathbb Z\}+I(0,d\kappa/q),\\
\Rightarrow\quad A\subset A'+I(0,\delta)&\subset\{\xi_0+a/q:a\in\mathbb Z\}+I(0,d\kappa/q+\delta).
\end{align*}
We impose the condition $\delta<d\kappa/q$, which will be satisfied if
\[\delta<d\kappa/(\kappa^{-d}\textstyle\prod h_i)\Longleftrightarrow\delta\prod h_i<d\kappa^{d+1}.\]
By \paperref{eq:1.7}{(\ref*{eq:1.7})},
\[\delta\prod h_i<\delta e^{K^2(\log K)^3}|A'|<e^{K^3}[\delta N(A,\delta)]\]
and $d\kappa^{d+1}>\kappa^{K^3}$, recall $d<K$.
Thus we need
\[e^{K^3}[\delta N(A,\delta)]<\kappa^{K^3}\Rightarrow\kappa\gtrsim[\delta N(A,\delta)]^{K^{-3}}.\]
Finally, taking $\sigma$ to be a sufficiently large absolute multiple of $K\kappa$ absorbs both $d\kappa/q$ and the rounding error $\delta$, and gives $\sigma/q>\delta$. This is consistent with \paperref{eq:1.8}{(\ref*{eq:1.8})}, and completes the proof. 
\end{proof}

\textbf{(3) \paperref{result:1.6}{Balog--Szemer\'edi--Gowers theorem}.}

\begin{theorem}\label{result:1.6} There exists an absolute constant $C$ such that the following holds.

Let $A,B$ be finite sets of an abelian group $Z$ and $G\subset A\times B$ such that, for some $K>1$,
\[|G|\geq\frac1K|A|\cdot|B|\]
and
\[|A\mathbin{\overset G+}B|:=|\{a+b:(a,b)\in G\}|\leq K|A|^{1/2}\cdot|B|^{1/2}.\]
Then there exist $A'\subset A$, $B'\subset B$ satisfying
\begin{gather*}
|A'|>K^{-C}|A|,\qquad |B'|>K^{-C}|B|,\\
|A'+B'|<K^C|A|^{1/2}|B|^{1/2},\\
\&\quad|G\cap(A'\times B')|>K^{-C}|A|\cdot|B|.
\end{gather*}
An equivalent energy formulation is the following. Suppose
\begin{align*}
E(A,B)&:=|\{(a,a',b,b')\in A\times A\times B\times B:a+b=a'+b'\}|\\
&\geq\frac1K|A|^{3/2}\cdot|B|^{3/2}.
\end{align*}
Then there are $A'\subset A$ and $B'\subset B$ with
\[
|A'|>K^{-C}|A|,\qquad |B'|>K^{-C}|B|,\qquad
|A'+B'|<K^C|A|^{1/2}|B|^{1/2}.
\]
\end{theorem}

``The \paperref{result:1.6}{Balog--Szemer\'edi--Gowers Theorem} allows us to replace a partial sumset by an honest sumset on large subsets. In other words, we can replace the dense subset $G\subset A\times B$ by a dense product set $A'\times B'$.''

\section{Tree structure with multiscale porosity}\label{sec:2}

We will construct a tree structure on $A$ using ``$A$ has small sumset'', which has multi-scale porosity. That is, at each level of the tree, it is a union of subsets contained in disjoint intervals of a certain size which is very small w.r.t. their mutual separation.

\subsection{Induce the porosity property at a single scale}\label{sec:2.1}
\textcolor{red}{
The quantitative version, \paperref{result:1.5}{Corollary~\ref*{result:1.5}}, has shown that a set $A$ with additive structure induces a separated property at a larger scale. However, a sparse structure does not necessarily imply that the pores are small. The problem is that the double constant $K$ is too large. So our aim is to relax the restriction on the doubling constant $K$ in \paperref{eq:1.8}{(\ref*{eq:1.8})}. More precisely, we may discretize $A$ at some slightly large scale, despite it is not contained in some arithmetic progression, it will have smaller and more flexible ``ratio'', $\sigma$, at this scale.} \textcolor{blue}{(ratio $=$ size/separation).}

First, we show that for the most pieces of $A$, their doubling constants are well-behaved.

\begin{lemma}\label{result:2.1} Assume $N(A+A,\delta)=:N(2A,\delta)<KN(A,\delta)$. Let $\eta>\delta$ and
\[
 A=\bigcup_{j\in\mathcal G}A_j,
\]
where $A_j=A\cap I_j\ne\varnothing$ for $j\in\mathcal G$ and $\{I_j\}$ are $\eta$-interval.
Denoted by $j^*\in\mathcal G$ satisfying
\[
 N(A_{j^*},\delta)=\max_j N(A_j,\delta)
\]
and let
\[
 \mathcal G_1=\{j\in\mathcal G:N(A_j,\delta)>\tau K^{-1}N(A_j+A_{j^*},\delta)\}.
\]

Then
\begin{equation}
 \sum_{j\notin\mathcal G_1}N(A_j,\delta)<4\tau N(A,\delta)\tag{2.1}\label{eq:2.1}
\end{equation}
and for each $j\in\mathcal G_1$,
\begin{equation}
 \frac{N(2A_j,\delta)}{N(A_j,\delta)}<\left(\frac K\tau\right)^3,\tag{2.2}\label{eq:2.2}
\end{equation}
where $0<\tau<1$ is a given parameter.
\end{lemma}

\begin{proof} Since $\{A_j+A_{j^*}\}_{j\in\mathcal G}$ are essentially disjoint, we have
\[
\begin{aligned}
 \sum_{j\notin\mathcal G_1}N(A_j,\delta)
 &\le\tau K^{-1}\sum_{j\in\mathcal G}N(A_j+A_{j^*},\delta)
 \lesssim\tau K^{-1}N(2A,\delta)\\
 &<4\tau N(A,\delta).
\end{aligned}
\]
To prove \paperref{eq:2.2}{(\ref*{eq:2.2})}, we use the sumset estimate \paperref{eq:1.2}{(\ref*{eq:1.2})}:
\[
 \frac{N(2A_j,\delta)}{N(A_j,\delta)}
 \le\left(\frac{N(A_j+A_{j^*},\delta)}{N(A_j,\delta)}\right)^3
 <\left(\frac K\tau\right)^3.
\]
\end{proof}
Next, assume $A\subset\mathbb R$ bounded with
\begin{equation}
 N(2A,\delta)<KN(A,\delta)\tag{2.3}\label{eq:2.3}
\end{equation}
\begin{equation}
 \text{\&}\quad \delta N(A,\delta)<(\kappa K^{-3})^t\quad(<\kappa),\tag{2.4}\label{eq:2.4}
\end{equation}

where $\kappa<10^{-3}$ and
\begin{equation}
 t\sim\frac{\log K}{\log\log K}\gg1.\tag{2.5}\label{eq:2.5}
\end{equation}
We will construct a initial tree structure of subsets of $A$.

\annot{\textbf{Clarification.} Our aim is to reduce the doubling constant by exploiting the multiplicative property of it at various scales. Until it suffices to give a small ``ratio'' in \paperref{result:1.5}{Corollary~\ref*{result:1.5}}. Roughly speaking, we extract from $A$ a tree consisting of "good" scales. The pieces at each level of the tree possess a consistent double constant, and the product of double constants across all levels is bounded by the initial double constant $K$. Consequently, provided the tree has a sufficient number of levels, we can always find a scale for which the corresponding double constant at that level is small enough.}

Take $\delta_0=\delta$. Let $\delta_0<\delta_1<1$ be such that
\begin{equation}
 \max_{|I|=\delta_1}\delta_0N(A\cap I,\delta_0)\sim\kappa\delta_1
 \quad\Longleftrightarrow\quad
 \max_{|I|=\delta_1}\frac{|A\cap I|}{|I|}\sim\kappa.\tag{2.6}\label{eq:2.6}
\end{equation}
This is well-defined, since taking $\delta_1=1$ and by \paperref{eq:2.4}{(\ref*{eq:2.4})},
\[
 \max_{|I|=\delta_1}\delta_0N(A\cap I,\delta_0)
 =\delta N(A,\delta)<\kappa\cdot1.
\]
On the other hand, taking $\delta_1=\delta_0$, we have trivially
\[
 \max_{|I|=\delta_0}\delta_0N(A\cap I,\delta_0)=\delta_0>\kappa\cdot\delta_0.
\]

Apply \paperref{result:2.1}{Lemma~\ref*{result:2.1}} with $(\delta,\eta)=(\delta_0,\delta_1)$ and $\tau=1/10$, thus if $A'=\bigcup_{j\in\mathcal G_1}A_j$, we get
\[
 \paperref{eq:2.1}{(\ref*{eq:2.1})}\ \Longrightarrow\ N(A',\delta)=N(A,\delta)-\sum_{j\notin\mathcal G_1}N(A_j,\delta)>\frac12N(A,\delta).
\]
And for each $j\in\mathcal G_1$, by definition of $\delta_1$,
\[
 \frac1{10K}\kappa\frac{\delta_1}{\delta_0}<N(A_j,\delta_0)\lesssim\kappa\frac{\delta_1}{\delta_0}.
\]
(The lower bound comes from, suppose $A_k$ attend the maximum \paperref{eq:2.6}{(\ref*{eq:2.6})}, then by the definition of $j^*$, we have
\[
 N(A_j,\delta_0)>\frac1{10K}N(A_j+A_{j^*},\delta)
 \ge\frac1{10K}N(A_k,\delta_0)
 =\frac1{10K}\kappa\frac{\delta_1}{\delta_0}.
\]
)
\[
 \paperref{eq:2.2}{(\ref*{eq:2.2})}\ \Longrightarrow\ 1\le\frac{N(2A_j,\delta_0)}{N(A_j,\delta_0)}\lesssim K^3.
\]

By this dyadic pigeonholing, we may extract a further subset $A^{(1)}\subset A'$, obtained by retaining equally many $\delta_0$-intervals in the selected intersections $A\cap I^{(\delta_1)}$, such that for some
\begin{equation}
 \frac\kappa{10K}\le\kappa_1\le\kappa
 \quad\text{\&}\quad 1\le R_1\lesssim K^3,\tag{2.7}\label{eq:2.7}
\end{equation}
we have
\begin{equation}
 N(A^{(1)}\cap I^{(\delta_1)},\delta_0)\sim\kappa_1\frac{\delta_1}{\delta_0}\tag{2.8}\label{eq:2.8}
\end{equation}

\begin{equation}
 \text{\&}\quad
 \frac{N(2(A^{(1)}\cap I^{(\delta_1)}),\delta_0)}{N(A^{(1)}\cap I^{(\delta_1)},\delta_0)}\sim R_1\tag{2.9}\label{eq:2.9}
\end{equation}
for $I^{(\delta_1)}$ such that $A^{(1)}\cap I^{(\delta_1)}\ne\varnothing$. Moreover,
\begin{equation}
 N(A^{(1)},\delta_0)\gtrsim(\log K)^{-2}N(A,\delta_0).\tag{2.10}\label{eq:2.10}
\end{equation}
Next let $\delta_1<\delta_2<1$ satisfying
\[
 \max_{|I|=\delta_2}\delta_1N(A^{(1)}\cap I,\delta_1)\sim\kappa\delta_2.
\]
This also make sense, since
\begin{itemize}
\item $\delta_2=1$:
\[
 \delta_0N(A^{(1)},\delta_0)
 \sim\delta_0\sum_{|I|=\delta_1}N(A^{(1)}\cap I,\delta_0)
 \overset{\paperref{eq:2.8}{(\ref*{eq:2.8})}}\sim\delta_0\cdot\kappa_1\frac{\delta_1}{\delta_0}N(A^{(1)},\delta_1)
\]
\[
\begin{aligned}
 \Longrightarrow\quad\delta_1N(A^{(1)},\delta_1)
 &\sim\frac{\delta_0}{\kappa_1}N(A^{(1)},\delta_0)\\
 &\overset{\paperref{eq:2.7}{(\ref*{eq:2.7})}}\lesssim10K\frac1\kappa\delta_0N(A,\delta_0)\\
 &\overset{\paperref{eq:2.4}{(\ref*{eq:2.4})}}<(\kappa K^{-3})^{t-1}<\kappa\cdot1.
\end{aligned}
\]
\item $\delta_2=\delta_1$:
\[
 \delta_1N(A^{(1)}\cap I^{(\delta_1)},\delta_1)=\delta_1>\kappa\cdot\delta_1.
\]
\end{itemize}
And the new doubling constant is
\[
 \frac{N(2A^{(1)},\delta_1)}{N(A^{(1)},\delta_1)}
 \sim\kappa_1\frac{\delta_1}{\delta_0}\frac{N(2A^{(1)},\delta_1)}{N(A^{(1)},\delta_0)}.
\]

\[
\begin{aligned}
 \frac{N(2A^{(1)},\delta_1)}{N(A^{(1)},\delta_1)}
 &\sim\frac{N(2A^{(1)},\delta_1)\,N(A^{(1)}\cap I^{(\delta_1)},\delta_0)}{N(A^{(1)},\delta_0)}\\
 &\le\frac{N(2A^{(1)},\delta_1)\,N((2A^{(1)})\cap I^{(\delta_1)},\delta_0)}{N(A^{(1)},\delta_0)}\\
 &<C\frac{N(2A^{(1)},\delta_0)}{N(A^{(1)},\delta_0)}\\
 &\overset{\paperref{eq:2.10}{(\ref*{eq:2.10})}}<C(\log K)^2\frac{N(2A,\delta_0)}{N(A,\delta_0)}\\
 &<CK(\log K)^2.
\end{aligned}
\]
\noindent{[Annotations to the two occurrences of $I^{(\delta_1)}$ above: the first interval intersects $A^{(1)}$''; the second interval intersects $2A^{(1)}$.]}

Apply \paperref{result:2.1}{Lemma~\ref*{result:2.1}} with $(\delta,\eta)=(\delta_1,\delta_2)$ and $\tau=1/10$, we obtain, by the similar argument as before,
\[
 \frac\kappa{10K(\log K)^2}<\kappa_2<\kappa
 \quad\text{\&}\quad1\le R_2\lesssim[K(\log K)^2]^3,
\]
and $A^{(2)}$, obtained by retaining equally many $\delta_1$-intervals in the selected intersections $A^{(1)}\cap I^{(\delta_2)}$, such that
\[
 N(A^{(2)}\cap I^{(\delta_2)},\delta_1)\sim\kappa_2\frac{\delta_2}{\delta_1}
\]
\[
 \text{\&}\quad\frac{N(2(A^{(2)}\cap I^{(\delta_2)}),\delta_1)}{N(A^{(2)}\cap I^{(\delta_2)},\delta_1)}\sim R_2,
\]

if $A^{(2)}\cap I^{(\delta_2)}\ne\varnothing$, and
\[
 N(A^{(2)},\delta_1)\gtrsim(\log K)^{-2}N(A^{(1)},\delta_1).
\]
Furthermore, the construction retains or discards each complete $\delta_1$-interval. Thus, for every surviving $I^{(\delta_1)}$, we have $A^{(2)}\cap I^{(\delta_1)}=A^{(1)}\cap I^{(\delta_1)}$. Consequently, properties \paperref{eq:2.8}{(\ref*{eq:2.8})}, \paperref{eq:2.9}{(\ref*{eq:2.9})} and the exact lower-level branching numbers are preserved.

After $t$ steps, we obtain a set $A^{(t)}$ satisfies
\begin{equation}
 N(A^{(t)}\cap I^{(\delta_s)},\delta_{s-1})\sim\kappa_s\frac{\delta_s}{\delta_{s-1}}\tag{2.11}\label{eq:2.11}
\end{equation}
\begin{equation}
 \text{\&}\quad\frac{N(2(A^{(t)}\cap I^{(\delta_s)}),\delta_{s-1})}{N(A^{(t)}\cap I^{(\delta_s)},\delta_{s-1})}\sim R_s\tag{2.12}\label{eq:2.12}
\end{equation}
if $A^{(t)}\cap I^{(\delta_s)}\ne\varnothing$, since in each step, we preserve
\[
 A^{(s)}\cap I^{(\delta_j)}=A^{(j)}\cap I^{(\delta_j)}\quad\text{for all }0\le j\le s\le t.
\]
Moreover,
\[
 \frac{N(A^{(t)},\delta)}{N(A,\delta)}
 =\prod_{s=1}^t\frac{N(A^{(s)},\delta)}{N(A^{(s-1)},\delta)}.
\]
(Note that by \paperref{eq:2.11}{(\ref*{eq:2.11})},

\[
\begin{aligned}
 N(A^{(s)},\delta)
 &=\sum_{|I|=\delta_1}N(A^{(s)}\cap I^{(\delta_1)},\delta_0)\\
 &\sim\kappa_1\frac{\delta_1}{\delta_0}N(A^{(s)},\delta_1)\\
 &\sim\cdots\sim\left(\prod_{i=1}^{s-1}\kappa_i\frac{\delta_i}{\delta_{i-1}}\right)N(A^{(s)},\delta_{s-1}).
\end{aligned}
\]
)
\begin{equation}
\begin{aligned}
 \frac{N(A^{(t)},\delta)}{N(A,\delta)}
 &=\prod_{s=1}^t\frac{N(A^{(s)},\delta_{s-1})}{N(A^{(s-1)},\delta_{s-1})}\\
 &>\bigl(C(\log K)^{-2}\bigr)^t>\frac1K,
\end{aligned}\tag{2.13}\label{eq:2.13}
\end{equation}
by our choice of $t$.

\begin{center}
\begin{tikzpicture}[x=0.9cm,y=0.8cm, every node/.style={font=\small}]
\coordinate (r) at (4,0);
\coordinate (a) at (2.4,-0.8); \coordinate (b) at (5.6,-0.8);
\coordinate (c) at (2.4,-2.2); \coordinate (d) at (5.6,-2.2);
\coordinate (c1) at (1.8,-3); \coordinate (c2) at (3,-3);
\coordinate (d1) at (5,-3); \coordinate (d2) at (6.2,-3);
\foreach \v in {a,b}\draw[thick] (r)--(\v);
\draw[thick,dotted] (a)--(c);\draw[thick,dotted] (b)--(d);
\foreach \v in {c1,c2}\draw[thick] (c)--(\v);
\foreach \v in {d1,d2}\draw[thick] (d)--(\v);
\foreach \x in {1.8,3,5,6.2}{
 \draw[thick,dotted] (\x,-3)--(\x,-4.5);
 \draw[thick] (\x,-4.5)--(\x-0.3,-5.2) (\x,-4.5)--(\x+0.3,-5.2);
 \fill (\x,-4.5) circle (2pt);
 \fill (\x-0.3,-5.2) circle (2pt);\fill (\x+0.3,-5.2) circle (2pt);
}
\foreach \v in {r,a,b,c,d,c1,c2,d1,d2}\fill (\v) circle (2pt);
\draw[red,thick] (1.95,-1.85) rectangle (6.05,-2.55);
\draw[cyan!70!blue,thick] (1.45,-2.7) rectangle (3.35,-3.35);
\filldraw[fill=cyan,draw=black] (c) circle (3pt);
\node[cyan!70!blue,right] at (2.45,-2.25) {$I^{(\delta_s)}$};
\node[red] at (4,-1.45) {$\#=N(A^{(t)},\delta_s)$};
\foreach \y/\l in {0/{\delta_t},-0.8/{\delta_{t-1}},-2.2/{\delta_s},-3/{\delta_{s-1}},-4.5/{\delta_1},-5.2/{\delta_0=\delta}}
 \node[right] at (7,\y) {$\l$};
\node at (4,-0.8) {$\cdots$};\node at (4,-2.2) {$\cdots$};
\node at (4,-3) {$\cdots$};\node at (4,-4.5) {$\cdots$};\node at (4,-5.2) {$\cdots$};
\node[cyan!70!blue,align=center,anchor=east] at (1,-2) {$\#\text{ of children}=$\\$N(A^{(t)}\cap I^{(\delta_s)},\delta_{s-1})$\\$\displaystyle\sim\kappa_s\frac{\delta_s}{\delta_{s-1}}$};
\node at (4,-6) {tree structure of $A^{(t)}$};
\end{tikzpicture}
\end{center}

Then, we will find some $s$ at which $R_s\lessapprox1$. Denote $A'=A^{(t)}$. Note that for $\delta_0<\delta_1<1$,
\[
 A'=\bigcup A'\cap I^{(\delta_1)}=: \bigcup_j A_j^{(\delta_1)}
 \quad\text{\&}\quad
 \bigcup A_j^{(\delta_1)}+\bigcup A_j^{(\delta_1)}
 =\bigcup_{j,j'}(A_j^{(\delta_1)}+A_{j'}^{(\delta_1)}),
\]
we have
\[
\begin{aligned}
 N(2A',\delta)
 &=\sum_{|I|=\delta_1}N((2A')\cap I,\delta)\\
 &=\sum_{|I|=\delta_1}N\left(\bigcup_{j,j'}((A_j^{(\delta_1)}+A_{j'}^{(\delta_1)})\cap I),\delta\right)
\end{aligned}
\]

\[
 N(2A',\delta)\ge N(2A',\delta_1)\min_{j,j'}N(A_j^{(\delta_1)}+A_{j'}^{(\delta_1)},\delta),
\]
since $A_j^{(\delta_1)}+A_{j'}^{(\delta_1)}$ has size $\sim\delta_1$. Recall $N(A_j^{(\delta_1)},\delta_0)\sim\kappa_1\delta_1/\delta_0$ and \paperref{eq:2.12}{(\ref*{eq:2.12})}, by sumset estimate \paperref{eq:1.2}{(\ref*{eq:1.2})}, for any $j,j'$,
\[
 \frac{N(A_j^{(\delta_1)}+A_{j'}^{(\delta_1)},\delta_0)}{N(A_j^{(\delta_1)},\delta_0)}
 >\left[\frac{N(2A_j^{(\delta_1)},\delta_0)}{N(A_j^{(\delta_1)},\delta_0)}\right]^{1/3}\sim R_1^{1/3}
\]
\[
 \Longrightarrow\quad N(2A',\delta)\gtrsim N(2A',\delta_1)\left(\kappa_1\frac{\delta_1}{\delta_0}\right)R_1^{1/3}.
\]
Iteration implies therefore,
\[
\begin{aligned}
 N(2A',\delta)
 &\ge c^tN(2A',\delta_t)\prod_{s=1}^t\left(\kappa_s\frac{\delta_s}{\delta_{s-1}}\right)R_s^{1/3}\\
 &\ge c^tN(A',\delta_t)\frac{\delta_t}{\delta}\prod_{s=1}^t\kappa_sR_s^{1/3},
\end{aligned}
\]
while
\[
 N(A',\delta)<\left(C\kappa_1\frac{\delta_1}{\delta_0}\right)N(A',\delta_1)
 <\cdots<\left(C^t\frac{\delta_t}{\delta}\prod_{s=1}^t\kappa_s\right)N(A',\delta_t).
\]
Combining these two estimates and \paperref{eq:2.13}{(\ref*{eq:2.13})}, we get
\[
 c^t\prod_{s=1}^tR_s^{1/3}
 <\frac{N(2A',\delta)}{N(A',\delta)}
 <K\frac{N(2A,\delta)}{N(A,\delta)}<K^2.
\]

By pigeonholing principle, there is some $1\le s\le t$ such that
\[
 R_s<CK^{6/t}<(\log K)^C.
\]
Fix $I=I^{(\delta_s)}$ such that $A'\cap I\ne\varnothing$. Then by \paperref{eq:2.11}{(\ref*{eq:2.11})} and \paperref{eq:2.12}{(\ref*{eq:2.12})},
\[
 \delta_{s-1}N(A'\cap I,\delta_{s-1})\sim\kappa_s\delta_s<\kappa\delta_s,
\]
\[
 \text{\&}\quad\frac{N(2(A'\cap I),\delta_{s-1})}{N(A'\cap I,\delta_{s-1})}\sim R_s<(\log K)^C.
\]
Next, we apply \paperref{result:1.5}{Corollary~\ref*{result:1.5}} to the rescaled set $\delta_s^{-1}(A'\cap I)$ with $\delta$ replaced by $\delta_{s-1}/\delta_s$ and $K$ by $(\log K)^C$. Thus,
\[
\begin{aligned}
 \sigma&\sim(\log K)^C\left[\frac{\delta_{s-1}}{\delta_s}
 N\left(\delta_s^{-1}(A'\cap I),\frac{\delta_{s-1}}{\delta_s}\right)\right]^{(\log K)^{-C}}\\
 &<(\log K)^C\kappa^{(\log K)^{-C}},
\end{aligned}
\]
and
\[
 \delta_s^{-1}(A'\cap I)\subset\left\{\xi_0'+\frac a{q'}:a\in\mathbb Z\right\}+I\left(0,\frac\sigma{q'}\right)
\]
for some $q'\in\mathbb Z_+$ with $\sigma/q'>\delta_{s-1}/\delta_s$. Let $\xi_0=\delta_s\xi_0'$ and $q=q'/\delta_s$, we have
\[
 A'\cap I\subset\left\{\xi_0+\frac aq:a\in\mathbb Z\right\}+I\left(0,\frac\sigma q\right)
\]

\[
 \text{\&}\quad\delta_{s-1}<\frac\sigma q<\frac1q<\delta_s.
\]
Here we in addition assume $\sigma<1$, which is acceptable whenever
\[
 C(\log\log K)(\log K)^C<\log\kappa^{-1}
 \quad\Longleftarrow\quad C(\log K)^C<\log\kappa^{-1}.
\]
For example, let $\kappa=Q^{-C}$ and $K=\log Q$ for some sufficiently large $Q$.

Furthermore, one may notice that $q$, $1/\delta_s<q<1/\delta_{s-1}$, depends on the interval $I=I^{(\delta_s)}$. But the value of $q$ (in dyadic number) only takes $\sim\log(\delta_s/\delta_{s-1})\le\log(1/\delta)$. Hence by dyadic pigeonholing, we may replace $A'$ by a further subset $A''$ such that
\[
 N(A'',\delta)>(\log\delta^{-1})^{-1}N(A',\delta)>
 \frac1{K\log\delta^{-1}}N(A,\delta)
\]
\[
 \text{\&}\quad A''\cap I\subset\left\{\xi_0+\frac aq:a\in\mathbb Z\right\}+I\left(0,\frac\sigma q\right)
\]
for fixed $q$, $\delta_s>1/q>\sigma/q>\delta_{s-1}$.

Hence, we prove the following fact which will substitute \paperref{result:1.5}{Corollary~\ref*{result:1.5}}.

\begin{lemma}\label{result:2.2} Let $A\subset\mathbb R$ be a bounded set, $\delta>0$, and assume
\[
 \frac{N(2A,\delta)}{N(A,\delta)}<K
 \quad\text{\&}\quad\delta N(A,\delta)<\kappa
\]
for some suitable $\kappa$ and $K$. Then there are $\sigma>0$, $\sigma\delta'>\delta$ and a subset $A'$ of $A$ satisfying the conditions:
\begin{itemize}
\item $\sigma<(\log K)^C\kappa^{(\log K)^{-C}}$.
\item $A'$ is contained in a union of intervals of size $\sigma\delta'$ and separation at least $\delta'$. (So the ``ratio'' $=(\sigma\delta')/\delta'=\sigma$ is smaller.)
\item $N(A',\delta)>\left(K\log\frac1\delta\right)^{-1}N(A,\delta)$.
\end{itemize}
\end{lemma}
Now, \paperref{result:2.2}{Lemma~\ref*{result:2.2}} provides the ``porosity property'' at a single scale $\sigma\delta'>\delta$. However it is still too weak for application to discretized ring problem. Next, by multiscale analysis, we will construct a subtree with ``multi-scale porosity property''.

From now on, we assume $A$ is a $(\delta,\sigma)_1$-set, $\delta^{1-\sigma+\varepsilon}<|A|<\delta^{1-\sigma-\varepsilon}$ and

\[
 |A+A|<\delta^{1-\sigma-\varepsilon}\quad(<\delta^{-\varepsilon}|A|).
\]

\subsection{Find good scales at which the double constant is small enough}\label{sec:2.2}
\textcolor{red}{
In this section, we first fix a sequence of scales. To apply \paperref{result:2.2}{Lemma~\ref*{result:2.2}}, we need the doubling constant for logarithmic growth. By once again using the regular tree structure and the submultiplicativity of the double constant, we can ensure that, provided the tree has a sufficient number of levels, the double constants for the majority of these levels remain sufficiently small. Finally, the non-concentration condition ensures that the smallest "good" scale must appear at a scale slightly larger than $\delta$.
}

More precisely, fix a large integer $Q$, which depends on $\varepsilon$ and we assume
\[
 \log Q<\frac1\varepsilon.
\]
Let
\[
 t\sim\frac{\log(1/\delta)}{\log Q}\quad\Longleftrightarrow\quad Q^t\sim\frac1\delta.
\]
The doubling constant is $Q^{t\varepsilon}<Q^{t/\log Q}$, so the first tree construction is aimed to bringing down the doubling constant to $(\log Q)^C$. Let
\[
 \delta_s=Q^s\delta,\quad0\le s\le t,\qquad
 \delta=\delta_0<\delta_1<\cdots<\delta_{t-1}<\delta_t\sim1.
\]
Similar as before, we will exploit multiplicative properties of the doubling constant at various scales $\{\delta_s\}_{s=0}^t$.

Start at $\delta_1$, if $A\cap I^{(\delta_1)}\ne\varnothing$, then $1\le N(A\cap I^{(\delta_1)},\delta)\le\delta_1/\delta=Q$, trivially. Furthermore, we have trivially
\[
 1\le\frac{N(2(A\cap I^{(\delta_1)}),\delta)}{N(A\cap I^{(\delta_1)},\delta)}
 \le\frac{N(A\cap I^{(\delta_1)},\delta)^2}{N(A\cap I^{(\delta_1)},\delta)}\le Q.
\]
Thus, by the exact regularization described above, there is $A^{(1)}\subset A$, obtained by retaining equally many $\delta$-intervals inside each selected $I^{(\delta_1)}$, such that

\[
 N(A^{(1)}\cap I^{(\delta_1)},\delta)=N_1
 \quad\text{\&}\quad\frac{N(2(A^{(1)}\cap I^{(\delta_1)}),\delta)}{N(A^{(1)}\cap I^{(\delta_1)},\delta)}\sim K_1,
\]
if $A^{(1)}\cap I^{(\delta_1)}\ne\varnothing$. Moreover, $|A^{(1)}|>(\log Q)^{-2}|A|$. Repeat this process then with $(A,\delta_0,\delta_1)$ replaced by $(A^{(1)},\delta_1,\delta_2)$ etc.

After $t$-steps, we obtain $A^{(t)}\subset A$ satisfying the conditions:
\begin{equation}
 N(A^{(t)}\cap I^{(\delta_s)},\delta_{s-1})=N_s\tag{2.14}\label{eq:2.14}
\end{equation}
\noindent{[Annotation: each vertex at level $s$ has $N_s$ children.]}
\begin{equation}
 \text{\&}\quad\frac{N(2(A^{(t)}\cap I^{(\delta_s)}),\delta_{s-1})}{N(A^{(t)}\cap I^{(\delta_s)},\delta_{s-1})}\sim K_s\tag{2.15}\label{eq:2.15}
\end{equation}
if $A^{(t)}\cap I^{(\delta_s)}\ne\varnothing$, $1\le s\le t$. And
\begin{equation}
 |A^{(t)}|>(\log Q)^{-2t}|A|.\tag{2.16}\label{eq:2.16}
\end{equation}
Denote $A'=A^{(t)}$. Observe that for any $1\le r\le s\le t$,
\[
 N(A'\cap I^{(\delta_s)},\delta_{r-1})=N_sN_{s-1}\cdots N_r\equiv\text{constant}.
\]
In particular, $|A'\cap I^{(\delta_s)}|=\delta N(A'\cap I^{(\delta_s)},\delta)$ are the same for all intervals $I^{(\delta_s)}$. So considering the sumset, we have again by \paperref{result:1.1}{Lemma~\ref*{result:1.1}}:
\[
\begin{aligned}
 N(2A',\delta)&=\sum_{|I|=\delta_1}N((2A')\cap I,\delta)\\
 &\ge N(2A',\delta_1)\min_{j,j'}N(A_j^{(\delta_1)}+A_{j'}^{(\delta_1)},\delta).
\end{aligned}
\]

\[
\begin{aligned}
 N(2A',\delta)
 &\ge N(2A',\delta_1)\cdot
 \left[\frac{N(2A_j^{(\delta_1)},\delta)}{N(A_j^{(\delta_1)},\delta)}\right]^{1/3}\cdot N(A_j^{(\delta_1)},\delta)\\
 &\gtrsim N(2A',\delta_1)N_1K_1^{1/3}.
\end{aligned}
\]
Iterating implies by \paperref{eq:2.16}{(\ref*{eq:2.16})}
\[
\begin{aligned}
 N(2A',\delta)
 &\ge c^tN(2A',\delta_t\sim1)\prod_{s=1}^tN_sK_s^{1/3}\\
 &=c^tN(A',\delta)\prod_{s=1}^tK_s^{1/3}\\
 &>c^t\frac{(\log Q)^{-2t}}\delta\cdot|A|\cdot\prod_{s=1}^tK_s^{1/3}\\
 &=N(A,\delta)\prod_{s=1}^t\frac{K_s^{1/3}}{C(\log Q)^2}.
\end{aligned}
\]
Therefore,
\begin{equation}
 K:=\frac{N(2A,\delta)}{N(A,\delta)}>\prod_{s=1}^t\frac{K_s^{1/3}}{C(\log Q)^2}.\tag{2.17}\label{eq:2.17}
\end{equation}
Recall we have assumed $A$ has small sumset:
\begin{equation}
 K\le\delta^{-\varepsilon}\tag{2.18}\label{eq:2.18}
\end{equation}
\[
 \text{\&}\quad\delta^{-\sigma+\varepsilon}<N(A,\delta)<\delta^{-\sigma-\varepsilon}.
\]
By carefully pigeonholing argument, we may choose a set of "good" index $S\subset$
$\{1,\ldots,t\}$ such that
\begin{equation}
 t_1=|S|>\frac{1-\sigma}4t,\tag{2.19}\label{eq:2.19}
\end{equation}
and for all $s\in S$,
\begin{equation}
 N_s<Q^{\sigma/10+9/10+\varepsilon}\tag{2.20}\label{eq:2.20}
\end{equation}
\noindent{[Annotation to $\sigma/10+9/10$: slightly large than $\sigma$ but still small than $1$.]}
\begin{equation}
 \text{\&}\quad K_s<K^{C/t}(\log Q)^C\quad\text{with }C=C(\sigma).\tag{2.21}\label{eq:2.21}
\end{equation}
(Indeed, let $S'=\{1\le s\le t:N_s<Q^{\sigma/10+9/10+\varepsilon}\}$. Then we have
\[
 Q^{t(\sigma+\varepsilon)}\sim N(A,\delta)>\prod_{s=1}^tN_s>
 \left(Q^{\sigma/10+9/10+\varepsilon}\right)^{t-|S'|}
\]
\[
 \Longrightarrow\quad t(\sigma+\varepsilon)\ge
 \left(\frac\sigma{10}+\frac9{10}+\varepsilon\right)(t-|S'|)
\]
\[
 \Longrightarrow\quad|S'|\ge\frac9{20}(1-\sigma)t\ge\frac{1-\sigma}3t.
\]
Next, let $S=\{s\in S':K_s<K^{C/t}(\log Q)^C\}$, where $C=C(\sigma)$ is chosen sufficiently large. Write $C_0$ for the absolute constant in \paperref{eq:2.17}{(\ref*{eq:2.17})}. Taking logarithms of the full product gives
\[
 \sum_{s=1}^t\log K_s
 \le 3\log K+6t\log\log Q+3t\log C_0.
\]
For $J=\{1\le s\le t:K_s\ge K^{C/t}(\log Q)^C\}$, it follows that
\[
 |J|\,C\left(\frac{\log K}{t}+\log\log Q\right)
 \le 3\log K+6t\log\log Q+3t\log C_0.
\]
Hence $|J|\le7t/C$ for $Q$ sufficiently large, and
\[
 |S|\ge |S'|-\frac{7t}{C}
 \ge\left(\frac9{20}(1-\sigma)-\frac7C\right)t
 >\frac{1-\sigma}{3}t>\frac{1-\sigma}{4}t
\]
when, for example, $C\ge84/(1-\sigma)$.)

Recall we assume $\varepsilon<1/\log Q$ and $t\sim\log\delta^{-1}/\log Q$, it follows that

\begin{equation}
 K_s<\delta^{-C/(t\log Q)}(\log Q)^C<(\log Q)^C\tag{2.22}\label{eq:2.22}
\end{equation}
for all $s\in S$.

For later use, we make the following remark.

Let $1\le s\le t$, $s\le\min S$, consider $A'\cap I^{(\delta_s)}\ne\varnothing$, we have
\[
 N(A'\cap I^{(\delta_s)},\delta)=\prod_{1\le s'\le s}N_{s'}.
\]
By the same argument as before, $K_s<(\log Q)^C$ for all but at $7t/C$ values of $s$. Hence, by our construction of $S$, restricting $s'\notin S$, we have $N_{s'}>Q^{\sigma/10+9/10+\varepsilon}$ except for at most $7t/C$ values. Therefore,
\[
\begin{aligned}
 N(A'\cap I^{(\delta_s)},\delta)
 &>Q^{(\sigma/10+9/10+\varepsilon)(s-7t/C)}\\
 &\sim\left(\frac{\delta_s}\delta\right)^{\sigma/10+9/10+\varepsilon}\delta^{(\sigma/10+9/10+\varepsilon)7/C}\\
 &>\left(\frac{\delta_s}\delta\right)^{\sigma/10+9/10+\varepsilon}\delta^{7/C}.
\end{aligned}
\]
Since $A$ is $(\delta,\sigma)_1$-set, we have
\[
 N(A\cap I^{(r)},\delta)\le\delta^{-\varepsilon}\left(\frac r\delta\right)^\sigma,
\]
for any $r$-interval $I^{(r)}$. Choosing $r=\delta_s$, we obtain

\[
 \left(\frac{\delta_s}\delta\right)^\sigma>
 \left(\frac{\delta_s}\delta\right)^{\sigma/10+9/10+\varepsilon}\delta^{7/C+\varepsilon}.
\]
Hence,
\begin{equation}
 \delta_s<\delta^{1-\varepsilon'}\quad\text{for some }\varepsilon'=\varepsilon'(C)=\varepsilon'(\sigma)\text{ small enough}.\tag{2.23}\label{eq:2.23}
\end{equation}

\subsection{Iteration to obtain the multiscale porous tree}\label{sec:2.3}
\textcolor{red}{
Now, we start constructing the desired tree. In the previous section, we already found a good set of scales. Between each pair of adjacent scales, we apply \paperref{result:2.2}{Lemma~\ref*{result:2.2}} to obtain an intermediate scale, which exhibits porosity—provided that the double constant is sufficiently small. This process ultimately yields a multiscale porous tree.}

Let $S=\{s_1,s_2,\ldots,s_{t_1}\}$ with $s_1>s_2>\cdots>s_{t_1}$. Take $s=s_1$ and fix $I^{(\delta_s)}$ such that $A'\cap I^{(\delta_s)}\ne\varnothing$. Then by \paperref{eq:2.20}{(\ref*{eq:2.20})} and \paperref{eq:2.22}{(\ref*{eq:2.22})}, we have
\[
 N(A'\cap I^{(\delta_s)},\delta_{s-1})=N_s<Q^{\sigma/10+9/10+\varepsilon}
\]
\[
 \text{\&}\quad\frac{N(2(A'\cap I^{(\delta_s)}),\delta_{s-1})}{N(A'\cap I^{(\delta_s)},\delta_{s-1})}
 \sim K_s<(\log Q)^C.
\]
Apply \paperref{result:2.2}{Lemma~\ref*{result:2.2}} to the rescaled set $\delta_s^{-1}(A'\cap I^{(\delta_s)})$ with $\delta$ replaced by $\delta_{s-1}/\delta_s=1/Q$, $\kappa=Q^{-(1-\sigma)/10+\varepsilon}$ and $K=(\log Q)^C$. We obtain subsets $\{B_k\}$ of $A'\cap I^{(\delta_{s_1})}$ and $\delta_{s_1-1}<Q^{-1}\lambda_1$ \& $\lambda_1\le\delta_{s_1}$ satisfying

(In fact, we can only obtain
\[
 \delta_{s_1-1}<Q^{-(\log\log Q)^{-C}}\lambda_1\quad(>Q^{-1}\lambda_1).
\]
But in future applications, the role of these two factors will be the same. So, WLOG, we assume $\delta_{s_1-1}<Q^{-1}\lambda_1$.)

\begin{enumerate}[label=(\roman*)]
\item\label{item:2-i} Each $B_k$ is a union of sets $A'\cap I^{(\delta_{s_1-1})}$.
\end{enumerate}

\begin{enumerate}[label=(\roman*),start=2]
\item\label{item:2-ii} $B_k$ is contained in an interval of size $Q^{-(\log\log Q)^{-C}}\lambda_1$.
\item\label{item:2-iii} $\operatorname{dist}(B_k,B_{k'})>\lambda_1$ for $k\ne k'$.
\item\label{item:2-iv} $\displaystyle\sum_kN(B_k,\delta_{s_1-1})>\frac1{(\log Q)^C}N(A'\cap I^{(\delta_{s_1})},\delta_{s_1-1})$.
\end{enumerate}
One may note that the value of $\lambda_1$ depends on $I^{(\delta_{s_1})}$, but, after a dyadic pigeonholing, since its range is $\le\log(\delta_{s_1}/\delta_{s_1-1})\le\log Q$, it may be specified allowing a loss of an extra $(\log Q)^{-1}$-factor.

Thus we obtain subsets $B_k\subset A'$ satisfying \paperref{item:2-i}{(i)}--\paperref{item:2-iii}{(iii)} and
\[
 \sum_kN(B_k,\delta_{s_1-1})>(\log Q)^{-C}N(A',\delta_{s_1-1})
\]
\[
 \Longrightarrow\quad\sum_k|B_k|>(\log Q)^{-C}|A'|.
\]
Next, pass to $s_2<s_1$ and consider $A'\cap I^{(\delta_{s_2})}\cap(\bigcup B_k)$ for which we repeat the procedure. Thus we obtain subsets $B_{k_1,k_2}\subset B_{k_1}$ for each $B_{k_1}$, and $\delta_{s_2-1}<\lambda_2<\delta_{s_2}$ such that \paperref{item:2-i}{(i)} $\cdots$ \paperref{item:2-ii}{(ii)} $\cdots$ \paperref{item:2-iii}{(iii)} $\cdots$ and \paperref{item:2-iv}{(iv)}:
\[
 \sum_{k_1,k_2}|B_{k_1,k_2}|>(\log Q)^{-C}\sum_{k_1}|B_{k_1}|>(\log Q)^{-2C}|A'|.
\]
Iteration over $s_1,\ldots,s_{t_1}$ gives a tree of subsets $B_{k_1,\ldots,k_{t_1}}\subset A$ satisfying
\begin{itemize}
\item Each $B_{k_1,\ldots,k_u}$ is a union of sets $A'\cap I^{(\delta_{s_u-1})}$.
\end{itemize}

\begin{itemize}
\item Each $B_{k_1,\ldots,k_u}$ is contained in an interval of size
\begin{equation}
 Q^{-(\log\log Q)^{-C}}\lambda_u.\tag{2.24}\label{eq:2.24}
\end{equation}
\item
\begin{equation}
 \operatorname{dist}(B_{k_1,\ldots,k_u},B_{\ell_1,\ldots,\ell_u})>\lambda_u
 \quad\text{for }\overline{k}\ne\overline{\ell}.\tag{2.25}\label{eq:2.25}
\end{equation}
\end{itemize}
(Here $1\le u\le t_1$ and we identify $B_{k_1,\ldots,k_u}$, the subsets obtained in $u$-steps, with
\[
 \bigcup_{k_{u+1},\ldots,k_{t_1}}B_{k_1,\ldots,k_u,k_{u+1},\ldots,k_{t_1}}.
\]
)
\begin{itemize}
\item
\begin{equation}
 \left|\bigcup_{\overline{k}}B_{k_1,\ldots,k_{t_1}}\right|
 >(\log Q)^{-Ct_1}|A'|>(\log Q)^{-Ct}|A|
 =\delta^{C\frac{\log\log Q}{\log Q}}|A|.\tag{2.26}\label{eq:2.26}
\end{equation}
\item
\begin{equation}
 \lambda_{u+1}<\delta_{s_{u+1}}=Q^{s_{u+1}}\delta
 \le Q^{s_u-1}\delta=\delta_{s_u-1}<\frac1Q\lambda_u.\tag{2.27}\label{eq:2.27}
\end{equation}
\end{itemize}
The tree $(B_{k_1,\ldots,k_{t_1}})$ need not yet be regular/uniform, in other word, at each given level of the tree, the number of branchings is the same for each node. But it may be regularized by dyadic pigeonholing from bottom to top with loss of an additional factor:
\[
 \frac{1}{\log(Q^{s_{t_1}})}
 \prod_{u=1}^{t_1-1}\frac1{\log(Q^{s_u-s_{u+1}+1})}
 \ge \frac{1}{(\log Q)^t}
 =\delta^{\log\log Q/\log Q}.
\]
This loss is absorbed in \paperref{eq:2.26}{(\ref*{eq:2.26})}. We continue to write $B_{k_1,\ldots,k_u}$ for the retained sets, each redefined as the union of its surviving leaves.

Let $K_u$ denote the branching number from level $u-1$ to level $u$, that is,
\[
 K_u=\#\{B_{k_1,\ldots,k_{u-1},k_u}:B_{k_1,\ldots,k_{u-1},k_u}\subset B_{k_1,\ldots,k_{u-1}}\}.
\]
Hence,
\begin{equation}
 B_{k_1,\ldots,k_{u-1}}=\bigcup_{k_u\le K_u}B_{k_1,\ldots,k_{u-1},k_u}.\tag{2.28}\label{eq:2.28}
\end{equation}
And in particular,
\begin{equation}
 \delta^{C\frac{\log\log Q}{\log Q}}|A|
 <|B_{\varnothing}|=\sum_{\overline{k}}|B_{k_1,\ldots,k_{t_1}}|
 =\left(\prod_{s\le t_1}K_s\right)\cdot\rho\tag{2.29}\label{eq:2.29}
\end{equation}
with $\rho=|B_{k_1,\ldots,k_{t_1}}|$ independent of the leaf.

\annot{[Annotations: $\prod_{s\le t_1}K_s$ is ``\# of all branchings''; $\rho$ is ``volume of each leaf''.]}

Recall \paperref{eq:2.23}{(\ref*{eq:2.23})}, for $1\le s\le\min S$, we have $\delta_s<\delta^{1-\varepsilon'}$, so $\delta_{s_{t_1}}<\delta^{1-\varepsilon'}$. Note that $\rho<\lambda_{t_1}<\delta_{s_{t_1}}$, we get
\begin{equation}
\begin{aligned}
 \prod_{s\le t_1}K_s
 &>\delta^{\varepsilon'-1}\cdot\delta^{C\frac{\log\log Q}{\log Q}}|A|\\
 &>\delta^{2\varepsilon'}\delta^{-1}\delta^{1-\sigma+\varepsilon}\\
 &>\delta^{-\sigma+2\varepsilon'}.
\end{aligned}\tag{2.30}\label{eq:2.30}
\end{equation}

\subsection{Merge levels with few branches to obtain lots of effective levels}\label{sec:2.4}
\textcolor{red}{
Although by \paperref{eq:2.30}{(\ref*{eq:2.30})}, we know for ``some'' $s$, $K_s$ is ``large''. But it is still possible to happen that the branchings which are too small are not negligible. These levels, which have a small number of branches, mean that we cannot achieve a positive gain at this level when applying the geometric argument presented in Section 4. To take care of this, we make a further construction, identifying consecutive levels to a single one. The non-concentration condition ensures that this process can run almost entirely from the bottom to the top.}

Define
\[
M=e^{\sqrt{\log Q}},\qquad \eta_Q=Q^{-(\log\log Q)^{-C}}.
\]
\annot{\textbf{Clarification.} In the following grouping argument, $r_i$ denotes the number of branching factors in a block. A block ending at level $m$ consists of the factors $K_{m-r_i+1},\ldots,K_m$, and the next block ends at level $m-r_i$. Thus no branching factor is used twice.}

Starting from $m=t_1$, let $r_1\ge1$ be the smallest integer such that
\[
 \prod_{j=t_1-r_1+1}^{t_1}K_j>M^{r_1}.
\]
This is possible by taking $r_1=t_1$ and using \paperref{eq:2.30}{(\ref*{eq:2.30})}, since
\[
 M^{t_1}\le M^t=\delta^{-1/\sqrt{\log Q}}
 <\delta^{-\sigma+2\varepsilon'}<\prod_{s\le t_1}K_s.
\]
If $r_1\le1000$, merge these $r_1$ branchings into a single branching $L>M^{r_1}$. By \paperref{eq:2.24}{(\ref*{eq:2.24})} and \paperref{eq:2.25}{(\ref*{eq:2.25})}, the diameter-to-separation ratio of the new children satisfies
\[
 \sigma_1<\eta_Q.
\]
\textcolor{cyan}{$\sigma_1$ is the ratio at the new lowest branching of the grouped tree.}
For $Q$ sufficiently large,
\[
 \sigma_1^{-1}>M^{\frac{\sqrt{\log Q}}{(\log\log Q)^C}}>M^{100r_1}.
\]

In this case the measure of the tree is unchanged.
\begin{center}
\includegraphics[width=.9\linewidth,height=.65\textheight,keepaspectratio]{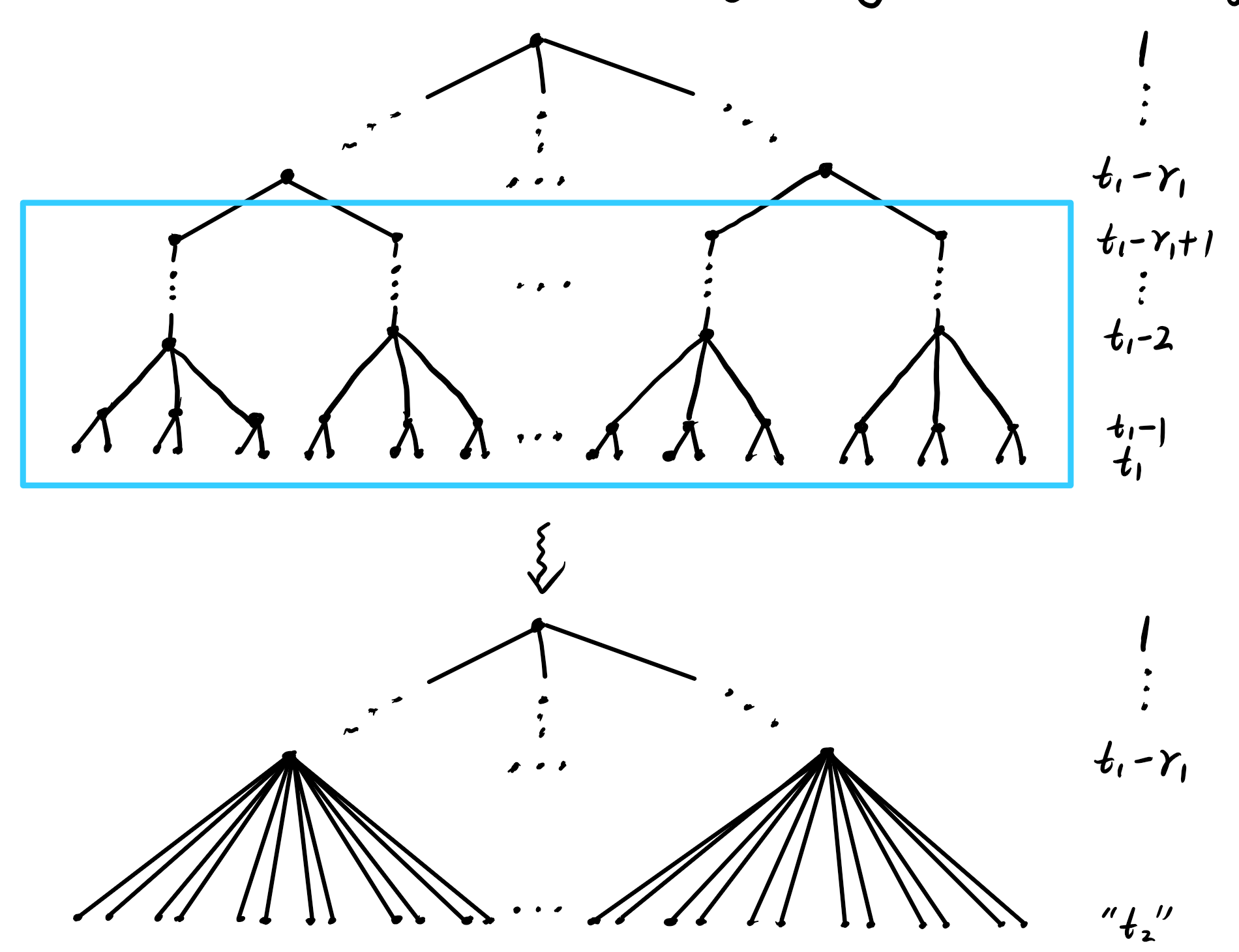}
\end{center}

If $r_1>1000$, take $h_1=\lfloor r_1/100\rfloor$. Retain the upper $r_1-h_1$ branchings of the block. By the minimality of $r_1$,
\[
 \prod_{j=t_1-h_1+1}^{t_1}K_j\le M^{h_1},
\]
and therefore the retained branching satisfies
\[
 L=\prod_{j=t_1-r_1+1}^{t_1-h_1}K_j
 >M^{r_1-h_1}\ge M^{99r_1/100}.
\]
\annot{\textbf{Clarification.} We retain the original union shorthand
$\bigcup_{t_1-\frac{r_1}{100}+1,\ldots,t_1}I^{(\delta_{t_1})}$
for piecing together the small intervals over the indicated levels. These are the label for original levels; the precise block endpoints in the proof are specified by $h_1=\lfloor r_1/100\rfloor$ and the block-length convention above.}

For each vertex at level $t_1-h_1$, retain one descendant at level $t_1$. This keeps all $L$ branchings above level $t_1-h_1$, while replacing each of their children by one small leaf. The selected sets have mutual separation at least $\lambda_{t_1-h_1}$ and diameter at most $\eta_Q\lambda_{t_1}$.

The actual weaker form of \paperref{eq:2.27}{(\ref*{eq:2.27})} is $\lambda_{u+1}\le\eta_Q\lambda_u$. Consequently,
\[
 \sigma_1\le\eta_Q\frac{\lambda_{t_1}}{\lambda_{t_1-h_1}}
 \le\eta_Q^{h_1+1}\le\eta_Q^{r_1/100}
 =e^{-\frac{\log Q}{(\log\log Q)^C}\frac{r_1}{100}}.
\]

Equivalently,
\[
 \sigma_1\le M^{-\frac{\sqrt{\log Q}}{(\log\log Q)^C}\frac{r_1}{100}},
 \qquad \sigma_1^{-1}>M^{100r_1}
\]
for $Q$ sufficiently large. This verifies that the simplified $Q^{-1}$ factor used in part (3) can throughout be replaced by the actual factor $\eta_Q$.

Since all subtrees at the same level have equal mass, the retained measure is at least
\[
 \frac{|B_\varnothing|}{\prod_{j=t_1-h_1+1}^{t_1}K_j}
 \ge M^{-h_1}|B_\varnothing|
 \ge M^{-r_1/100}|B_\varnothing|
 >M^{-r_1/90}|B_\varnothing|.
\]
\begin{center}
\includegraphics[width=.85\linewidth,height=.38\textheight,keepaspectratio]{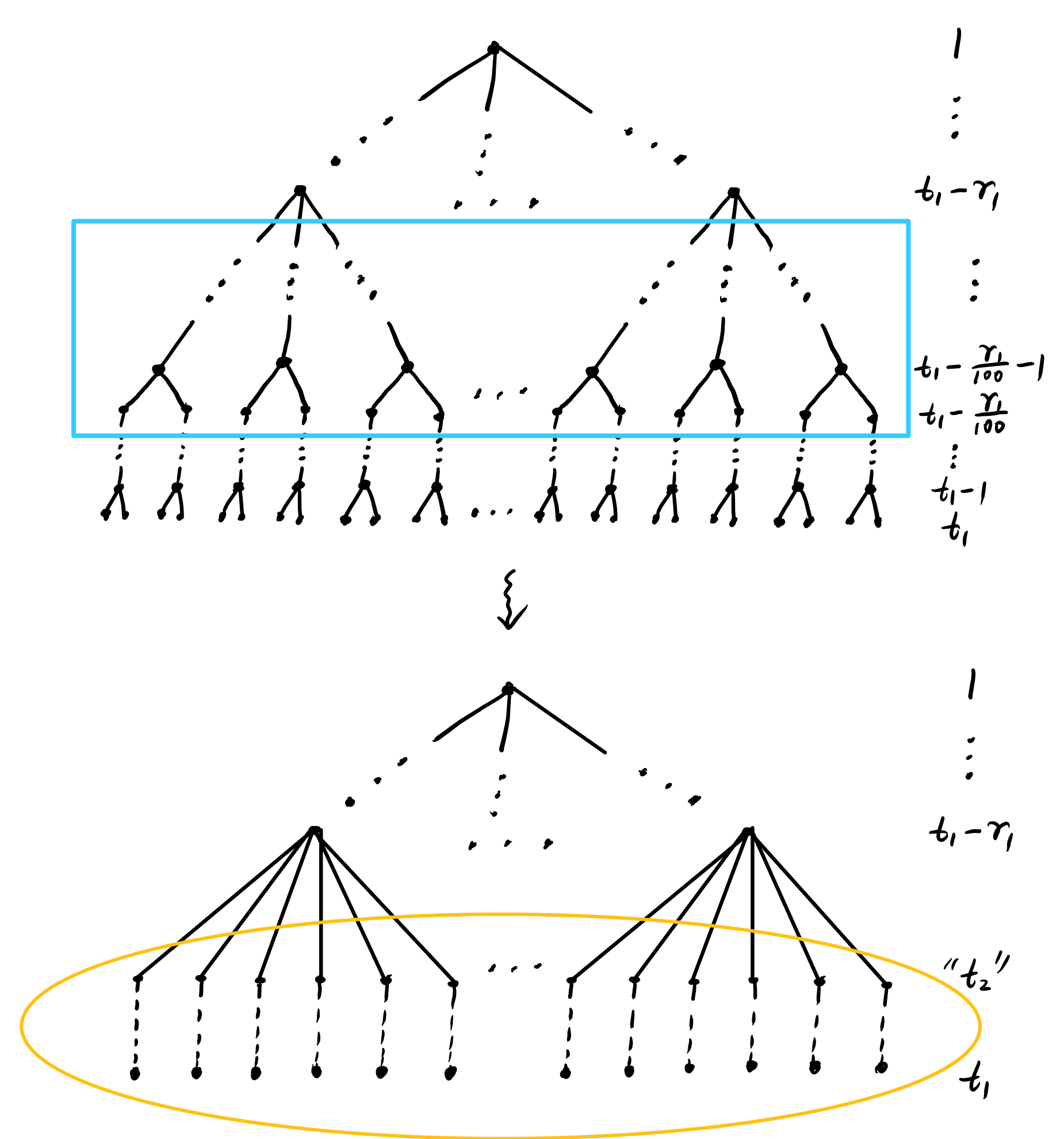}
\end{center}

\textcolor{orange}{For each vertex at level $t_1-h_1$, retain one descendant at level $t_1$. The branchings above level $t_1-h_1$ within the block are all retained.}

Next repeat the procedure with bottom level $m=t_1-r_1$, and then $m=t_1-r_1-r_2$, and so on. In a later long block, each selected level-$m$ descendant is retained together with its entire already constructed subtree below it. The same branching and mass estimates apply because those subtrees have equal mass. The blocks are disjoint and their lengths sum to $t_1-t'$ if the procedure stops at level $t'$.

If $t'=0$, all levels have been processed and $t_1>(1-\sigma)t/3$ already gives \paperref{eq:2.31}{(\ref*{eq:2.31})}. Otherwise, failure to find a further block implies

\[
K_{t'}\cdots K_1\leq M^{t'}<e^{t_1\sqrt{\log Q}}
 <\left(\frac1\delta\right)^{1/\sqrt{\log Q}}.
\]
Hence,
\[
|B_{k_1,\ldots,k_{t'}}|
 =\frac{|B_\varnothing|}{K_1\cdots K_{t'}}
 >\delta^{\varepsilon'}|B_\varnothing|
 >\delta^{\varepsilon'}|A|.
\]
Assume $B_{k_1,\ldots,k_{t'}}$ appears at stage $s_0\in S$, thus
$B_{k_1,\ldots,k_{t'}}\subset B_\varnothing\cap I^{(\delta_{s_0})}$
for some interval $I^{(\delta_{s_0})}$. Then by construction,
$\#\{s\in S:s<s_0\}=t_1-t'$.
\[
s_1>s_2>\cdots>s_0>
\underbrace{s_{t'+1}>\cdots>s_{t_1}}_{{\sum r_i=t_1-t'}}.
\]
{The final block consists of the $r_1$ indices $s_{t_1-r_1+1},\ldots,s_{t_1}$; the preceding block has length $r_2$.}

Next we show this procedure can run long enough. Arguing as in the end of (2),
\begin{center}
\begin{tikzpicture}[x=1cm,y=1cm]
\draw[->] (0,0)--(11.6,0);
\foreach \x/\lab in {0.5/\delta,1.2/\delta_1,9.6/\delta_t,10.6/1}{\draw (\x,-.08)--(\x,.08);\node[below] at (\x,-.08) {$\lab$};}
\foreach \x/\lab in {2.5/\delta_{s_{t_1}},3.6/\delta_{s_{t_1-1}},5.4/\delta_{s_0},7.2/\delta_{s_2},8.4/\delta_{s_1}}{\draw[cyan] (\x,-.08)--(\x,.08);\node[below,text=cyan] at (\x,-.08) {$\lab$};}
\node at (1.9,.15) {$\cdots$};\node[cyan] at (4.5,.15) {$\cdots$};\node[cyan] at (6.3,.15) {$\cdots$};\node at (9,.15) {$\cdots$};
\draw[cyan] (2.5,.35)--(2.5,.5)--(8.4,.5)--(8.4,.35);\node[above,text=cyan] at (5.45,.5) {$S$};
\end{tikzpicture}
\end{center}
\[
\begin{aligned}
N(A'\cap I^{(\delta_{s_0})},\delta)
 &=\prod_{s=1}^{s_0}N_s\\
 &\geq Q^{(\frac{\sigma}{10}+\frac9{10}+\varepsilon)
 (s_0-\#\{s\in S:s<s_0\}-\frac{7t}{C})}
\end{aligned}
\]
\noindent{Annotations: $7t/C$ counts exceptional high-doubling indices. The other term counts the indices in $[1,s_0]$ outside $S$, up to one endpoint.}

\[
N(A'\cap I^{(\delta_{s_0})},\delta)
 >\left(\frac{\delta_{s_0}}\delta\right)^{\frac{\sigma}{10}+\frac9{10}}
 Q^{-(t_1-t')}\delta^{\varepsilon'}.
\]
Since $A'$ is $(\delta,\sigma)$-set, we have
\[
N(A'\cap I^{(\delta_{s_0})},\delta)
 <\delta^{-\varepsilon}\left(\frac{\delta_{s_0}}\delta\right)^\sigma
\]
\[
\Longrightarrow\quad
\delta_{s_0}<\delta^{1-\varepsilon'}Q^{\frac{10}{9(1-\sigma)}(t_1-t')}
\]
\[
\begin{aligned}
\Longrightarrow\quad \delta^{-\sigma+\varepsilon'}
 &<N(A,\delta)<\delta^{-\varepsilon'}N(B_{k_1,\ldots,k_{t'}},\delta)\\
 &<\delta^{-\varepsilon'}\frac{\delta_{s_0}}\delta\\
 &<\delta^{-\varepsilon'}Q^{\frac{10}{9(1-\sigma)}(t_1-t')}
\end{aligned}
\]
\[
\Longrightarrow\quad Q^{t(\sigma-\varepsilon')}
 <Q^{\frac{10}{9(1-\sigma)}(t_1-t')}.
\]
Hence,
\begin{equation}
r_1+r_2+\cdots=t_1-t'>\frac13\sigma(1-\sigma)t,
\tag{2.31}\label{eq:2.31}
\end{equation}
recalling we chose $\varepsilon'$ small enough in terms of $\sigma$.

Thus we obtain a sub-tree $(C_{l_1,\ldots,l_{t_2}})$ with the following properties.

\begin{equation}
\operatorname{dist}(C_{l_1,\ldots,l_s},C_{l'_1,\ldots,l'_s})>\mu_s
\quad\text{for }\bar l\ne\bar l'.
\tag{2.32}\label{eq:2.32}
\end{equation}
\begin{equation}
\operatorname{diam}C_{l_1,\ldots,l_s}<\sigma_s\mu_s.
\tag{2.33}\label{eq:2.33}
\end{equation}
\begin{equation}
\sigma_s<Q^{-(\log\log Q)^{-C}}.
\tag{2.34}\label{eq:2.34}
\end{equation}
\begin{equation}
C_{l_1,\ldots,l_{s-1}}=\bigcup_{l\leq L_s}C_{l_1,\ldots,l_{s-1},l}.
\tag{2.35}\label{eq:2.35}
\end{equation}
\begin{equation}
\prod_{s=1}^{t_2}(L_s\wedge\sigma_s^{-1/100})
 >M^{\frac{99}{100}(r_1+r_2+\cdots)}=M^{\gamma t},
\tag{2.36}\label{eq:2.36}
\end{equation}
where $\frac{99}{100}(r_1+r_2+\cdots)=:\gamma t$ and $\gamma>c\sigma(1-\sigma)$ for an absolute constant $c>0$. The preliminary regularization loss is $M^{-O(t\log\log Q/\sqrt{\log Q})}$ and is absorbed in the slack of the following estimate by choosing $Q$ sufficiently large:
\begin{equation}
|C_\varnothing|=\left|\bigcup_{\bar l}C_{l_1,\ldots,l_{t_2}}\right|
 >M^{-\frac1{90}(r_1+r_2+\cdots)}|B_\varnothing|
 >M^{-\frac1{88}\gamma t}|A|.
\tag{2.37}\label{eq:2.37}
\end{equation}

\section{An alternative multiscale porous tree construction}\label{sec:3}

In this section, we give a different way to handle the additive part. We will use the \paperref{result:1.3}{Pl\"unnecke--Ruzsa inequality} to replace the \paperref{result:1.4}{Freiman's theorem} in \paperref{sec:2}{Section~2}. And finally we still obtain a sub-tree $C_{l_1,\ldots,l_{t_2}}$ with the similar properties as before.

\subsection{Initial regularization of the set}\label{sec:3.1}

Here we don't impose any conditions on $A\subset\mathbb R$. We will extract from $A$ a large subset with a ``regular tree structure''. Furthermore, at each level, we retain a pair of separated nodes.

\begin{lemma}\label{result:3.1} Suppose $\delta=2^{-m}$ ($m$ large) and let $A\subset[1,2]$ be a collection of $\delta$-separated points and $|A|=\delta^{-\sigma}$ for some fixed $0<\sigma<1$. Let $T$ be a large constant and assume $m=Tm_1$. Thus we have a sequence of scales:
\[
\delta=\Delta^{m_1}<\Delta^{m_1-1}<\cdots<\Delta^1<\Delta^0=1,
\]
where $\Delta=2^{-T}$. Then there are a subset $A_1\subset A$ and a subset
$S\subset\{0,\ldots,m_1-1\}$, i.e., a subsequence of scales $\{\Delta^s\}_{s\in S}$,
such that for each

$s\in S$, there is a dyadic integer $2\leq R_s<2^{n_s-sT}$, where
$sT<n_s<(s+1)T-4$, with the following properties:

\begin{enumerate}[label=(\roman*)]
\item\label{item:3-i} ``$s\notin S$ is bad'': if $I$ is a $2^{-sT}=\Delta^s$-interval, then there is $\leq1$ $2^{-(s+1)T}=\Delta^{s+1}$-interval $J\subset I$ such that $A_1\cap J\ne\varnothing$. In other word,
\[
|A_1\cap I|_{\Delta^{s+1}}=1\quad\text{for }I\in\mathcal D_{\Delta^s}(A_1),
\]
so the branching number at level $s$ is exactly $1$.

(Here $|E|_\Delta=N(E,\Delta)$ and
$\mathcal D_\Delta(E):=\{\Delta\text{-interval }I:E\cap I\ne\varnothing\}$.)

\item\label{item:3-ii} ``$s\in S$ is structured'': if $I$ is a $2^{-sT}=\Delta^s$-interval,
$A_1\cap I\ne\varnothing$, then the number of $2^{-n_s}$-intervals $J\subset I$ such that $A_1\cap J\ne\varnothing$ is $R_s$ and each such $2^{-n_s}$-interval $J$ contains a single $2^{-(s+1)T}=\Delta^{s+1}$-interval $J'$ that intersects $A_1$. Hence, also
\begin{equation}
N(A_1\cap I,2^{-(s+1)T})=|A_1\cap I|_{\Delta^{s+1}}\sim R_s.
\tag{3.1}\label{eq:3.1}
\end{equation}
So the branching number at level $s$ is $\sim R_s$. Furthermore, there exists a pair of $2^{-n_s}$-interval $J_1,J_2\subset I$ such that $A_1\cap J_1\ne\varnothing$, $A_1\cap J_2\ne\varnothing$, and
\begin{equation}
2\cdot2^{-n_s}<\operatorname{dist}(J_1,J_2)<10\cdot2^{-n_s}.
\tag{3.2}\label{eq:3.2}
\end{equation}

\item
\begin{equation}
|A_1|=\prod_{s\in S}R_s>(cT^{-2})^{m_1}|A|>
\delta^{\frac{\log T}{T}}\cdot\delta^{-\sigma}.
\tag{3.3}\label{eq:3.3}
\end{equation}
\end{enumerate}
\end{lemma}

\begin{proof} The construction is the standard pigeonholing from bottom to top.
We start at the bottom of the tree, considering for each $\Delta^{m_1-1}$-interval $I$ such that $I\cap A\ne\varnothing$ the number of $\Delta^{m_1}=\delta$-intervals $J\subset I$ with $J\cap A\ne\varnothing$, i.e., $|A\cap I|_{\Delta^{m_1}}$. If $|A\cap I|_{\Delta^{m_1}}\leq10^3$, fix one such $\Delta^{m_1}$-interval $J=J_1$.

If $|A\cap I|_{\Delta^{m_1}}>10^3$, introduce the largest integer (depending on $I$)
\[
(m_1-1)T<n<m_1T-4
\]
for which there exists a pair of $2^{-n}$-intervals $J_1,J_2\subset I$, $J_1\cap A\ne\varnothing$, $J_2\cap A\ne\varnothing$ that satisfies
\[
2\cdot2^{-n}<\operatorname{dist}(J_1,J_2)<10\cdot2^{-n}.
\]
Notice that this condition implies that each $2^{-n}$-interval $J\subset I$ contains $\sim1$ $\Delta^{m_1}$-intervals. Otherwise, if there are too many $\Delta^{m_1}$-intervals contained in some $J$, then we may find a larger $n'>n$, and a pair of $2^{-n'}$-intervals $J'_1,J'_2$

$\subset J\subset I$ satisfies $\operatorname{dist}(J'_1,J'_2)\sim2^{-n'}$, which contradicts to the maximization of $n$. (Such $n$ exists by the same argument when $n\sim(m_1-1)T$.) Therefore,
\[
N(A\cap I,2^{-n})\sim N(A\cap I,\Delta^{m_1}).
\]
\begin{center}
\begin{tikzpicture}[x=1cm,y=1cm]
\draw (0,1)--(6,1);\draw (0,1)--(0,1.15);\draw (6,1)--(6,1.15);
\node[above] at (3,1) {$J$};\node[right] at (6,1) {$2^{-n}$};
\draw[cyan,thick] (0,.45)--(1.8,.45);\draw[cyan,thick] (4.2,.45)--(6,.45);
\node[above,text=cyan] at (.9,.45) {$J'_1$};\node[above,text=cyan] at (5.1,.45) {$J'_2$};
\draw[cyan,<->] (1.8,.35)--(4.2,.35);\node[below,text=cyan] at (3,.35) {$\operatorname{dist}\sim2^{-n'}$};
\node[right,text=cyan] at (6,.45) {$2^{-n'}$};
\foreach \x in {0,.5,1,1.5,2.5,3.4,4.4,4.9,5.4,5.9}{\draw (\x,-.5)--(\x+.2,-.5);\draw (\x,-.5)--(\x,-.4);\draw (\x+.2,-.5)--(\x+.2,-.4);}
\node[right] at (6.2,-.5) {$\Delta^{m_1}$};
\end{tikzpicture}
\end{center}
Next we uniformize the branching number $N(A\cap I,2^{-n})$. Define $R$ to be the dyadic integer such that
\[
N(A\cap I,2^{-n})\sim R<2^{n-(m_1-1)T}.
\]
Obviously, the integers $n$ and $R$ (depending on $I$) take $\lesssim T$ values. By pigeonholing principle, we may introduce a subset $A^{(m_1-1)}\subset A$,
\[
|A^{(m_1-1)}|>10^{-3}T^{-2}|A|
\]
\noindent{Annotation to $10^{-3}$: this follows from the first case.}

satisfying one of the following alternatives:
\begin{itemize}
\item If $I$ is a $\Delta^{m_1-1}$-interval, there is $\leq1$ $\Delta^{m_1}$-interval $J\subset I$ with $J\cap A^{(m_1-1)}\ne\varnothing$. In this case, $m_1-1\notin S$.
\item If $I$ is a $\Delta^{m_1-1}$-interval and $A^{(m_1-1)}\cap I\ne\varnothing$, then \paperref{item:3-ii}{(ii)} holds for some $n_{m_1-1}:=n$ and $R_{m_1-1}:=R$. In this case, $m_1-1\in S$.
\end{itemize}

Next, repeat the construction for the set $A^{(m_1-1)}$ and considering
$2^{-(m_1-2)T}=\Delta^{m_1-2}$-intervals $I$ with $A^{(m_1-1)}\cap I\ne\varnothing$ and
$|A^{(m_1-1)}\cap I|_{\Delta^{m_1-1}}$.
At this stage every $\Delta^{m_1-1}$-interval contains the same number of retained leaves. We therefore perform the same selections, and obtain $A^{(m_1-2)}$ as the intersection of $A^{(m_1-1)}$ and a collection of $\Delta^{m_1-1}$-intervals with
\[
|A^{(m_1-2)}|>10^{-3}T^{-2}|A^{(m_1-1)}|,
\]
and $A^{(m_1-2)}$ satisfies either \paperref{item:3-i}{(i)} or \paperref{item:3-ii}{(ii)} with $s=m_1-2$, for some
\[
(m_1-2)T<n_{m_1-2}<(m_1-1)T-4
\quad\text{and}\quad R_{m_1-2}<2^{n_{m_1-2}-(m_1-2)T}.
\]
Moreover, in second step we don't change the ``leaf'' attached to the vertex at level $m_1-1$, in other word,
\[
A^{(m_1-2)}\cap I=A^{(m_1-1)}\cap I
\]
for $\Delta^{m_1-1}$-interval $I$ intersecting $A^{(m_1-2)}$. Hence properties \paperref{item:3-i}{(i)} and \paperref{item:3-ii}{(ii)} at level $m_1-1$ remain preserved.
\begin{center}
\includegraphics[width=\linewidth]{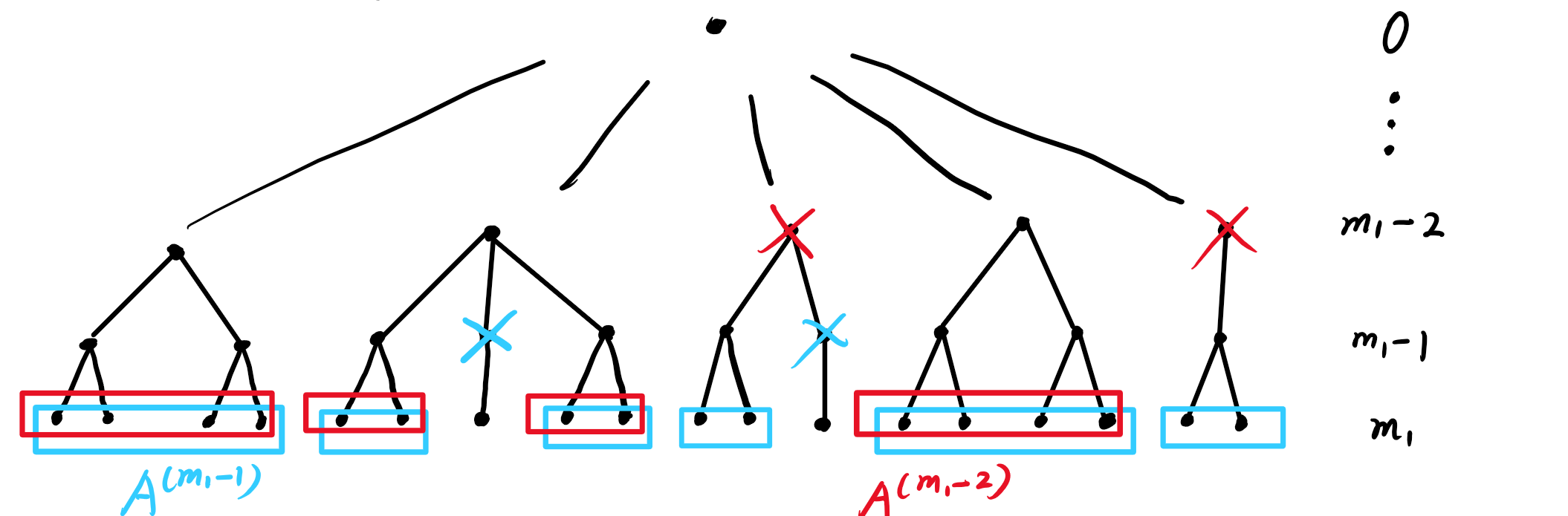}
\end{center}
\centerline{At each stage, we remove some vertices and their leaves.}

The continuation of the process is clear. We obtain
\[
A\supset A^{(m_1-1)}\supset A^{(m_1-2)}\supset\cdots\supset A^{(s)}\supset A^{(s-1)}\supset\cdots\supset A^{(1)}\supset A^{(0)},
\]
where
\begin{equation}
|A^{(s-1)}|>10^{-3}T^{-2}|A^{(s)}|
\tag{3.4}\label{eq:3.4}
\end{equation}
and
\[
A^{(s-1)}\cap I=A^{(s)}\cap I
\]
for any $2^{-sT}=\Delta^s$-interval $I$ intersecting $A^{(s-1)}$, and each $A^{(s)}$ satisfies either \paperref{item:3-i}{(i)} or \paperref{item:3-ii}{(ii)}. Let $A_1=A^{(0)}$, we get $A_1$ satisfying the desired properties \paperref{item:3-i}{(i)} and \paperref{item:3-ii}{(ii)}. Iteration \paperref{eq:3.4}{(\ref*{eq:3.4})} gives \paperref{eq:3.3}{(\ref*{eq:3.3})}:
\[
|A_1|=\prod_{s\in S}R_s>(10^{-3}T^{-2})^{m_1}|A|.
\]
\end{proof}

\subsection{Use of the sumset property}\label{sec:3.2}

Let first sketch the idea of the construction.
\begin{center}
\fcolorbox{black!55}{white}{\begin{minipage}{.9\linewidth}
\centering
Tree structure from \paperref{result:3.1}{Lemma~\ref*{result:3.1}} $+$ have small sumset
\[
\Downarrow
\]
Tree with ``partially maximal'' (in some sense) branching numbers
\[
\Downarrow
\]
Tree with multiscale porosity property
\[
\Downarrow
\]
Multiscale porosity sub-tree with large enough branching numbers
\end{minipage}}
\end{center}

\textcolor{red}{
In this subsection, we first employ the \paperref{result:3.1}{regularization lemma} to construct a tree, and then use the \paperref{result:1.3}{Pl\"unnecke–Ruzsa inequality} combined with binary tree encoding to obtain a subtree with partially saturated branching number. At this stage, the tree already exhibits "porousity" at levels with a small number of nodes, and there is a substantial number of such "good levels".
}

From now on, we assume $A$ satisfies the non-concentration property:
\[
|A\cap I|<\rho^\kappa|A|\quad\text{if }\delta<\rho<\delta^{\varepsilon_0}
\]
for any $\rho$-interval, some $\kappa>0$ and $\varepsilon_0=\varepsilon_0(\sigma,\kappa)>0$ sufficiently small. We also assume $|A|=\delta^{-\sigma}$, $\delta=\Delta^{m_1}$ as before. In addition, we assume a doubling estimate
\[
N(A+A,\delta)<\delta^{-\sigma-\varepsilon_1}
\quad(=\delta^{-\varepsilon_1}|A|),
\]
with $\varepsilon_1>0$ sufficiently small depending on $\sigma$ and $\kappa$ and to be specified.

\noindent\textbf{Step 1.} By \paperref{result:3.1}{Lemma~\ref*{result:3.1}}, we obtain a dense subset $A_1\subset A$ with regular tree structure. Let
\[
\begin{gathered}
S=\{s_1<s_2<\cdots<s_t\}\\
\Updownarrow\\
\delta\leq\Delta^{s_t}<\Delta^{s_{t-1}}<\cdots<\Delta^{s_2}<\Delta^{s_1}\leq1.
\end{gathered}
\]
For each $c\in\{0,1\}^t$, we will introduce an element $x_c\in A_1\subset A$.
Since $s_1\in S$, it follows from \paperref{item:3-ii}{(ii)} that, there exists a pair of $2^{-n_{s_1}}$-intervals $I_0,I_1$ intersecting $A_1$ such that
\[
2\cdot2^{-n_{s_1}}<\operatorname{dist}(I_0,I_1)<10\cdot2^{-n_{s_1}}.
\]
Moreover, each $I_i$ contains a single $\Delta^{s_1+1}$ $(\geq\Delta^{s_2})$-interval which intersects $A_1$. Denoted by $I'_0\subset I_0$ and $I'_1\subset I_1$ such $\Delta^{s_1+1}$-interval. In addition, we may assume $I_0$ is located to the left of $I_1$, denoted by $I_0\prec I_1$.
Again by \paperref{item:3-ii}{(ii)}, there exist pairs $I_{0,0},I_{0,1}\subset I'_0$ and $I_{1,0},I_{1,1}\subset I'_1$ of $2^{-n_{s_2}}$-intervals satisfying
\[
\begin{gathered}
2\cdot2^{-n_{s_2}}<\operatorname{dist}(I_{0,0},I_{0,1})<10\cdot2^{-n_{s_2}},\\
2\cdot2^{-n_{s_2}}<\operatorname{dist}(I_{1,0},I_{1,1})<10\cdot2^{-n_{s_2}},\\
\text{and } I_{0,0}\prec I_{0,1}\prec I_{1,0}\prec I_{1,1}.
\end{gathered}
\]
Moreover, each $I_{i,j}$ contains a single $\Delta^{s_2+1}$ $(\geq\Delta^{s_3})$-interval. Continuing the construction, we obtain eventually $2^{-n_{s_t}}$-intervals $I_c=I_{c_1,\ldots,c_t}$ for $c=(c_1,\ldots,c_t)\in\{0,1\}^t$ such that $I_c\cap A_1\ne\varnothing$ and with the following properties:
\begin{equation}
I_c\prec I_{c'}\quad\text{whenever }c\prec c'\text{ in lexicographical order},
\tag{3.5}\label{eq:3.5}
\end{equation}
and, if $c_1=c'_1,\ldots,c_\tau=c'_\tau$ but $c_{\tau+1}\ne c'_{\tau+1}$, then
\begin{equation}
2\cdot2^{-n_{s_{\tau+1}}}<\operatorname{dist}(I_c,I_{c'})<10\cdot2^{-n_{s_{\tau+1}}}.
\tag{3.6}\label{eq:3.6}
\end{equation}
Take

\begin{equation}
x_c\in I_c\cap A_1\quad\text{for }c\in\{0,1\}^t.
\tag{3.7}\label{eq:3.7}
\end{equation}
\textcolor{blue}{
Let's explain this construction intuitively. Firstly, we find a regular tree structure in $A$. Moreover, we may obtain a binary tree with $t$ levels, which is a little porous at each level
($2\cdot\text{diameter}\leq\text{separation}\leq10\cdot\text{diameter}$).
Hence, for any $t$-dimensional Boolean vector $c$, it could uniquely determines a path in the binary tree: at level $\tau$, a value $1$ means go right, $0$ means go left. And therefore $c$ uniquely determines a leaf $x_c$ in $A$.}
\begin{center}
\includegraphics[width=\linewidth]{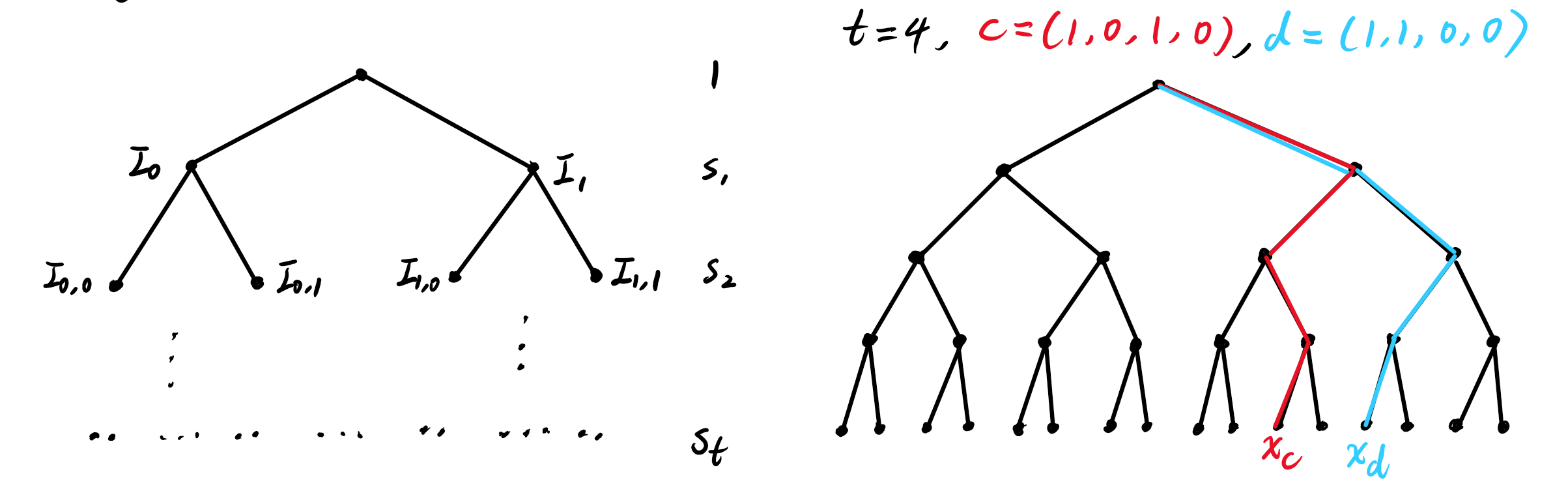}
\end{center}
Denoted by $\widetilde A=2^TA=A+\cdots+A$ the $2^T$-fold sumset of $A$. It follows from the \paperref{result:1.3}{Pl\"unnecke--Ruzsa inequality}, \paperref{result:3.1}{Lemma~\ref*{result:3.1}}, we get
\begin{equation}
N(\widetilde A,\delta)<(\delta^{-\varepsilon_1})^{2^T}N(A,\delta)
 =\delta^{-\sigma-2^T\varepsilon_1}.
\tag{3.8}\label{eq:3.8}
\end{equation}

Next, we construct a map
\[
\varphi:\Omega:=\prod_{s\in S}\{1,\ldots,2^{n_s-sT}\}
 \longrightarrow\widetilde A_1=\underbrace{A_1+\cdots+A_1}_{2^T}
\]
which maps $\Omega$ to $\delta$-separated points in $\widetilde A_1$.
Let $\bar k=(k_1,\ldots,k_t)\in\Omega$, $1\leq k_\tau\leq2^{n_{s_\tau}-s_\tau T}$. Take
\[
\varphi(\bar k)=x_{c^{(1)}}+x_{c^{(2)}}+\cdots+x_{c^{(2^T)}},
\]
where for $u=1,\ldots,2^T$, the element $c^{(u)}\in\{0,1\}^t$ is defined as follows:
\[
c^{(u)}=(c^{(u)}_1,\ldots,c^{(u)}_t)
\quad\text{and}\quad
c^{(u)}_\tau=\begin{cases}
1,&\text{if }u\leq k_\tau,\\
0,&\text{if }k_\tau<u\leq2^T.
\end{cases}
\]
\par\noindent{\rule{\linewidth}{.4pt}}\par
Before continuing, let's explain the geometric meaning of $\varphi$.
Our aim is to select many levels at which the branching number $R_s\sim2^{n_s-sT}$.
Recall $\prod R_s\sim|A|\approx|2^TA|$ and $\prod R_s\leq\prod2^{n_s-sT}$,
we need to define an injective
\[
\varphi:\prod\{1,\ldots,2^{n_s-sT}\}\longrightarrow2^TA.
\]
For any $t$-dimensional vector $\bar k=(k_1,\ldots,k_t)$, where
$1\leq k_\tau\leq2^{n_{s_\tau}-s_\tau T}$, the map $\varphi$ uses $\bar k$ to encode
$2^T$ $t$-dim Boolean vectors $c^{(u)}\Longleftrightarrow2^T$ paths in the binary tree
$\Longleftrightarrow2^T$ leaves $x_c$ in $A$, and then sums them up.
The encoding is as following: at level $\tau$, the first $k_\tau$ vector go right,
and the remaining vectors go left. For a schematic example with five summands (and hence five Boolean vectors), write
\[
2^T=5,\qquad t=5,\qquad\bar k=(3,2,4,1,2).
\]
\begin{center}
\textcolor{red}{$c^{(1)}$}\quad
\textcolor{cyan}{$c^{(2)}$}\quad
\textcolor{orange}{$c^{(3)}$}\quad
\textcolor{green!65!black}{$c^{(4)}$}\quad
\textcolor{violet}{$c^{(5)}$}.
\end{center}

\begin{center}
\includegraphics[width=\linewidth]{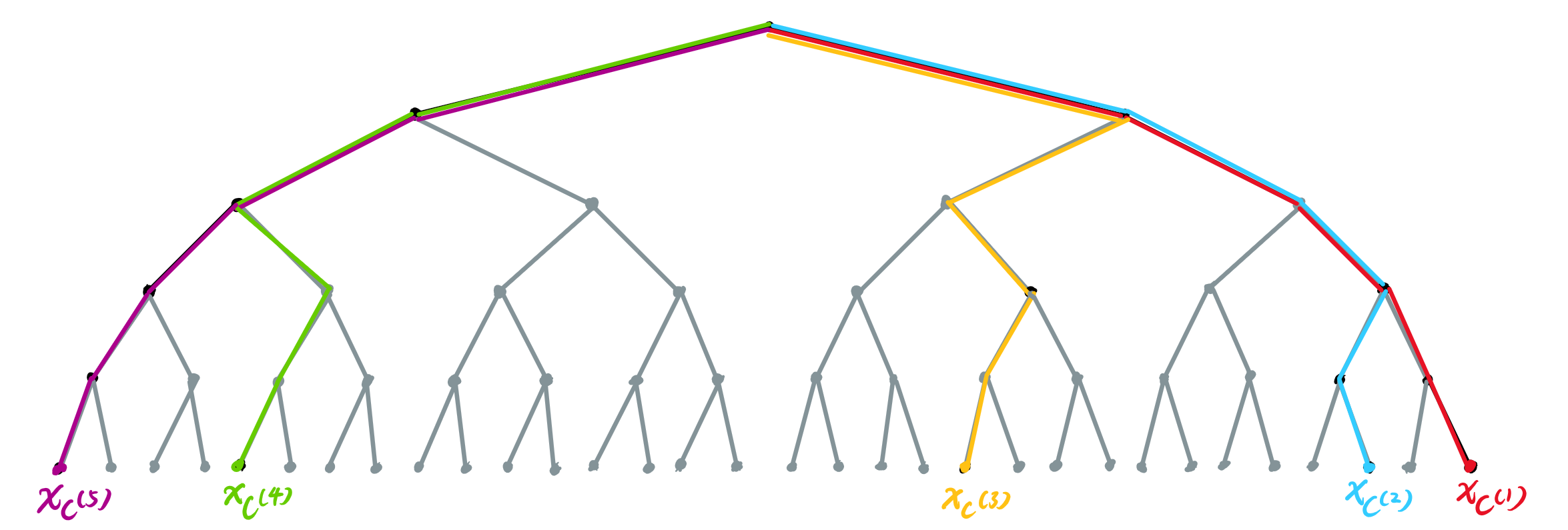}
\end{center}
Finally, we show why $\varphi$ is injective. For two $t$-dim vectors $\bar k\prec\bar l$ which first differ at level $\tau+1$, let the corresponding encoded families of $2^T$ vectors be $c^{(u)}$ and $d^{(u)}$. Then, each pair of $c^{(u)}$ and $d^{(u)}$ follow the same path up to the level $\tau$. At level $\tau+1$, the situation differs:
the \textcolor{red}{first $k_{\tau+1}$} vectors still \textcolor{red}{go right},
and the vectors \textcolor{cyan}{after the first $l_{\tau+1}$ go left};
but for the \textcolor{orange}{middle $l_{\tau+1}-k_{\tau+1}$} vectors,
\textcolor{orange}{$c^{(u)}$ go left} while \textcolor{orange}{$d^{(u)}$ go right}.
\begin{center}
\includegraphics[width=.45\linewidth]{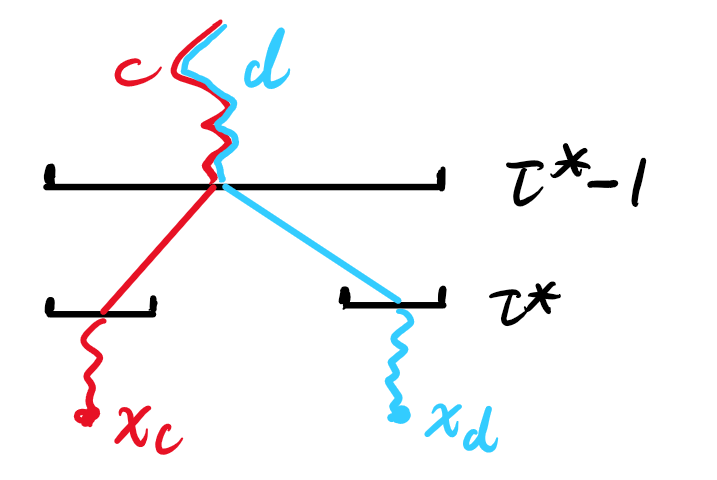}
\end{center}
Observed that: if two paths \textcolor{red}{$c$} and \textcolor{cyan}{$d$} first diverge at level $\tau^*$, then the distance between $x_c$ and $x_d$ satisfies:
\begin{itemize}
\item $|x_c-x_d|\geq\text{separation at level }\tau^*$,
\item $|x_c-x_d|\leq 2\cdot\text{diameter}+\text{separation at level }\tau^*
\leq 12\cdot\text{diameter at level }\tau^*$,
\end{itemize}
since the separation-to-diameter ratio lies in $[2,10]$ at each branching level by our construction.
Take $U_0=\{k_{\tau+1}+1,\ldots,l_{\tau+1}\}$. Thus,
\[
|\varphi(\bar k)-\varphi(\bar l)|
\geq\left|\sum_{u\in U_0}(x_{c^{(u)}}-x_{d^{(u)}})\right|
-\left|\sum_{u\notin U_0}(x_{c^{(u)}}-x_{d^{(u)}})\right|.
\]
The terms are divided into two classes:
\begin{itemize}
\item The first class consists of those $u$ whose paths (for $c^{(u)},d^{(u)}$) separate at level $\tau+1$, and $\#u=l_{\tau+1}-k_{\tau+1}\geq1$. Moreover,
\[
|x_{c^{(u)}}-x_{d^{(u)}}|\geq\text{separation at level }\tau+1.
\]
\item The second class consists of the remaining $u$. Group them by their first differing level $\tau'\geq\tau+2$, omitting identical paths. Note that
\[
\begin{aligned}
\text{the second term}
&\lesssim(\text{maximal branching number at level }\tau+2)\\
&\qquad\times(\text{diameter at level }\tau+2)\\
&\lesssim\text{diameter at level }\tau+1.
\end{aligned}
\]
\end{itemize}

\[
\text{the second term}
 \leq15\cdot2^{-s_{\tau+2}T}
 <\frac{15}{16}2^{-n_{s_{\tau+1}}}
 <\text{the first term}.
\]
So, $|\varphi(\bar k)-\varphi(\bar l)|>0$ and $\varphi$ is injective.
\par\noindent{\rule{\linewidth}{0.4pt}}

Let $\bar k=(k_1,\ldots,k_t)$ and $\bar l=(l_1,\ldots,l_t)\in\Omega$ be distinct elements and
\[
k_1=l_1,\ldots,k_\tau=l_\tau\quad\text{but}\quad k_{\tau+1}<l_{\tau+1}.
\]
Thus
\[
\varphi(\bar k)=x_{c^{(1)}}+\cdots+x_{c^{(2^T)}},
\qquad
\varphi(\bar l)=x_{d^{(1)}}+\cdots+x_{d^{(2^T)}},
\]
with
\[
c_{\tau'}^{(u)}=\begin{cases}1,&\text{if }u\leq k_{\tau'},\\0,&\text{if }u>k_{\tau'},\end{cases}
\qquad
d_{\tau'}^{(u)}=\begin{cases}1,&\text{if }u\leq l_{\tau'},\\0,&\text{if }u>l_{\tau'}.\end{cases}
\]
Hence,
\begin{itemize}
\item
\[
\begin{aligned}
c_{\tau'}^{(u)}&=d_{\tau'}^{(u)}&&\text{if }\tau'\leq\tau,\quad1\leq u\leq2^T,\\
c_{\tau+1}^{(u)}&=d_{\tau+1}^{(u)}&&\text{if }u\leq k_{\tau+1}\text{ or }u>l_{\tau+1},\\
c_{\tau+1}^{(u)}&<d_{\tau+1}^{(u)}&&\text{if }k_{\tau+1}<u\leq l_{\tau+1}.
\end{aligned}
\]
\item Since $x_{c^{(u)}}\in I_{c^{(u)}}$, $x_{d^{(u)}}\in I_{d^{(u)}}$, it follows that
\end{itemize}

\[
\begin{aligned}
|x_{c^{(u)}}-x_{d^{(u)}}|
&\leq\operatorname{diam}\bigl(I_{(c_1^{(u)},\ldots,c_{\tau^*}^{(u)})}\bigr)
+\operatorname{diam}\bigl(I_{(d_1^{(u)},\ldots,d_{\tau^*}^{(u)})}\bigr)\\
&\quad+\operatorname{dist}\bigl(I_{(c_1^{(u)},\ldots,c_{\tau^*}^{(u)})},I_{(d_1^{(u)},\ldots,d_{\tau^*}^{(u)})}\bigr)\\
&<14\cdot2^{-n_{s_{\tau^*}}},
\end{aligned}
\]
where
\[
\tau^*=\min\{\tau'\geq\tau+1:c_{\tau'}^{(u)}\neq d_{\tau'}^{(u)}\}.
\]
Therefore, we obtain
\[
\begin{aligned}
|\varphi(\bar k)-\varphi(\bar l)|
&=\left|\sum_{u=1}^{2^T}(x_{c^{(u)}}-x_{d^{(u)}})\right|\\
&\geq\left|\sum_{u=k_{\tau+1}+1}^{l_{\tau+1}}(x_{c^{(u)}}-x_{d^{(u)}})\right|
-\left|\sum_{\substack{u\leq k_{\tau+1}\\\text{or }u>l_{\tau+1}}}(x_{c^{(u)}}-x_{d^{(u)}})\right|\\
&>(l_{\tau+1}-k_{\tau+1})\cdot2\cdot2^{-n_{s_{\tau+1}}}
-\sum_{\tau'\geq\tau+2}\sum_{\substack{u\\\tau^*(u)=\tau'}}|x_{c^{(u)}}-x_{d^{(u)}}|.
\end{aligned}
\]
\textcolor{cyan}{For the first sum: $c_1^{(u)}=d_1^{(u)},\ldots,c_\tau^{(u)}=d_\tau^{(u)}$ and $c_{\tau+1}^{(u)}<d_{\tau+1}^{(u)}$. For the second sum: $\tau^*\geq\tau+2$.}

(To estimate $\#\{u:\tau^*(u)=\tau'\}$, note that by definition $\tau^*$, $c_{\tau'}^{(u)}\neq d_{\tau'}^{(u)}$ for such $u$. Thus
\[
\text{either }u\leq k_{\tau'}\ \&\ u>l_{\tau'}
\quad\text{or}\quad u>k_{\tau'}\ \&\ u\leq l_{\tau'},
\]
\[
\Rightarrow\quad\#\{u:\tau^*(u)=\tau'\}\leq|k_{\tau'}-l_{\tau'}|\leq\max\{k_{\tau'},l_{\tau'}\}.)
\]

Thus, 
\[
\begin{aligned}
|\varphi(\bar k)-\varphi(\bar l)|
&>2(l_{\tau+1}-k_{\tau+1})2^{-n_{s_{\tau+1}}}
-14\sum_{\tau'\geq\tau+2}\max\{k_{\tau'},l_{\tau'}\}2^{-n_{s_{\tau'}}}\\
&>2^{-n_{s_{\tau+1}}+1}-14\sum_{\tau'\geq\tau+2}2^{n_{s_{\tau'}}-s_{\tau'}T}\cdot2^{-n_{s_{\tau'}}}\\
&>2^{-n_{s_{\tau+1}}+1}-15\cdot2^{-s_{\tau+2}T}.
\end{aligned}
\]
(Recall $n_{s_{\tau+1}}<(s_{\tau+1}+1)T-4\leq s_{\tau+2}T-4$
\[
\Rightarrow\quad 2^{-n_{s_{\tau+1}}}>2^{-s_{\tau+2}T+4}=16\cdot2^{-s_{\tau+2}T}.)
\]
\[
|\varphi(\bar k)-\varphi(\bar l)|>2^{-n_{s_{\tau+1}}}>\delta.
\]
Thus, $\varphi$ maps $\Omega$ to $\delta$-separated points in $\widetilde A_1$, and then
\[
N(\widetilde A_1,\delta)\geq N(\varphi(\Omega),\delta)=|\Omega|=\prod_{s\in S}2^{n_s-sT}.
\]
Assume $2^T<\varepsilon_1^{-1/2}$. We obtain from \paperref{eq:3.8}{(\ref*{eq:3.8})} that
\[
\prod_{s\in S}2^{n_s-sT}\leq N(\widetilde A,\delta)<\delta^{-\varepsilon_1\cdot2^T}|A|.
\]

\[
\prod_{s\in S}2^{n_s-sT}
\overset{\paperref{eq:3.3}{(\ref*{eq:3.3})}}{<}\delta^{-2^{-T}}\cdot\delta^{-\frac{\log T}{T}}\cdot|A_1|
<\delta^{-\frac{\log T}{T}}|A_1|,
\]
where $|A_1|=\prod_{s\in S}R_s$. Define
\[
S_1=\{s\in S:2^{n_s-sT}>R_s>2^{(1-T^{-1/2})(n_s-sT)}\}.\tag{3.9}\label{eq:3.9}
\]
It follows that,
\[
\begin{aligned}
\prod_{s\in S\setminus S_1}2^{(1-T^{-1/2})(n_s-sT)}
&\geq\prod_{s\in S\setminus S_1}R_s\\
&>\prod_{s\in S}R_s\cdot\prod_{s\in S_1}2^{-(n_s-sT)}\\
&=|A_1|\cdot\prod_{s\in S_1}2^{-(n_s-sT)}\\
&>\delta^{\frac{\log T}{T}}\cdot\prod_{s\in S\setminus S_1}2^{n_s-sT}.
\end{aligned}
\]
\[
\Rightarrow\quad\prod_{s\in S\setminus S_1}2^{n_s-sT}<\delta^{-\frac{\log T}{T^{1/2}}}.\tag{3.10}\label{eq:3.10}
\]
Now, we replace $A_1$ further with another subset $A_2$ satisfying again properties \paperref{item:3-i}{(i)} and \paperref{item:3-ii}{(ii)} in \paperref{result:3.1}{Lemma~\ref*{result:3.1}}, but with $S$ replaced with $S_1$ and satisfying
\[
|A_2|>\delta^{2(\log T)/T^{1/2}}|A_1|.
\]
Proceed as follows.

Start from $A_1$, if $s_1\in S_1$, take $A_1^{(1)}=A_1$. If $s_1\notin S_1$, by \paperref{item:3-ii}{(ii)} and pigeonholing principle, we may choose in each $\Delta^{s_1}$-interval $I$ with $I\cap A_1\neq\varnothing$ a $\Delta^{s_1+1}$-interval $J_I\subset I$ such that $J_I\cap A_1\neq\varnothing$ and $|J_I\cap A_1|\geq R_{s_1}^{-1}|I\cap A_1|$. Let
\[
A_1^{(1)}=\bigcup_I(A_1\cap J_I),
\]
hence
\[
|A_1^{(1)}|=|A_1|_{\Delta^{s_1}}|A_1\cap J_I|
\geq\frac1{R_{s_1}}|A_1|_{\Delta^{s_1}}|I\cap A_1|
=\frac1{R_{s_1}}|A_1|.
\]
Next we consider $s_2>s_1$, $s_2\in S$. If $s_2\in S_1$, let $A_1^{(2)}=A_1^{(1)}$. If $s_2\notin S_1$, replace $A_1^{(1)}$ with
\[
A_1^{(2)}=\bigcup_{I\in\mathcal D_{\Delta^{s_2}}(A_1^{(1)})}(A_1^{(1)}\cap J_I),
\quad\text{where }|J_I|=\Delta^{s_2+1},
\]
again, $|A_1^{(2)}|\geq\dfrac1{R_{s_2}}|A_1^{(1)}|$.
The continuation of the process is clear. Set $A_2=A_1^{(t)}$. The set $A_2\subset A_1$ eventually obtained satisfies, by \paperref{eq:3.10}{(\ref*{eq:3.10})},
\[
\begin{aligned}
|A_2|&\geq\left(\prod_{s\in S\setminus S_1}\frac1{R_s}\right)\cdot|A_1|\\
&>\delta^{(\log T)/T^{1/2}}\cdot\delta^{(\log T)/T}\cdot|A|\\
&>\delta^{(\log T)/T^{1/2}}|A|,
\end{aligned}\tag{3.11}\label{eq:3.11}
\]
and the following properties:
\begin{itemize}
\item If $s\notin S_1$ and $I$ is a $\Delta^s$-interval with $A_2\cap I\ne\varnothing$, there is $\leq1$ $\Delta^{s+1}$-interval $J\subset I$ such that $A_2\cap J\neq\varnothing$, i.e.,
\[
|A_2\cap I|_{\Delta^{s+1}}=N(A_2\cap I,\Delta^{s+1})=1.\tag{3.12}\label{eq:3.12}
\]
\item If $s\in S_1$ and $I$ is a $\Delta^s$-interval with $A_2\cap I\neq\varnothing$, there are $R_s$ $2^{-n_s}$-intervals $J\subset I$ intersecting $A_2$, where
\[
2^{(1-T^{-1/2})(n_s-sT)}\leq R_s<2^{n_s-sT}.\tag{3.13}\label{eq:3.13}
\]
Furthermore, each $2^{-n_s}$-interval $J\subset I$ such that $J\cap A_2\neq\varnothing$ contains a single $\Delta^{s+1}$-interval $J'$ intersecting $A_2$, i.e.,
\end{itemize}

\[
|A_2\cap I|_{\Delta^{s+1}}=N(A_2\cap I,\Delta^{s+1})\sim R_s.\tag{3.14}\label{eq:3.14}
\]
From now on, extend the notation by setting $R_s=1$ and $n_s=sT$ for $s\notin S_1$. 

\textcolor{blue}{\paperref{eq:3.13}{(\ref*{eq:3.13})} and \paperref{eq:3.14}{(\ref*{eq:3.14})} tell us that at the intermediate scale $2^{-n_s}$, the number of branches is already near saturation/partial maximal; furthermore, within each interval at this intermediate scale, the points in $A_2$ are, in fact, concentrated at a finer scale. Consequently, thanks to the separation property at the intermediate scale—provided that this intermediate scale is sufficiently far from the fine scale—we can obtain porosity property at fine scale.}

As we mentioned above, the set $A_2$ has been equipped with porosity property at certain levels. Fix some $\beta=\beta(\sigma)\in(0,1)$ to be specified and define
\[
S_2=\{s\in S_1:(s+1)T-n_s<\beta T\}.
\]
\textcolor{cyan}{$ n_s-sT>(1-\beta)T \iff(s+1)T-n_s<\beta T$, i.e., $2^{-n_s}$ is closer to the fine scale for $s\in S_2$.}
Hence, from \paperref{eq:3.13}{(\ref*{eq:3.13})},
\[
2^{T-4}\geq2^{n_s-sT}>R_s>2^{(1-T^{-1/2})(1-\beta)T}\quad\text{if }s\in S_2.
\]
Also, by \paperref{eq:3.12}{(\ref*{eq:3.12})} and \paperref{eq:3.14}{(\ref*{eq:3.14})}, we get if $s\in\{0,\ldots,m_1-1\}\setminus S_2$, the $\Delta^{s+1}$-intervals $J$ intersecting $A_2$ can be partitioned into $2$ classes each of which consists of $2^{-(s+1-\beta)T}=\Delta^{s+1-\beta}$-separated intervals. Now, we obtain the porosity property at each level $s\notin S_2$.

Finally, we show the number of these ``good'' levels is large enough. Let $\{0,\ldots,m_1-1\}\setminus S_2=\{s_1<s_2<\cdots<s_{t_1}\}$. Then we have

\[
\begin{aligned}
\delta^{-1}=2^{m_1T}
&=\prod_{s=0}^{m_1-1}2^T
=\prod_{s\notin S_2}2^T\cdot\prod_{s\in S_2}2^T\\
&<2^{t_1T}\left(\prod_{s\in S_2}R_s\right)^{\frac1{(1-\beta)(1-T^{-1/2})}}\\
&<2^{t_1T}\cdot|A|^{1/[(1-\beta)(1-T^{-1/2})]}\\
&=2^{t_1T}\cdot\delta^{-\sigma/[(1-\beta)(1-T^{-1/2})]}\\
&<2^{t_1T}\cdot\delta^{-\frac{1+\sigma}{2}},
\end{aligned}
\]
whenever letting $\beta\ll1$ first, and $T=T(\beta,\sigma)=T(\sigma)\gg1$, such that
\[
\frac{2\sigma}{1+\sigma}<(1-\beta)(1-T^{-1/2}).
\]
For example, we can choose
\[
1-\beta>\frac12\left(1+\frac{2\sigma}{1+\sigma}\right)=\frac{1+3\sigma}{2(1+\sigma)},
\]
\[
\&\quad1-T^{-1/2}>\left(\frac{2\sigma}{1+\sigma}\right)\bigg/\left(\frac{1+3\sigma}{2(1+\sigma)}\right)=\frac{4\sigma}{1+3\sigma}.
\]

It follows that
\[
2^{t_1T}>\delta^{-\frac{1-\sigma}{2}}=2^{(\frac{1-\sigma}{2})m_1T}
\]
\[
\Rightarrow\quad m_1\geq t_1>\frac{1-\sigma}{2}m_1.\tag{3.15}\label{eq:3.15}
\]
In other word, $\{0,\ldots,m_1-1\}\setminus S_2$ is ``dense'' in $\{0,\ldots,m_1-1\}$.

\subsection{Hide bad levels}\label{sec:3.3}
\textcolor{red}{
Now, we construct from $A_2$ a new system of sets $B_{k_1,\ldots,k_j}$, $j\leq t_1$, with ``porosity property'' at each level. The approach is to record only the "good" levels, while hiding the "bad" levels that lie between any two good ones; doing so simply increases the number of branches.}

Set $s_0=-1$. Before imposing separation, the number of children between consecutive good levels is
\[
K_j^{\mathrm{old}}=\prod_{s_{j-1}<s\leq s_j}R_s,
\qquad 1\leq j\leq t_1.
\]
At the first level $s_1$, take the non-empty intersections $A_2\cap I$ with $2^{-n_{s_1}}$-intervals $I$. Each such set is contained in a single $\Delta^{s_1+1}$-interval. Partition these intervals into the two classes described above and retain exactly
\[
K_1=\max\{1,K_1^{\mathrm{old}}/2\}
\]
children, keeping their entire descendant subtrees.(The factor $1/2$ arises because, strictly speaking, while we know two separated points, the intervals containing these two points might end up being adjacent; in such cases, we need to discard one of the intervals.) Denote the retained sets by $B_{k_1}$, $1\leq k_1\leq K_1$. Then
\[
\operatorname{diam}B_{k_1}\leq\Delta^{s_1+1},\qquad
\operatorname{dist}(B_{k_1},B_{k_1'})>\Delta^{s_1+1-\beta}
\quad(k_1\ne k_1').
\]

Within each retained first-level set, repeat the same selection at level $s_2$, keeping exactly $K_2=\max\{1,K_2^{\mathrm{old}}/2\}$ children and their entire descendant subtrees. Continue in this way. At stage $j$, every parent retains exactly
\[
K_j=\max\{1,K_j^{\mathrm{old}}/2\}
\]
children.

After all steps, redefine each earlier node as the union of its attached nodes. The separation between children in different parents follows from the already imposed coarser separation, since
$\Delta^{s_{j-1}+1-\beta}>\Delta^{s_j+1-\beta}$. Thus we obtain a regular system $B_{k_1,\ldots,k_j}$, $1\leq k_i\leq K_i$, $j\leq t_1$, with the following properties. Write $B_\varnothing$ for the union of all its leaves.
\begin{itemize}
\item
\[
|B_\varnothing|=\left|\bigcup_{\bar k}B_{\bar k}\right|
\geq2^{-t_1}|A_2|\geq2^{-m_1}|A_2|.\tag{3.16}\label{eq:3.16}
\]
\item
\[
B_{k_1,\ldots,k_j}=\bigcup_{1\leq k\leq K_{j+1}}B_{k_1,\ldots,k_j,k}.\tag{3.17}\label{eq:3.17}
\]
\item
\[
\operatorname{diam}B_{k_1,\ldots,k_j}\leq\Delta^{s_j+1}=2^{-T\beta}\lambda_j.\tag{3.18}\label{eq:3.18}
\]
\item
\[
\operatorname{dist}(B_{k_1,\ldots,k_j},B_{k_1',\ldots,k_j'})>\lambda_j
\quad\text{if }\bar k\ne\bar k',\tag{3.19}\label{eq:3.19}
\]
where
\[
\lambda_j=2^{-(s_j+1-\beta)T}=\Delta^{s_j+1-\beta}.\tag{3.20}\label{eq:3.20}
\]
\end{itemize}

\subsection{Merge levels with few branches to obtain lots of effective levels}\label{sec:3.4}
From this point on, we follow the procedure in \paperref{sec:2}{Section~2}. Now we introduce a new system $C_{l_1,\ldots,l_s}$, $1\leq l_i\leq L_i$ and $1\leq s\leq t_2$. (Here, $\Delta^{-1}$ plays the same role as $Q$ before.) Define
\[
M=[2^{T^{2/3}}]\tag{3.21}\label{eq:3.21}
\]
and let
\[
10^3<D\in\mathbb Z_+,\tag{3.22}\label{eq:3.22}
\]
which is yet to be specified.

Starting from $t_1$, let $r_1\geq1$ be the minimal integer such that
\[
\prod_{j=t_1-r_1+1}^{t_1}K_j>M^{r_1}.
\]
Here $r_1$ is the number of branching levels in the block. This make sense, and we will explain the existence of $r_1$ later. We distinglish two cases, as before.

If $r_1\leq10^3$, identify the branching levels $t_1-r_1+1,\ldots,t_1$ to a single one, with branching
\[
L_{t_2}=\prod_{j=t_1-r_1+1}^{t_1}K_j>M^{r_1}.
\]
Use the original level-$t_1$ sets as the new children. Their separation is greater than $\mu_{t_2}=\lambda_{t_1}$, and their diameter is at most $\sigma_{t_2}\mu_{t_2}$, where
\[
\sigma_{t_2}=2^{-\beta T},\qquad
\sigma_{t_2}^{-1}=2^{\beta T}>M^{Dr_1},
\]
provided
\[
T>10^{10}\left(\frac D\beta\right)^3.\tag{3.23}\label{eq:3.23}
\]

In this case there is no additional loss of mass.

If $r_1>10^3$, take $q_1=\lfloor r_1/100\rfloor$. Minimality of $r_1$ gives
\[
\prod_{j=t_1-q_1+1}^{t_1}K_j\leq M^{q_1}.
\]
Identify the branching levels $t_1-r_1+1,\ldots,t_1-q_1$ to a single one, with branching
\[
\begin{aligned}
L_{t_2}
&=\prod_{j=t_1-r_1+1}^{t_1-q_1}K_j\\
&>\frac{M^{r_1}}{\prod_{j=t_1-q_1+1}^{t_1}K_j}
\geq M^{r_1-q_1}\geq M^{\frac{99}{100}r_1}.
\end{aligned}
\]
To obtain sufficiently small diameter, for every old vertex $v$ at level $t_1-q_1$, choose one descendant $w(v)$ at level $t_1$ and retain the set $B_{w(v)}$ as a new child. This is the formal meaning of the descending schematic labels
\[
B_{k_{t_1},\ldots,k_{t_1-r_1/100}}
\quad\longmapsto\quad
B_{k_{t_1},\ldots,k_{t_1-\frac{r_1}{100}},
\underbrace{1,\ldots,1}_{r_1/100}}
=C_{l_1,\ldots,l_{t_2}}.
\]

By \paperref{eq:3.18}{(\ref*{eq:3.18})} and \paperref{eq:3.19}{(\ref*{eq:3.19})}, distinct new children have separation greater than
$\mu_{t_2}=\lambda_{t_1-q_1}$ and diameter at most
$2^{-T\beta}\lambda_{t_1}=\sigma_{t_2}\mu_{t_2}$. Hence
\[
\begin{aligned}
\sigma_{t_2}^{-1}
&=2^{T\beta}\frac{\lambda_{t_1-q_1}}{\lambda_{t_1}}
=2^{T\beta}2^{(s_{t_1}-s_{t_1-q_1})T}\\
&\geq2^{T(q_1+\beta)}>M^{Dr_1},
\end{aligned}
\]
where the last inequality follows from \paperref{eq:3.23}{(\ref*{eq:3.23})}, since $r_1>10^3$ and $q_1\geq r_1/100-1$.

All reference-tree vertices of a fixed level carry equal mass. Thus this pruning retains a fraction at least
\[
\left(\prod_{j=t_1-q_1+1}^{t_1}K_j\right)^{-1}
\geq M^{-q_1}\geq M^{-r_1/100}>M^{-r_1/90}
\]
of the mass. Repeat, if possible, from old level $b_1=t_1-r_1$, and then from
$b_2=b_1-r_2$, etc. At a later stage, the selected level-$b_i$ descendant keeps the entire already constructed lower subtree attached to it. Since the same numbers of children are retained under every parent, the mass fractions remain uniform.

Suppose the procedure stops at level $t'=t_1-r_1-r_2-\cdots$, which means
\[
K_1\cdots K_{t'}\leq M^{t'}
\leq2^{T^{2/3}m_1}=\delta^{-T^{-1/3}}.
\]
Its equal-mass nodes and \paperref{eq:3.16}{(\ref*{eq:3.16})} give
\[
2^{-m_1}|A_2|\leq |B_\varnothing|
=(K_1\cdots K_{t'})|B_v|
\]
for any reference vertex $v$ of level $t'$. Such a node is contained in some $\Delta^{s_{t'}}$-interval $I$, so
\[
\begin{aligned}
|A\cap I|&\geq|A_2\cap I|\geq|B_v|
\geq\delta^{T^{-1/3}}2^{-m_1}|A_2|\\
&\geq\delta^{T^{-1/3}+T^{-1}+C(\log T)/T^{1/2}}|A|
>\delta^{2T^{-1/3}}|A|
\end{aligned}
\]
for $T$ sufficiently large. Here we used \paperref{eq:3.11}{(\ref*{eq:3.11})} with its absolute constants made explicit.

The non-concentration condition therefore implies
\[
\begin{aligned}
\text{either }\delta^{2T^{-1/3}}\leq(\Delta^{s_{t'}})^\kappa
&\quad\Longrightarrow\quad s_{t'}\leq2\kappa^{-1}T^{-1/3}m_1,\\
\text{or }\Delta^{s_{t'}}>\delta^{\varepsilon_0}
&\quad\Longrightarrow\quad s_{t'}<\varepsilon_0m_1.
\end{aligned}
\]
Since the good levels are distinct non-negative integers, $t'\leq s_{t'}+1$. Thus, for $m_1$ sufficiently large,
\[
t'\leq1+\max\{\varepsilon_0,2\kappa^{-1}T^{-1/3}\}m_1
<\frac{1-\sigma}{4}m_1,
\]
provided
\[
\varepsilon_0<\frac{1-\sigma}{8},\qquad
T^{-1/3}<\frac{1-\sigma}{16}\kappa.\tag{3.24}\label{eq:3.24}
\]
Together with \paperref{eq:3.15}{(\ref*{eq:3.15})}, this gives
\[
t'<\frac{1-\sigma}{4}m_1<\frac12t_1,\tag{3.25}\label{eq:3.25}
\]
and hence
\[
r_1+r_2+\cdots=t_1-t'>\frac12t_1>\frac{1-\sigma}{4}m_1.\tag{3.26}\label{eq:3.26}
\]
The same argument rules out stopping at the initial level $t_1$, and thus proves existence of the initial $r_1$.

Consequently, we obtain a tree $C_{l_1,\ldots,l_s}$, $1\leq s\leq t_2$, with the following properties:

\begin{itemize}
\item
\[
C_{l_1,\ldots,l_s}=\bigcup_{1\leq l\leq L_{s+1}}C_{l_1,\ldots,l_s,l}.\tag{3.27}\label{eq:3.27}
\]
\item
\[
\operatorname{dist}(C_{l_1,\ldots,l_s},C_{l_1',\ldots,l_s'})>\mu_s
\quad\text{if }\bar l\ne\bar l'.\tag{3.28}\label{eq:3.28}
\]
\item
\[
\operatorname{diam}C_{l_1,\ldots,l_s}\leq\sigma_s\mu_s,
\qquad\sigma_s\leq2^{-T\beta}.\tag{3.29}\label{eq:3.29}
\]
\item
\[
\prod_{s=1}^{t_2}(L_s\wedge\sigma_s^{-1/D})
>M^{\frac{99}{100}(r_1+r_2+\cdots)}=M^{\gamma m_1},\tag{3.30}\label{eq:3.30}
\]
where
\[
\gamma=\frac{99}{100m_1}(r_1+r_2+\cdots)>\frac{1-\sigma}{5}.
\]
\item
\[
\begin{aligned}
|C_\varnothing|
=\left|\bigcup_{\bar l}C_{\bar l}\right|
&>M^{-\frac1{90}(r_1+r_2+\cdots)}2^{-m_1}|A_2|\\
&>M^{-\frac\gamma{80}m_1}|A|.
\end{aligned}\tag{3.31}\label{eq:3.31}
\]
\end{itemize}

In what follows, it is convenient to replace our discrete sets with unions of $\delta$-intervals by defining
\[
A'=\{x\in\mathbb R:\operatorname{dist}(x,A)<\tfrac\delta2\}=A+I(0,\tfrac\delta2)
\]
\[
\&\quad C_{\bar l}'=\{x\in\mathbb R:\operatorname{dist}(x,C_{\bar l})<\tfrac\delta2\}=C_{\bar l}+I(0,\tfrac\delta2).
\]
Hence
\[
|A'|=\delta|A|=\delta^{1-\sigma}\quad\&\quad|C_{\bar l}'|=\delta|C_{\bar l}|.
\]
After harmless absolute changes in the diameter and separation constants, the thickened sets $C_{\bar l}'$ still satisfy \paperref{eq:3.27}{(\ref*{eq:3.27})}--\paperref{eq:3.31}{(\ref*{eq:3.31})}.

\section{Amplification by multiplication}\label{sec:4}

Let $C_{\bar l}'$ be the tree obtained in \paperref{sec:2}{Section~2} \& \paperref{sec:3}{3}. In this section, we start using the multiplicative structure and consider the form
\[
vC_\varnothing'+xC_\varnothing'
\]
where $v$ is fixed and $x$ is selected randomly from $[0,v]$. We will use the multiscale porosity property of $C_\varnothing'$ to obtain the averaged expansion estimate \paperref{eq:4.12}{(\ref*{eq:4.12})}.
\textcolor{red}{
A key observation is that, at a fixed scale, the porosity property implies that if the sums of two branches overlap, this overlap forces $x$ to be concentrated within a short interval. Meanwhile, the non-concentration property ensures that the probability of such an occurrence is very small. Finally, by leveraging the multiscale structure—even if some levels fail to yield a "power saving"—aggregating across all levels allows one to achieve a gain with high probability.}

Let's begin the proof. Assume $\delta^{1/4}<v\leq1$ and $\mu\in\mathcal P([0,v])$ satisfying the condition: for some $\tau>0$,
\[
\mu(I)<c\mu([0,\Delta])<c\left(\frac\Delta v\right)^\tau\tag{4.1}\label{eq:4.1}
\]
whenever $I\subset[0,v]$ is a $\Delta$-interval, $\delta<\Delta<v$. We show in the \paperref{sec:5}{next section} how to convert the measure $\mu$ from \paperref{result:0.7}{Theorem~\ref*{result:0.7}} to one satisfying the above property \paperref{eq:4.1}{(\ref*{eq:4.1})}.

Let $x\in[0,v]$ and $\mathbb E_x$ refer to the $\mu$-expectation. Our aim is to estimate from below
\[
\mathbb E_x[|vC_\varnothing'+xC_\varnothing'|].
\]
We use the \paperref{result:4.1}{following observation}.

\begin{lemma}\label{result:4.1}
Let $(J_r)_{1\leq r\leq R}$ and $(J_r')_{1\leq r\leq R}$ be two sequences of disjoint intervals in $\mathbb R$. Then
\[
\sum_{r,s=1}^R\chi_{J_r+J_s'}\leq CR.\tag{4.2}\label{eq:4.2}
\]
In other word, the sumsets $(\bigcup J_r)+(\bigcup J_r')=\bigcup(J_r+J_s')$ have $\lesssim R$ overlapping, compared to their total number $R^2$.
\end{lemma}

\begin{proof} Write $r\prec s$ if $|J_r|\leq|J_s'|$. Clearly, if we fix $r$, then the sets $(J_r+J_s')_{s\succ r}$ have bounded multiplicity. Indeed, if $x\in(J_r+J_{s_1}')\cap(J_r+J_{s_2}')$ for $s_1\succ r$, $s_2\succ r$, then
\[
x=x_r+x_{s_1}=x_r'+x_{s_2}\quad\Rightarrow\quad|x_{s_1}-x_{s_2}|\leq|J_r|.
\]

But $(J_s')_{1\leq s\leq R}$ are disjoint and $|J_r|\leq|J_{s_1}'|,|J_{s_2}'|$. Thus, such $s_2$ are finite independent of $R$. Therefore,
\[
\sum_{s\succ r}\chi_{J_r+J_s'}<C
\]
\[
\Rightarrow\quad\sum_{\substack{r,s\\r\prec s}}\chi_{J_r+J_s'}
=\sum_r\sum_{\substack{s\\s\succ r}}\chi_{J_r+J_s'}<CR.
\]
Similarly, we show that
\[
\sum_{\substack{r,s\\r\succ s}}\chi_{J_r+J_s'}<CR.
\]
Hence,
\[
\sum_{r,s}\chi_{J_r+J_s'}\leq\left(\sum_{r\succ s}+\sum_{r\prec s}\right)\chi_{J_r+J_s'}<CR.
\]
\end{proof}
Recall \paperref{eq:3.30}{(\ref*{eq:3.30})}, let $M_s=L_s\wedge\sigma_s^{-1/D}$ $(1\leq s\leq t_2)$ satisfy
\[
\prod_{s=1}^{t_2}M_s>M^{\gamma m_1}.\tag{4.3}\label{eq:4.3}
\]
We choose the parameter $D$ to satisfy
\[
D>\frac{20}{\tau}.\tag{4.4}\label{eq:4.4}
\]

\textcolor{blue}{The role of $M_s$}: The best-case scenario is none of the intervals overlap. We divide the $L_s$ intervals into several blocks, with each block containing $M_s\le L_s$ intervals. Fix two blocks and $x$, we can view each interval $vC_{k}'+xC_{l}'$ as a single point $v\xi_k+x\xi_l$; if two intervals overlap, it means the distance between their corresponding points is small. It follows that $x$ lies in an interval of length $\sim v\sigma_s$. The non-concentration of $\mu$ ensures that the probability of $x$ falling into such a "bad" interval is very small: $\lesssim\sigma_s^{\tau}$. Moreover, there are $\sim M_s^4$ pairs of interval. So the probability that the two intervals overlap is $\lesssim M_s^4\sigma_s^{\tau}$. This is the reason why we need to restrict $M_s$ to $\le\sigma^{-1/D}$.

\annot{Thus, blocks helps us capture the gains that arise when intervals are non-overlapping with high probability, while the losses incurred from overlaps between blocks are controled by \paperref{result:4.1}{lemma~\ref*{result:4.1}}.
}

Starting from $s=1$, we thus have $C_\varnothing'=\bigcup_{l_1\leq L_1}C_{l_1}'$, where $C_{l_1}'$ are contained in consecutive intervals of size $\sigma_1\mu_1$ and separation $\mu_1$. Write further
\[
C_\varnothing'=\bigcup_{r\leq L_1/M_1}D_r
\quad\&\quad D_r=\bigcup_{l_1\in\mathcal E_r}C_{l_1}'\quad\text{with }\#\mathcal E_r=M_1,
\]
i.e., $D_r$ are obtained as a union for $l_1\in\mathcal E_r\subset\{1,\ldots,L_1\}$ of $M_1$ consecutive $C_{l_1}'$-sets, hence they are also contained in disjoint intervals. Thus,
\[
vC_\varnothing'+xC_\varnothing'=\bigcup_{r,r'\leq L_1/M_1}(vD_r+xD_{r'}),
\]
and \paperref{result:4.1}{Lemma~\ref*{result:4.1}} implies the following bound,
\[
\sum_{r,r'}|vD_r+xD_{r'}|
=\int_{vC_\varnothing'+xC_\varnothing'}\sum_{r,r'}\chi_{vD_r+xD_{r'}}
\leq C\frac{L_1}{M_1}|vC_\varnothing'+xC_\varnothing'|
\]
\[
\Rightarrow\quad|vC_\varnothing'+xC_\varnothing'|>c\frac{M_1}{L_1}\sum_{r,r'}|vD_r+xD_{r'}|.\tag{4.5}\label{eq:4.5}
\]
Write further, for fixed $r,r'$,
\[
vD_r+xD_{r'}=\bigcup_{\substack{l_1\in\mathcal E_r\\l_1'\in\mathcal E_{r'}}}(vC_{l_1}'+xC_{l_1'}')
\]
and again by the multiplicity bound \paperref{eq:4.2}{(\ref*{eq:4.2})}, we get

\begin{equation}
 |vD_r+xD_{r'}|>\frac{c}{M_1}
 \sum_{\substack{l_1\in\mathcal E_r\\l_1'\in\mathcal E_{r'}}}
 |vC'_{l_1}+xC'_{l_1'}|.\tag{4.6}\label{eq:4.6}
\end{equation}
On the other hand, fixing points $\xi_l\in C_l$ so that
$|\xi_l-\xi_{l'}|>\frac12\mu_1$ for $l\ne l'$, it follows \paperref{eq:3.29}{(\ref*{eq:3.29})} that,
for $x\in[0,v]$,
\[
 vC'_{l_1}+xC'_{l_1'}\subset
 v\xi_{l_1}+x\xi_{l_1'}+[-2v\sigma_1\mu_1,2v\sigma_1\mu_1].
\]
Therefore, all sets
$\bigl(vC'_{l_1}+xC'_{l_1'}\bigr)_{l_1\in\mathcal E_r,\,l_1'\in\mathcal E_{r'}}$
are disjoint if
\begin{equation}
 \min_{\substack{k,l\in\mathcal E_r\\k',l'\in\mathcal E_{r'}\\(k,k')\ne(l,l')}}
 \bigl|(v\xi_k+x\xi_{k'})-(v\xi_l+x\xi_{l'})\bigr|
 >10v\sigma_1\mu_1.\tag{4.7}\label{eq:4.7}
\end{equation}
Note that if $k'=l'$, $|v\xi_k-v\xi_l|>\frac12v\mu_1>10v\sigma_1\mu_1$
since $\sigma_1<2^{-T\beta}\ll1$.
So we may assume $k'\ne l'$. In this case,
$|\xi_{k'}-\xi_{l'}|>\frac12\mu_1$, and the condition \paperref{eq:4.7}{(\ref*{eq:4.7})} is satisfied if
\begin{equation}
 \min_{\substack{k,l\in\mathcal E_r\\k',l'\in\mathcal E_{r'}\\k'\ne l'}}
 \left|x+v\frac{\xi_k-\xi_l}{\xi_{k'}-\xi_{l'}}\right|>20v\sigma_1.
 \tag{4.8}\label{eq:4.8}
\end{equation}

Note that the exceptional set
\[
 \{x\in[0,v]:\paperref{eq:4.8}{(\ref*{eq:4.8})}\text{ fails for }x\}
 \subset\bigcup_{\substack{k,l\in\mathcal E_r\\k',l'\in\mathcal E_{r'}\\k'\ne l'}}
 I\left(-v\frac{\xi_k-\xi_l}{\xi_{k'}-\xi_{l'}},20v\sigma_1\right).
\]
By \paperref{eq:4.1}{(\ref*{eq:4.1})}, \paperref{eq:4.4}{(\ref*{eq:4.4})} and the definition of $M_1$,
\[
 \begin{aligned}
 \mu(\{x\in[0,v]:\paperref{eq:4.8}{(\ref*{eq:4.8})}\text{ fails for }x\})
 &\lesssim M_1^4\left(\frac{v\sigma_1}{v}\right)^\tau\\
 &\lesssim\sigma_1^{\tau-4/D}<\sigma_1^{\tau/2},
 \end{aligned}
\]
provided $T$ is sufficiently large. The same estimate holds at every level.
If $x$ satisfies \paperref{eq:4.8}{(\ref*{eq:4.8})}, hence \paperref{eq:4.7}{(\ref*{eq:4.7})}, we have
\begin{equation}
 |vD_r+xD_{r'}|=
 \sum_{\substack{l_1\in\mathcal E_r\\l_1'\in\mathcal E_{r'}}}
 |vC'_{l_1}+xC'_{l_1'}|.\tag{4.9}\label{eq:4.9}
\end{equation}
From \paperref{eq:4.6}{(\ref*{eq:4.6})} and \paperref{eq:4.9}{(\ref*{eq:4.9})}, we may write
\[
 |vD_r+xD_{r'}|>
 \sum_{\substack{l_1\in\mathcal E_r\\l_1'\in\mathcal E_{r'}}}
 \varphi_{l_1,l_1'}(x)|vC'_{l_1}+xC'_{l_1'}|,
\]
where $\varphi_{l_1,l_1'}:[0,v]\to\{c/M_1,1\}$ and
$\mu[\varphi_{l_1,l_1'}=c/M_1]<\sigma_1^{\tau/2}$.
Recalling \paperref{eq:4.5}{(\ref*{eq:4.5})}, we see that

\begin{equation}
 |vC'_{\varnothing}+xC'_{\varnothing}|>
 \sum_{l_1,l_1'}\psi_{l_1,l_1'}(x)|vC'_{l_1}+xC'_{l_1'}|,
 \tag{4.10}\label{eq:4.10}
\end{equation}
where $\psi_{l_1,l_1'}(x)=c\frac{M_1}{L_1}\varphi_{l_1,l_1'}(x)$ satisfying
\[
\begin{aligned}
 \mathbb E_x[\log\psi_{l_1,l_1'}]
 &=\int\log\psi_{l_1,l_1'}(x)\,d\mu(x)\\
 &=\left(\log\frac{cM_1}{L_1}\right)
 \left(1-\mu[\varphi_{l_1,l_1'}=c/M_1]\right)
 +\left(\log\frac{c^2}{L_1}\right)\mu[\varphi_{l_1,l_1'}=c/M_1]\\
 &=\log\frac{cM_1}{L_1}
 +\underbrace{\left(\log\frac{c}{M_1}\right)}_{{<0}}
 \underbrace{\mu[\varphi_{l_1,l_1'}=c/M_1]}_{{<\sigma_1^{\tau/2}}}\\
 &\geq-\log L_1+\frac34\log M_1+\frac54\log c\\
 &\geq\log\left(\frac1{L_1}M_1^{1/2}\right).
\end{aligned}
\]
Here we take $\sigma_1^{\tau/2}<1/4$ and $M_1\geq c^{-5}$. In particular, we require
\begin{equation}
T>10/(\beta\tau),\tag{4.11}\label{eq:4.11}
\end{equation}
since that $\sigma_s<2^{-\beta T}$.

Repeat with the sets for each $l_1,l_1'$,
\[
 vC'_{l_1}+xC'_{l_1'}=
 \bigcup_{l_2,l_2'\le L_2}
 (vC'_{l_1,l_2}+xC'_{l_1',l_2'})
\]
to obtain
\[
 |vC'_{l_1}+xC'_{l_1'}|>
 \sum_{l_2,l_2'}\psi_{l_1,l_2;l_1',l_2'}(x)
 |vC'_{l_1,l_2}+xC'_{l_1',l_2'}|,
\]
where
\[
 \mathbb E_x[\log\psi_{l_1,l_2;l_1',l_2'}]
 >\log\left(\frac1{L_2}M_2^{1/2}\right).
\]
Iteration clearly provides the following
\[
\begin{aligned}
 |vC'_{\varnothing}+xC'_{\varnothing}|
 &>\sum_{\substack{l_1,\ldots,l_{t_2}\\l_1',\ldots,l_{t_2}'}}
 \psi_{l_1,l_1'}(x)\psi_{l_1,l_2;l_1',l_2'}(x)\cdots
 \psi_{l_1,\ldots,l_{t_2};l_1',\ldots,l_{t_2}'}(x)\\
 &\hspace{22mm}\times
 |vC'_{l_1,\ldots,l_{t_2}}+xC'_{l_1',\ldots,l_{t_2}'}|\\
 &\ge v\sum_{\bar l,\bar l'}\psi_{l_1,l_1'}(x)\cdots
 \psi_{\bar l;\bar l'}|C'_{l_1,\ldots,l_{t_2}}|,
\end{aligned}
\]
where, by Jensen's inequality,

\[
\begin{aligned}
 \mathbb E_x[\psi_{l_1,l_1'}\psi_{l_1,l_2;l_1',l_2'}\cdots]
 &\ge e^{\mathbb E_x[\log(\psi_{l_1,l_1'}\psi_{l_1,l_2;l_1',l_2'}\cdots)]}\\
 &=e^{\mathbb E_x[\log\psi_{l_1,l_1'}]
 +\mathbb E_x[\log\psi_{l_1,l_2;l_1',l_2'}]+\cdots}\\
 &>\prod_{s=1}^{t_2}\frac{M_s^{1/2}}{L_s}.
\end{aligned}
\]
Therefore, by \paperref{eq:3.30}{(\ref*{eq:3.30})} and \paperref{eq:3.31}{(\ref*{eq:3.31})},
\begin{equation}
\begin{aligned}
 \mathbb E_x[|vA'+xA'|]
 &\ge\mathbb E_x[|vC'_{\varnothing}+xC'_{\varnothing}|]\\
 &>v\sum_{\bar l,\bar l'}
 \mathbb E_x[\psi_{l_1,l_1'}\psi_{l_1,l_2;l_1',l_2'}\cdots]
 |C'_{l_1,\ldots,l_{t_2}}|\\
 &>v\left(\prod_{s=1}^{t_2}\frac{M_s^{1/2}}{L_s}\right)
 \left(\prod_{s=1}^{t_2}L_s\right)
 \sum_{l_1,\ldots,l_{t_2}}\delta|C_{l_1,\ldots,l_{t_2}}|\\
 &>v\delta\cdot M^{\gamma m_1/2}\cdot M^{-\gamma m_1/80}\cdot|A|\\
 &>vM^{\gamma m_1/4}|A'|\\
 &>v\delta^{-\frac{1-\sigma}{20}T^{-1/3}}|A'|.
\end{aligned}\tag{4.12}\label{eq:4.12}
\end{equation}

\section[Proof of Theorem 0.7]{Proof of \paperref{result:0.7}{Theorem 0.7}}\label{sec:5}

In this section, our first aim is to extract from $\mu$ in \paperref{result:0.7}{Theorem~\ref*{result:0.7}}
a measure $\mu'$ satisfying condition \paperref{eq:4.1}{(\ref*{eq:4.1})}. Recall $\mu$ satisfies
\[
 \mu(I)<\rho^\kappa
\]
whenever $I\subset[0,1]$ is a $\rho$-interval,
$\delta<\rho<\delta^{\varepsilon_0}$. To compare \paperref{eq:4.1}{(\ref*{eq:4.1})}:
\[
 \mu'(I)<\left(\frac\rho v\right)^\tau
\]
whenever $I\subset[0,v]$ is a $\rho$-interval, $\delta<\rho<v$,
here $\delta^{1/4}<v<1$ and $\tau$ is some constant to be determined later.
The restriction $\rho<\delta^{\varepsilon_0}$ nevertheless does not allow
one to invoke directly the conclusion from \paperref{sec:4}{Section~4} for $\mu$, and an
additional argument is needed.

Let $\mu_-\in\mathcal P(\mathbb R)$ be defined by $\mu_-(A)=\mu(-A)$ and set
\[
 \mu_1=\mu*\mu_-.
\]
Thus,
\[
 \widehat\mu_1=|\widehat\mu(\xi)|^2\ge0,
 \qquad
 \operatorname{supp}(\mu_1)\subset
 \operatorname{supp}(\mu)-\operatorname{supp}(\mu).
\]

Let $F_\Delta$ be the function
\begin{center}
\begin{tikzpicture}[x=1.1cm,y=2cm,baseline=(current bounding box.center)]
 \draw[->] (-2.45,0)--(2.5,0);
 \draw[thick] (-2,0)--(0,1)--(2,0);
 \draw[dashed] (0,0)--(0,1);
 \draw[dashed] (-1,0)--(-1,.5)--(1,.5)--(1,0);
 \foreach \x/\lab in {-2/{-2\Delta},-1/{-\Delta},0/{0},1/{\Delta},2/{2\Delta}}
 {\draw (\x,-.025)--(\x,.025);\node[below] at (\x,0) {$\lab$};}
 \fill (0,1) circle (1pt);
 \node[above right] at (0,1) {$1$};
 \node[right] at (1,.5) {$\frac12$};
 \node[right,align=left] at (2.8,.65)
 {$\chi_{I(0,\Delta)}\le2F_\Delta$\\$F_\Delta\le\chi_{I(0,2\Delta)}$};
\end{tikzpicture}
\end{center}
Consequently $0\leq\check F_\Delta(\xi)\leq2\check F_{\Delta/2}(\xi)$, and we have
\[
\begin{aligned}
 \mu_1(I(a,\Delta))
 &\le\int2F_\Delta(x-a)\,d\mu_1(x)\\
 &=2\int\widehat\mu_1(\xi)\check F_\Delta(\xi)e^{2\pi ia\xi}\,d\xi\\
 &\le4\int\widehat\mu_1(\xi)\check F_{\Delta/2}(\xi)\,d\xi\\
 &=4\int F_{\Delta/2}(x)\,d\mu_1(x)\\
 &\le4\mu_1(I(0,\Delta))\\
 &=4\iint\chi_{I(0,\Delta)}(x-y)\,d\mu(x)\,d\mu(y)
\end{aligned}
\]

\begin{equation}
\begin{aligned}
 &=4\int\mu(I(y,\Delta))\,d\mu(y)\\
 &\le4\sup_a\mu(I(a,\Delta)).
\end{aligned}\tag{5.1}\label{eq:5.1}
\end{equation}
Recalling $\mu$ satisfies a non-concentration assumption, we see that
\begin{equation}
 \mu_1(I(a,\Delta))<4\Delta^\kappa
 \quad\text{for }\delta<\Delta<\delta^{\varepsilon_0}.
 \tag{5.2}\label{eq:5.2}
\end{equation}
Let $\tau=\kappa/2$ and, in this selection argument, use the closed symmetric intervals $I(0,r)=[-r,r]$. Define
\begin{equation}
 v=\min\{r\in[\delta,1]:\mu_1([-r,r])\geq r^\tau\}.\tag{5.3}\label{eq:5.3}
\end{equation}
The set is nonempty because $r=1$ works, and the minimum exists by upper semicontinuity of $r\mapsto\mu_1([-r,r])$. The non-concentration bound gives $\mu_1([-\delta,\delta])\lesssim\delta^\kappa<\delta^\tau$, and excludes a neighborhood immediately above $\delta$, so $v>\delta$.
It follows from \paperref{eq:5.2}{(\ref*{eq:5.2})} that,
\[
 \text{if }v<\delta^{\varepsilon_0},\text{ then }
 v^\tau<4v^\kappa\ \Rightarrow\ v>\left(\frac1{16}\right)^{1/\kappa}.
\]
Thus,
\begin{equation}
 v\ge\min\left\{\delta^{\varepsilon_0},\left(\frac1{16}\right)^{1/\kappa}\right\}
 =\delta^{\varepsilon_0},\tag{5.4}\label{eq:5.4}
\end{equation}
whenever $\delta$ is small enough in terms of $\kappa$ and $\varepsilon_0$.
Next, we claim that
\[
 v\in\operatorname{supp}(\mu)-\operatorname{supp}(\mu).
\]
Indeed, consider
\[
 v'=\max\{|x|:x\in I(0,v)\cap\operatorname{supp}(\mu_1)\}\le v,
\]

then
\[
 v'\in\operatorname{supp}(\mu_1)\quad\text{or}\quad
 v'\in-\operatorname{supp}(\mu_1),
\]
in either case,
\[
 v'\in\operatorname{supp}(\mu)-\operatorname{supp}(\mu).
\]
Moreover,
\[
 \mu_1\bigl(\underbrace{I(0,v)\setminus I(0,v')}_{ {\text{which is disjoint with }\operatorname{supp}(\mu_1)}}\bigr)=0
\]
\[
\begin{aligned}
 &\Rightarrow\mu_1(I(0,v'))=\mu_1(I(0,v))\ge v^\tau\ge(v')^\tau\\
 &\Rightarrow v'=v\quad\text{by the minimality of }v.
\end{aligned}
\]
The case $v'\leq\delta$ in this argument is excluded by $\mu_1([-\delta,\delta])<\delta^\tau\leq v^\tau$. Since $\mu_1$ is symmetric, the positive point $v$ belongs to $\operatorname{supp}(\mu_1)$.

Define
\[
 \mu'=\frac{\mu_1|_{[0,v]}}{\mu_1([0,v])}.
\]
Symmetry gives $\mu_1([0,v])\geq\tfrac12\mu_1([-v,v])\geq\tfrac12v^\tau$. It follows that for $\delta<\Delta<v$,
\[
\begin{aligned}
 \mu'(I(a,\Delta))
 &\le\frac1{\mu_1([0,v])}\mu_1(I(a,\Delta))\\
 &\mathrel{\mathop\le^{\paperref{eq:5.1}{(\ref*{eq:5.1})}\,\&\,\paperref{eq:5.3}{(\ref*{eq:5.3})}}}
 \frac2{v^\tau}\cdot4\mu_1(I(0,\Delta))\\
 &<8\left(\frac\Delta v\right)^\tau\quad\text{by definition of }v.
\end{aligned}
\]

Moreover, \paperref{eq:5.1}{(\ref*{eq:5.1})} and symmetry give $\mu'(I(a,\Delta))\lesssim\mu'([0,\Delta])$, so both inequalities in \paperref{eq:4.1}{(\ref*{eq:4.1})} hold for $\mu'$, with absolute constants.

We now combine these estimates with the subcritical tree construction of \paperref{sec:2}{Sections~2}--\paperref{sec:3}{3}. Recall \paperref{eq:4.12}{(\ref*{eq:4.12})},
\[
 \int|vA'+xA'|\,d\mu'(x)
 >v\delta^{-\frac{1-\sigma}{20}T^{-1/3}}|A'|
 \mathrel{\mathop>^{\paperref{eq:5.4}{(\ref*{eq:5.4})}}}
 \delta^{\varepsilon_0-\frac{1-\sigma}{20}T^{-1/3}}|A'|.
\]
By pigeonholing principle, there exists $x\in\operatorname{supp}(\mu')$ such that
\begin{equation}
 |vA'+xA'|>
 \delta^{\varepsilon_0-\frac{1-\sigma}{20}T^{-1/3}}|A'|.\tag{5.5}\label{eq:5.5}
\end{equation}
Recall the parameters introduced above and their constraints:
\[
\begin{gathered}
 2^T<\eta^{-1/2},\qquad
 (1-\beta)(1-T^{-1/2})>\frac{2\sigma}{1+\sigma},\qquad
 T>10^{10}\left(\frac D\beta\right)^3,\\
 \varepsilon_0<\frac{1-\sigma}{8},\qquad
 T^{-1/3}<\frac{1-\sigma}{16}\kappa,\qquad
 D>\frac{20}\tau,\qquad \tau=\frac\kappa2,\qquad
 T>\frac{10}{\beta\tau},\qquad
 \rho(T)<\frac{1-\sigma}{5000}.
\end{gathered}
\]
Taking $0<\kappa\leq1$ without loss of generality, choose $D=C_0/\kappa$, where $C_0$ is a sufficiently large absolute constant, and then $T=C(\sigma)\kappa^{-3}$ with $C(\sigma)$ sufficiently large.
Assume assumption $\paperref{eq:star}{(*)}$ holds with $\varepsilon_0<c(\sigma)\kappa$.
The conclusion \paperref{eq:5.5}{(\ref*{eq:5.5})} then gives
\[
 |vA'+xA'|>\delta^{-c(\sigma)\kappa}|A'|
\]
for some $v,x\in\operatorname{supp}(\mu)-\operatorname{supp}(\mu)$.
It follows that there exist elements $x_1,x_2,x_3,x_4\in\operatorname{supp}(\mu)$ such that
\[
 N(x_1A+x_2A-x_3A-x_4A,\delta)>\delta^{-c(\sigma)\kappa}|A|.
\]
Using both the difference and sum versions of the \paperref{result:1.2}{Ruzsa triangle inequality},
\[
\begin{aligned}
 &N(x_1A+x_2A-x_3A-x_4A,\delta)\\
 &\qquad\le\frac{N(x_1A+x_2A+A,\delta)N(x_3A+x_4A+A,\delta)}{|A|},\\
 &N(x_1A+x_2A+A,\delta)
 \le\frac{N(x_1A+A+A,\delta)N(x_2A+A,\delta)}{|A|},\\
 &N(x_1A+A+A,\delta)
 \le\frac{N(x_1A+A,\delta)N(A+A+A,\delta)}{|A|}.
\end{aligned}
\]
By \paperref{result:1.3}{Pl\"unnecke--Ruzsa inequality},
\[
 N(A+A+A,\delta)<(\delta^{-\eta})^3|A|.
\]
Combining all these above inequalities, we obtain
\[
 \delta^{-(c(\sigma)\kappa-6\eta)}|A|^4
 <\prod_{i=1}^4N(A+x_iA,\delta).
\]

Hence, we may find some $x\in\operatorname{supp}(\mu)$ such that
\[
 N(A+xA,\delta)
 >\delta^{-(\frac{c(\sigma)}4\kappa-\frac32\eta)}|A|
 >\delta^{-\frac{c(\sigma)}5\kappa}|A|,
\]
here we assume $\eta<\frac1{30}c(\sigma)\kappa$.

In summary, we have proved that assume $\varepsilon_0<c(\sigma)\kappa$,
$\eta<\min\{\frac{c(\sigma)}{30}\kappa,
4^{-T}=4^{-C(\sigma)\kappa^{-3}}\}$ and $A$ has small
sumset:
$N(A+A,\delta)<\delta^{-\eta}|A|$, then we can find some
$x\in\operatorname{supp}(\mu)$ such that
\[
 N(A+xA,\delta)>\delta^{-\frac{c(\sigma)}5\kappa}|A|.
\]
On the other hand, suppose $N(A+A,\delta)>\delta^{-\eta}|A|$. The non-concentration of $\mu$ gives $\mu([0,\delta^{\varepsilon_0}])\lesssim\delta^{\varepsilon_0\kappa}<1$. Hence we may choose $x\in\operatorname{supp}(\mu)$ with $x\geq\delta^{\varepsilon_0}$. Since $A$ is $\delta$-separated, each $\delta$-interval contains at most $O(1/x)$ points of $xA$, and therefore
\[
N(xA,\delta)\gtrsim x|A|\geq\delta^{\varepsilon_0}|A|.
\]
The sum triangle inequality now yields
\[
\delta^{-\eta}|A|
\lesssim\frac{N(A+xA,\delta)^2}{N(xA,\delta)}
\lesssim\delta^{-\varepsilon_0}\frac{N(A+xA,\delta)^2}{|A|},
\]
so
\[
N(A+xA,\delta)\gtrsim\delta^{-(\eta-\varepsilon_0)/2}|A|
>\delta^{-\eta/4}|A|
\]
for sufficiently small $\delta$, provided $\varepsilon_0<\eta/4$.

The parameters can be chosen in the following order: first fix $\beta,\tau,D,T$ as above, next take
\[
0<\eta<\min\{4^{-T},c(\sigma)\kappa/30\},
\]
and finally choose $\varepsilon_0>0$ smaller than all the preceding upper bounds and $\eta/4$. Both cases then give expansion with exponent $\varepsilon_1=\eta/4$. This proves \paperref{result:0.7}{Theorem~\ref*{result:0.7}}.

\section{Proof of the discretized sum-product theorem}\label{sec:6}

Let $A$ be a $\delta$-separated set, $A\subset[1,2]$, $|A|=\delta^{-\sigma}$
and satisfy the non-concentration property:
\[
 |A\cap I(a,\rho)|<\rho^\kappa|A|
 \quad\text{for any }a\in\mathbb R,\quad\delta<\rho<\delta^{\varepsilon_0}.
\]
Apply \paperref{result:0.7}{Theorem~\ref*{result:0.7}} with measure $\mu\in\mathcal P([0,1])$ given by
\[
 \mu=\frac1{|A|}\sum_{x\in A}\delta_{x/a},
\]
where $a=\max A$ (so that $a^{-1}A\subset[0,1]$).
It follows from \paperref{result:0.7}{Theorem~\ref*{result:0.7}} that for some $\varepsilon_2=\varepsilon_2(\sigma,\kappa)>0$
and $\widetilde x\in\operatorname{supp}(\mu)=a^{-1}A$, we have
\[
 N(A+\widetilde xA,\delta)>\delta^{-\varepsilon_2}|A|
 \ \Rightarrow\ N(aA+xA,\delta)>\delta^{-\varepsilon_2}|A|,
\]
for some $x=a\widetilde x\in A$. Therefore,
\begin{equation}
 N(A\cdot A+A\cdot A,\delta)>\delta^{-\varepsilon_2}|A|.\tag{6.1}\label{eq:6.1}
\end{equation}
Statement \paperref{eq:6.1}{(\ref*{eq:6.1})} doesn't immediately imply \paperref{eq:0.2}{(\ref*{eq:0.2})} and we rely on further
construction. We argue by contradiction: suppose that \paperref{eq:0.2}{(\ref*{eq:0.2})} fails, so
\begin{equation}
 N(A+A,\delta)<\delta^{-\varepsilon}|A|
 \quad\&\quad N(A\cdot A,\delta)<\delta^{-\varepsilon}|A|,\tag{6.2}\label{eq:6.2}
\end{equation}

where $\varepsilon>0$ is given and small. Our aim is to construct a large
subset $A_1$ of $A$ such that $N(A_1\cdot A_1+A_1\cdot A_1)$ is also well
behaved. The next result provides a stronger statement of independent interest.

\begin{proposition}\label{result:6.1}
Assume $A$ satisfies \paperref{eq:6.2}{(\ref*{eq:6.2})}. Then there exists a subset $A_1\subset A$ such that
\begin{equation}|A_1|>\delta^{2\varepsilon}|A|\tag{6.3}\label{eq:6.3}\end{equation}
and, for any given $k\in\mathbb Z_+$,
\begin{equation}
 N(kA_1^k,\delta)<\delta^{-c(k)\varepsilon}|A_1|,\tag{6.4}\label{eq:6.4}
\end{equation}
where $kA_1^k$ denotes the $k$-fold sumset of the $k$-fold product-set $A_1^k$ of $A_1$.
\end{proposition}

\begin{proof}
Taking $C=-\{x\in A:\operatorname{dist}(x,B)<\delta\}$ and $I=N(-C,\delta)$,
apply \paperref{result:1.2}{Ruzsa's triangle inequality} to $A,B,C$, we obtain
\begin{equation}
\begin{aligned}
 N(A-B,\delta)
 &\le\frac{N(A-C,\delta)N(B-C,\delta)}{N(C,\delta)}\\
 &\le\frac{N(A+A,\delta)N(B+B,\delta)}I.
\end{aligned}\tag{6.5}\label{eq:6.5}
\end{equation}

We first exploit the multiplicative structure of $A$.
We consider the multiplicative energy of $A$, by Cauchy--Schwarz inequality,
\begin{equation}
\begin{aligned}
 &|\{(x_1,x_2,x_3,x_4)\in A^4:|x_1x_2-x_3x_4|<\delta\}|\\
 &=\sum_{I\in\mathcal D_\delta(A\cdot A)}
 |\{(x_1,x_2,x_3,x_4)\in A^4:x_1x_2,x_3x_4\in I\}|\\
 &=\sum_{I\in\mathcal D_\delta(A\cdot A)}
 |\{(x_1,x_2)\in A^2:x_1x_2\in I\}|^2\\
 &\ge N(A\cdot A,\delta)^{-1}
 \left(\sum_{I\in\mathcal D_\delta(A\cdot A)}
 |\{(x_1,x_2)\in A^2:x_1x_2\in I\}|\right)^2\\
 &=N(A\cdot A,\delta)^{-1}|A^2|^2\\
 &\mathrel{\mathop>^{\paperref{eq:6.2}{(\ref*{eq:6.2})}}}\delta^\varepsilon|A|^3.
\end{aligned}\tag{6.6}\label{eq:6.6}
\end{equation}
Hence, by pigeonholing principle, we may fix an element $a\in A$ and a
subset $A_1\subset A$ such that
\begin{itemize}
 \item $|A_1|>\delta^{2\varepsilon}|A|$;
 \item $|\{(x,y)\in A^2:|x_1x-ay|<\delta\}|>\delta^{2\varepsilon}|A|$
 for all $x_1\in A_1$.
\end{itemize}

\par\noindent\rule{\linewidth}{0.4pt}
Indeed, for any $a\in A$, set
\[
 A_1(a)=\{x_1\in A:|\{(x,y)\in A^2:|x_1x-ay|<\delta\}|>
 \delta^{2\varepsilon}|A|\}.
\]
Notice that for any $x_1,a\in A$, we have the following
\[
\begin{aligned}
 &|\{(x,y)\in A^2:|x_1x-ay|<\delta\}|\\
 &=\sum_{x\in A}\left|\left\{y\in A:
 \left|y-\frac{x_1x}a\right|<\frac\delta a\right\}\right|\\
 &\le\sum_{x\in A}\left|\left\{y\in A:
 \left|y-\frac{x_1x}a\right|<\delta\right\}\right|\\
 &\lesssim N(A,\delta)=|A|,
\end{aligned}
\]
since $a\in[1,2]$ and the points in $A$ are $\delta$-separated.
Then, by \paperref{eq:6.6}{(\ref*{eq:6.6})},
\[
\begin{aligned}
 \delta^\varepsilon|A|^3
 &<|\{(x_1,x,a,y)\in A^4:|x_1x-ay|<\delta\}|\\
 &=\left(\sum_{a\in A}\sum_{x_1\in A_1(a)}
 +\sum_{a\in A}\sum_{x_1\notin A_1(a)}\right)
 |\{(x,y)\in A^2:|x_1x-ay|<\delta\}|\\
 &<\sum_{a\in A}|A_1(a)|\cdot C|A|
 +|A|\cdot|A|\cdot\delta^{2\varepsilon}|A|
\end{aligned}
\]
\[
 \Rightarrow\delta^{2\varepsilon}|A|^2<\sum_{a\in A}|A_1(a)|.
\]

Hence, there exists some $a\in A$ such that $|A_1(a)|>\delta^{2\varepsilon}|A|$,
and take $A_1=A_1(a)$, as desired.
\par\noindent\rule{\linewidth}{0.4pt}

For $x_1\in A_1$, we will apply \paperref{eq:6.5}{(\ref*{eq:6.5})} to $x_1A$ and $aA$, while
\[
\begin{aligned}
 I&=N(\{\widetilde x\in x_1A:
 \operatorname{dist}(\widetilde x,aA)<\delta\},\delta)\\
 &=|\{x\in A:|x_1x-ay|<\delta\text{ for some }y\in A\}|\\
 &\sim\sum_{y\in A}|\{x\in A:|x_1x-ay|<\delta\}|
 \quad\text{(since these points are $\delta$-separated)}\\
 &=|\{(x,y)\in A^2:|x_1x-ay|<\delta\}|\\
 &>\delta^{2\varepsilon}|A|.
\end{aligned}
\]
It follows from \paperref{eq:6.2}{(\ref*{eq:6.2})} and \paperref{eq:6.5}{(\ref*{eq:6.5})} that
\[
\begin{aligned}
 N(x_1A-aA,\delta)
 &\le\frac{N(x_1A+x_1A,\delta)N(aA+aA,\delta)}I\\
 &\lesssim\frac{N(A+A,\delta)^2}{\delta^{2\varepsilon}|A|}\\
 &\lesssim\delta^{-4\varepsilon}|A|.
\end{aligned}
\]
Application of \paperref{result:1.2}{Ruzsa's triangle inequality} gives
\[
 N(x_1A-x_2A,\delta)
 \le\frac{N(x_1A-aA,\delta)N(x_2A-aA,\delta)}{N(aA,\delta)}
\]

\begin{equation}
 \lesssim\delta^{-8\varepsilon}|A|\quad\text{if }x_1,x_2\in A_1.
 \tag{6.7}\label{eq:6.7}
\end{equation}
Next, let $x_1,x_2,x_3,x_4\in A_1$, it follows from \paperref{eq:6.7}{(\ref*{eq:6.7})} that
\[
\begin{aligned}
 N(x_1x_2A-x_3x_2A,\delta)&\lesssim\delta^{-8\varepsilon}|A|,\\
 N(x_3x_4A-x_3x_2A,\delta)&\lesssim\delta^{-8\varepsilon}|A|,
\end{aligned}
\]
and by \paperref{result:1.2}{Ruzsa's triangle inequality} again, that
\[
 N(x_1x_2A-x_3x_4A,\delta)\lesssim\delta^{-16\varepsilon}|A|.
\]
The continuation of the process is clear. If $x_1,\ldots,x_{2k}\in A_1$,
we have
\[
\begin{aligned}
 &N(x_1\cdots x_kA_1-x_{k+1}\cdots x_{2k}A_1,\delta)\\
 &\qquad\le N(x_1\cdots x_kA-x_{k+1}\cdots x_{2k}A,\delta)
 <\delta^{-C(k)\varepsilon}|A|.
\end{aligned}
\]
Hence, if $x_1,x_2\in A_1$ and $y_1,y_2\in A_1^k$, we have also
\begin{equation}
\begin{aligned}
 N(y_1x_1^{-1}A_1-y_2x_2^{-1}A_1,\delta)
 &\sim N(x_2y_1A_1-x_1y_2A_1,\delta)\\
 &<\delta^{-C(k+1)\varepsilon}|A|\\
 &<\delta^{-C(k)\varepsilon}|A_1|.
\end{aligned}\tag{6.8}\label{eq:6.8}
\end{equation}

We now estimate $N(A_1^k-A_1^k,\delta)$. From \paperref{eq:6.2}{(\ref*{eq:6.2})}, \paperref{result:1.2}{Ruzsa's triangle
inequality} and the \paperref{result:1.3}{Pl\"unnecke--Ruzsa inequality} (applied in multiplicative
form), we have
\[
\begin{aligned}
 N(A_1^kA_1^{-1},\delta)
 &\le N(A^kA^{-1},\delta)\\
 &\le\frac{N(A^k\cdot A,\delta)N(A\cdot A,\delta)}{N(A^{-1},\delta)}\\
 &\le\frac{(\delta^{-\varepsilon})^{k+1}|A|
 \cdot\delta^{-\varepsilon}|A|}{|A|}\\
 &=\delta^{-(k+2)\varepsilon}|A|\\
 &<\delta^{-(k+4)\varepsilon}|A_1|.
\end{aligned}
\]
Equivalently, there exists a $\delta$-separated subset
$\mathcal E\subset A_1^kA_1^{-1}$ satisfying
\begin{equation}
 |\mathcal E|<\delta^{-(k+4)\varepsilon}|A_1|
 \quad\&\quad
 A_1^kA_1^{-1}\subset\bigcup_{z\in\mathcal E}I(z,\delta).
 \tag{6.9}\label{eq:6.9}
\end{equation}
Fix now $y_1,y_2\in A_1^k$.

For arbitrary chosen $x_1,x_2\in A_1$, take $z_1,z_2\in\mathcal E$ such that

\[
\begin{aligned}
 &y_ix_i^{-1}\in I(z_i,\delta)
 \ \Leftrightarrow\ |y_ix_i^{-1}-z_i|<\delta\\
 &\Rightarrow |(y_1-y_2)-(x_1z_1-x_2z_2)|<4\delta\\
 &\Rightarrow y_1-y_2\in z_1A_1-z_2A_1+I(0,4\delta).
\end{aligned}
\]
Clearly, we can assume the union in \paperref{eq:6.9}{(\ref*{eq:6.9})} has bounded overlapping.
For fixed $y_i\in A_1^k\subset[1,2^k]$, the map $x_i\mapsto y_i/x_i$ is bi-Lipschitz on $[1,2]$, with constants depending only on $k$. Hence a chosen $z_i$ can arise from only $O_k(1)$ points of the $\delta$-separated set $A_1$. Thus, if we vary $x_1,x_2\in A_1$, the pair $(z_1,z_2)$ ranges over $\gtrsim_k|A_1|^2$ distinct elements of $\mathcal E\times\mathcal E$.
Therefore, we have
\[
 \sum_{z_1,z_2\in\mathcal E}
 \chi_{z_1A_1-z_2A_1+I(0,4\delta)}
 \gtrsim|A_1|^2\chi_{A_1^k-A_1^k}.
\]
It follows that
\[
\begin{aligned}
 N(A_1^k-A_1^k,\delta)
 &=\delta^{-1}|A_1^k-A_1^k|
 =\delta^{-1}\int\chi_{A_1^k-A_1^k}\\
 &\lesssim|A_1|^{-2}\delta^{-1}
 \int\sum_{z_1,z_2\in\mathcal E}
 \chi_{z_1A_1-z_2A_1+I(0,4\delta)}\\
 &=|A_1|^{-2}\sum_{z_1,z_2\in\mathcal E}
 \delta^{-1}|z_1A_1-z_2A_1+I(0,4\delta)|\\
 &\sim|A_1|^{-2}\sum_{z_1,z_2\in\mathcal E}N(z_1A_1-z_2A_1,\delta)\\
 &\mathrel{\mathop\lesssim^{\paperref{eq:6.8}{(\ref*{eq:6.8})}}}
 |A_1|^{-2}|\mathcal E|^2\delta^{-C(k)\varepsilon}|A_1|\\
 &\mathrel{\mathop<^{\paperref{eq:6.9}{(\ref*{eq:6.9})}}}\delta^{-C(k)\varepsilon}|A_1|.
\end{aligned}
\]
Finally, by \paperref{result:1.3}{Pl\"unnecke--Ruzsa inequality}, we get
\[
 N(kA_1^k,\delta)<(\delta^{-C(k)\varepsilon})^k|A_1|
 =\delta^{-c(k)\varepsilon}|A_1|,
\]
This proves \paperref{result:6.1}{Proposition~\ref*{result:6.1}}.
\end{proof}

Now we continue to prove \paperref{result:0.6}{Theorem~\ref*{result:0.6}}. Application of \paperref{result:6.1}{Proposition~\ref*{result:6.1}} with
$\varepsilon=\varepsilon_1$ and $k=2$ gives a set $A_1\subset A$ satisfying
\begin{itemize}
 \item $|A_1|>\delta^{2\varepsilon_1}|A|
 =\delta^{-(\sigma-2\varepsilon_1)}$;
 \item $N(A_1\cdot A_1+A_1\cdot A_1,\delta)
 <\delta^{-C\varepsilon_1}|A_1|$.
\end{itemize}
Hence,
\[
 |A_1|=\delta^{-\sigma_1}\quad\text{with }\sigma-2\varepsilon_1<\sigma_1\le\sigma
\]
and if $I$ is a $\rho$-interval, $\delta<\rho<\delta^{\varepsilon_0}$,
we still have
\[
\begin{aligned}
 |A_1\cap I|&\le|A\cap I|<\rho^\kappa|A|
 <\rho^\kappa\delta^{-2\varepsilon_1}|A_1|\\
 &<\rho^{\kappa-2\varepsilon_1/\varepsilon_0}|A_1|
 <\rho^{\kappa/2}|A_1|
\end{aligned}
\]
provided

\[
\varepsilon_1<\frac14\varepsilon_0\kappa.
\]
Now apply \paperref{result:0.7}{Theorem~\ref*{result:0.7}}, hence \paperref{eq:6.1}{(\ref*{eq:6.1})}, to the set $A_1$ with $\sigma$ replaced with $\sigma_1$ and $\kappa$ by $\kappa/2$, we obtain
\[
N(A_1\cdot A_1+A_1\cdot A_1,\delta)>\delta^{-\varepsilon_2}|A_1|.
\]
But this is contradictory if $\varepsilon_1<\min\{\frac14\varepsilon_0\kappa,\frac1C\varepsilon_2\}$. And we complete our proof of \paperref{result:0.6}{Theorem~\ref*{result:0.6}}.

\section{Projection theorems}\label{sec:7}

\subsection{Discretized projection theorem in plane}\label{sec:7.1}
We start with the discretized projection theorem, \paperref{result:0.8}{Theorem~\ref*{result:0.8}}. It is essentially another expression of \paperref{result:0.7}{Theorem~\ref*{result:0.7}}: If $\mathcal A\subset[1,2]\times[1,2]$ has the form $\mathcal A=A\times A$, where $A\subset[1,2]$ satisfies
\[
N(A,\delta)=\delta^{-\alpha/2}
\quad\&\quad
N(A\cap I,\delta)<\rho^{\beta/2}N(A,\delta)
\]
for any $\rho$-intervals $I$, $\delta<\rho<\delta^{\tau_0}$, then \paperref{result:0.7}{Theorem~\ref*{result:0.7}} implies that, there exists $x\in\operatorname{supp}(\mu_1)$ such that
\[
N(A+xA,\delta)>\delta^{-\varepsilon_1}N(A,\delta)
=\delta^{-\varepsilon_1-\alpha/2}.
\]
Choose $\theta=(1,x)\in\operatorname{supp}(\mu_1)$, notice that
\begin{align*}
\pi_\theta(\mathcal A)
&=\pi_{(1,x)}(A\times A)
=\{(1,x)\cdot(a,b):(a,b)\in A\times A\}\\
&=\{a+xb:a,b\in A\}\\
&=A+xA,
\end{align*}
we obtain
\[
|\pi_\theta(\mathcal A)|=\delta N(\pi_\theta(\mathcal A),\delta)
>\delta\cdot\delta^{-\varepsilon_1-\alpha/2}
=\delta^{1-(\alpha/2+\varepsilon_1)},
\]
as desired.

However, the set $\mathcal A$ couldn't be represented by this form in general. To handle this problem, we rely on the \paperref{result:1.6}{Balog--Szemer\'edi--Gowers theorem}.

\begin{proof}[{Proof of \paperref{result:0.8}{Theorem~\ref*{result:0.8}}}]
Choose $\eta>\alpha/2$ close to $\alpha/2$ (to be specified) and assume
\[
|\pi_\theta(\mathcal A)|\leq\delta^{1-\eta}
\quad\text{for all }\theta\in\operatorname{supp}(\mu_1).
\]
From our assumption on $\mu_1$, we may in particular specify $\theta_1,\theta_2,\theta_3\in S^1$ such that
\[
|\theta_i-\theta_j|>\delta^{\tau_0}\quad\text{for }i\ne j
\]
\[
\&\quad |\pi_{\theta_i}(\mathcal A)|\leq\delta^{1-\eta},\quad i=1,2,3.
\]
After a linear coordinate change, we let $\theta_1=0$, $\theta_2=\pi/2$ and $\theta_3=\pi/4$.

Let
\[
A_i=\pi_{\theta_i}(\mathcal A),\quad i=1,2,3
\]
then by assumption $N(A_i,\delta)\leq\delta^{-\eta}$. And notice that
\begin{align*}
\mathcal A&\subset\pi_{\theta_i}(\mathcal A)\times\pi_{\theta_j}(\mathcal A)
=A_i\times A_j\quad\text{for }i\ne j\\
\Rightarrow\quad\delta^{-\alpha}
&=N(\mathcal A,\delta)\leq N(A_i,\delta)N(A_j,\delta)
\leq\delta^{-\eta}N(A_i,\delta)\\
\Rightarrow\quad\delta^{\eta-\alpha}
&\leq N(A_i,\delta)\leq\delta^{-\eta}\quad\text{for }i=1,2,3.
\end{align*}
Note that $\mathcal A\subset A_1\times A_2$ and

\[
N(\mathcal A,\delta)=\delta^{-\alpha}
\geq\delta^{2\eta-\alpha}N(A_1,\delta)N(A_2,\delta),
\]
we consider the partial sumset
\begin{align*}
A_1\overset{\mathcal A}{+}A_2
&=\{x+y:(x,y)\in\mathcal A\}
=\{(1,1)\cdot(x,y)\in\mathcal A\}\\
&=\pi_{\theta_3}(\mathcal A)=A_3
\end{align*}
and
\begin{align*}
N(A_1\overset{\mathcal A}{+}A_2,\delta)
&=N(A_3,\delta)<\delta^{-\eta}\\
&\leq\delta^{-(2\eta-\alpha)}N(A_1,\delta)^{1/2}N(A_2,\delta)^{1/2}.
\end{align*}
Hence, applying the \paperref{result:1.6}{Balog--Szemer\'edi--Gowers theorem} gives subsets $B_1\subset A_1$, $B_2\subset A_2$ with the following properties:
\begin{equation}
N(B_i,\delta)>\delta^{-\gamma_1}\quad(i=1,2),\qquad
\gamma_1=\alpha-\eta-C(2\eta-\alpha)\nearrow\frac\alpha2.
\tag{7.1}\label{eq:7.1}
\end{equation}
\begin{equation}
N(B_1+B_2,\delta)<\delta^{-\gamma_2},\qquad
\gamma_2=\eta+C(2\eta-\alpha)\searrow\frac\alpha2.
\tag{7.2}\label{eq:7.2}
\end{equation}
\begin{equation}
|\mathcal A\cap(B_1\times B_2)|>\delta^{2-\alpha_1},\qquad
\alpha_1=2(\alpha-\eta)-C(2\eta-\alpha)\nearrow\alpha.
\tag{7.3}\label{eq:7.3}
\end{equation}
where
\begin{equation}
\left|\gamma_1-\frac\alpha2\right|
+\left|\gamma_2-\frac\alpha2\right|+|\alpha_1-\alpha|
<C\left(\left|\eta-\frac\alpha2\right|+\tau_0\right).
\tag{7.4}\label{eq:7.4}
\end{equation}
Now, roughly speaking, we could view $\mathcal A\approx B_1\times B_2$, and we consider $A\approx B_1$. We will apply \paperref{result:0.7}{Theorem~\ref*{result:0.7}} to $B_1$ to obtain some $x=\tan\theta$ with

$|B_1+xB_1|\gg|B_1|$. Even though
\[
\pi_\theta(\mathcal A)\approx B_1+xB_2\ne B_1+xB_1,
\]
we can still to show
\[
|\pi_\theta(\mathcal A)|\gtrsim|B_1+xB_1|
\quad(\gg|B_1|\approx\delta^{1-\alpha/2})
\]
by some sumset inequality and \paperref{eq:7.2}{(\ref*{eq:7.2})}: $|B_1+B_2|\lesssim|B_1|$. Next, let's make this process more rigorous.

Define
\[
A=\left\{x_1\in B_1:|B_2\cap\mathcal A(x_1)|>\frac12\delta^{1+\eta-\alpha_1}\right\},
\]
where $\mathcal A(x_1)=\{x_2:(x_1,x_2)\in\mathcal A\}$ is the $x_1$-section. By \paperref{eq:7.3}{(\ref*{eq:7.3})},
\begin{align*}
\delta^{2-\alpha_1}
&<|\mathcal A\cap(B_1\times B_2)|\\
&=|\mathcal A\cap(A\times B_2)|
+\underbrace{|\mathcal A\cap((B_1\setminus A)\times B_2)|}_{\displaystyle{=\int_{B_1\setminus A}|B_2\cap\mathcal A(x_1)|\,dx_1}}\\
&\leq |A|\cdot|B_2|+|A_1|\cdot\frac12\delta^{1+\eta-\alpha_1}\\
&\leq\delta^{1-\eta}|A|+\frac12\delta^{2-\alpha_1}
\end{align*}

\begin{align*}
\Rightarrow\quad |A|&>\frac12\delta^{1+\eta-\alpha_1}\\
\Rightarrow\quad N(A,\delta)&=\delta^{-\sigma}
\quad\text{with }\eta>\sigma>\alpha_1-\eta
\end{align*}
where,
\begin{equation}
\left|\sigma-\frac\alpha2\right|<C\left(\left|\eta-\frac\alpha2\right|+\tau_0\right).
\tag{7.5}\label{eq:7.5}
\end{equation}
Note that, since $0<\alpha<2$, \paperref{eq:7.5}{(\ref*{eq:7.5})} in particular ensures that
\[
\min\{\sigma,1-\sigma\}>\frac13\min\{\alpha,2-\alpha\}
\quad\Rightarrow\quad\sigma\in(0,1).
\]
Consider for the image measure of $\mu_1$ under the map $\theta\mapsto\tan\theta$, it is easy to check
\[
(\tan\mu_1)(I(a,\rho))<\rho^\kappa
\quad\text{for any }a\in\mathbb R\text{ and }\delta<\rho<\delta^{\tau_0}.
\]
In order to apply \paperref{result:0.7}{Theorem~\ref*{result:0.7}}, it suffices to prove $A$ satisfying the non-concentration property $\paperref{eq:star}{(*)}$. This is true for an appropriate choice of $\theta_1$ using the non-concentration property for $\mathcal A$. We omit this technical selection step, as in a proof sketch, and refer to Bourgain's argument.

Once \paperref{result:0.7}{Theorem~\ref*{result:0.7}} is applicable, we obtain some element $x=\tan\theta$ for $\theta\in\operatorname{supp}(\mu_1)$, satisfying \paperref{eq:0.3}{(\ref*{eq:0.3})}:
\[
N(A+xA,\delta)>\delta^{-\varepsilon_1}N(A,\delta)
\]

\begin{equation}
\Longleftrightarrow\quad |A+xA|>\delta^{-\varepsilon_1}|A|
\tag{7.6}\label{eq:7.6}
\end{equation}
with $\varepsilon_1=\varepsilon_1(\sigma,\kappa,\beta)>0$.

We now proceed as follows. Let
\[
\mathcal B=\mathcal A\cap(B_1\times B_2).
\]
By \paperref{eq:7.3}{(\ref*{eq:7.3})}, $\mathcal B$ satisfies
\begin{equation}
|\mathcal B|>\delta^{2-\alpha_1}.
\tag{7.7}\label{eq:7.7}
\end{equation}
Notice that for any set $E,F\subset\mathbb R^2$, any $\eta\in E$,
\begin{align*}
\chi_E(\eta)&=\frac1{|F-\eta|}\int\chi_{F-\eta}(\xi)\,d\xi\\
&\leq\frac1{|F|}\int_{F-E}\chi_{F-\xi}(\eta)\,d\xi
\end{align*}
\begin{equation}
\Rightarrow\quad\chi_E\leq\frac1{|F|}\int_{F-E}\chi_{F-\xi}\,d\xi.
\tag{7.8}\label{eq:7.8}
\end{equation}
(When $\eta\notin E$ this inequality is trivial.) Let $\pi(\xi)=\xi_1+x\xi_2$ for $\xi\in\mathbb R^2$, then $\pi(\xi)\approx\pi_\theta(\xi)$ up to a constant, since $x=\tan\theta$, and $|\pi(\mathcal B)|\leq|\pi_\theta(\mathcal A)|$. It clearly follows from \paperref{eq:7.8}{(\ref*{eq:7.8})} that
\[
\chi_{A+xA}=\chi_{A\times A}\circ\pi^{-1}
\leq\frac1{|\mathcal B|}\int_{\mathcal B-(A\times A)}
(\chi_{\mathcal B-\xi}\circ\pi^{-1})\,d\xi
\]

\[
\leq\frac1{|\mathcal B|}\int_{(B_1-A)\times(B_2-A)}
\chi_{\pi(\mathcal B)-\pi(\xi)}\,d\xi.
\]
Recall \paperref{eq:7.6}{(\ref*{eq:7.6})} and \paperref{eq:7.7}{(\ref*{eq:7.7})}, we get
\begin{align*}
\delta^{-\varepsilon_1}|A|&<|A+xA|\\
&\leq\frac1{|\mathcal B|}\int_{(B_1-A)\times(B_2-A)}
\int\chi_{\pi(\mathcal B)-\pi(\xi)}(\eta)\,d\eta\,d\xi\\
&<\delta^{-2+\alpha_1}|B_1-A|\cdot|B_2-A|\cdot|\pi(\mathcal B)|\\
&\leq\delta^{-2+\alpha_1}|B_1-B_1|\cdot|B_1-B_2|\cdot|\pi_\theta(\mathcal A)|.
\end{align*}
By \paperref{eq:7.1}{(\ref*{eq:7.1})}, \paperref{eq:7.2}{(\ref*{eq:7.2})}, \paperref{eq:7.4}{(\ref*{eq:7.4})} and \paperref{result:1.2}{Ruzsa's triangle inequality},
\[
|B_1-B_1|\leq\frac{|B_1+B_2|\cdot|B_1+B_2|}{|B_2|}
\]
\[
|B_1-B_2|\leq\frac{|B_1+B_1|\cdot|B_1+B_2|}{|B_1|}
\leq\frac{|B_1+B_2|\cdot|B_1+B_2|}{|B_2|}\cdot\frac{|B_1+B_2|}{|B_1|}
\]
\[
\Rightarrow\quad |B_1-B_1|\cdot|B_1-B_2|
\lesssim\frac{|B_1+B_2|^5}{|B_1||B_2|^2}<\delta^{2+3\gamma_1-5\gamma_2}
\]
\[
\Rightarrow\quad\delta^{-\varepsilon_1}|A|
<\delta^{-2+\alpha_1}\delta^{2(1-\gamma)}|\pi_\theta(\mathcal A)|
\]

with
\[
\left|\gamma-\frac\alpha2\right|<C\left|\eta-\frac\alpha2\right|.
\]
Recall $|A|=\delta^{1-\sigma}$ and \paperref{eq:7.5}{(\ref*{eq:7.5})}, it follows that
\[
|\pi_\theta(\mathcal A)|>\delta^{1-(\alpha_1-2\gamma+\sigma)-\varepsilon_1}.
\]
Since
\[
\left|\alpha_1-2\gamma+\sigma-\frac\alpha2\right|
<C\left|\eta-\frac\alpha2\right|
\]
we get, for $|\eta-\alpha/2|$ sufficiently small,
\[
|\pi_\theta(\mathcal A)|>\delta^{1-(\alpha/2+\varepsilon_1/2)}>\delta^{1-\eta}
\]
and hence a contradiction.
\end{proof}

\begin{remark} As we mentioned about \paperref{result:0.7}{Theorem~\ref*{result:0.7}}, \paperref{eq:0.3}{(\ref*{eq:0.3})}, the conclusion \paperref{eq:0.4}{(\ref*{eq:0.4})} in \paperref{result:0.8}{Theorem~\ref*{result:0.8}} also holds for all $\theta\in D\subset\operatorname{supp}(\mu_1)$, where
\[
\mu_1(D)>1-\delta^{\frac1{10}\kappa\tau_0}.
\]
Moreover, for $\theta\in D$, we have not only
\[
|\pi_\theta(\mathcal A)|>\delta^{1-\eta},
\]
but also
\[
|\pi_\theta(\mathcal A')|>\delta^{1-\eta}
\]

whenever $\mathcal A'\subset\mathcal B$ satisfies
\[
|\mathcal A'|>\delta^\varepsilon|\mathcal B|,
\]
for some $\varepsilon=\varepsilon(\alpha,\beta,\kappa)>0$.

For this fixed structured set $\mathcal B$, the subset estimate follows by using $F=\mathcal A'$ in the averaging argument, with the same $A$ and $B_i$.
\end{remark}

\subsection{Dimension of projection set in plane}\label{sec:7.2}
Now we sketch the proof of \paperref{result:0.9}{Theorem~\ref*{result:0.9}}, the projection theorem in terms of the Hausdorff dimension.

Set
\[
E=\{\theta\in S^1:\dim_H\pi_\theta(\mathcal A)<\eta\},
\]
assume
\[
\dim_H E>\kappa.
\]
By Frostman's lemma, there is a measure $\mu_1$ on $E$ satisfying
\[
\mu_1(I(a,\rho))<\rho^\kappa
\quad\text{for any }a\in\mathbb R\text{ and }\rho>0.
\]
And then
\[
\dim_H\pi_\theta(\mathcal A)<\eta
\quad\text{for all }\theta\in\operatorname{supp}(\mu_1).
\]
Since $\dim_H\mathcal A>\alpha$, apply Frostman's lemma again produce a measure $\mathbb P$ on $\mathcal A$ satisfying

\[
\mathbb P(B(x,\rho))<\rho^\alpha
\quad\text{for all }x\in\mathbb R^2\text{ and }\rho>0.
\]
For $\theta\in\operatorname{supp}(\mu_1)$,
\[
\dim_H\pi_\theta(\mathcal A)<\eta
\quad\Rightarrow\quad
\pi_\theta(\mathcal A)\subset\bigcup_{k\geq k_0}\bigcup_{I\in\mathcal G_{\theta,k}}I
\quad\text{for any }k_0\text{ large},
\]
where $\mathcal G_{\theta,k}$ consists of $\leq2^{\eta k}$ $2^{-k}$-intervals. Hence,
\begin{align*}
1 &=\mathbb P(\mathcal A)
\leq\sum_{k\geq k_0}\sum_{I\in\mathcal G_{\theta,k}}
\mathbb P(\mathcal A\cap\pi_\theta^{-1}(I))\\
\Rightarrow\quad
&\sum_{k\geq k_0}\int\sum_{I\in\mathcal G_{\theta,k}}
\mathbb P(\mathcal A\cap\pi_\theta^{-1}(I))\,d\mu_1(\theta)\geq1
\end{align*}
\begin{equation}
\Rightarrow\quad\exists k>k_0\text{ s.t. }\int\sum_{I\in\mathcal G_{\theta,k}}
\mathbb P(\mathcal A\cap\pi_\theta^{-1}(I))\,d\mu_1(\theta)\gtrsim\frac1{k^2}.
\tag{7.9}\label{eq:7.9}
\end{equation}
At this stage, we discretized the problem at scale $\delta=2^{-k}$. Suppose $\mathcal A_\theta=\bigcup_{I\in\mathcal G_{\theta,k}}\mathcal A\cap\pi_\theta^{-1}(I)$, then
\begin{align*}
N(\pi_\theta(\mathcal A_\theta),\delta)&=|\mathcal G_{\theta,k}|\leq\delta^{-\eta}
\quad\text{for all }\theta\in\operatorname{supp}(\mu_1)\\
\Rightarrow\quad |\pi_\theta(\mathcal A_\theta)|&\leq\delta^{1-\eta}
\quad\text{for all }\theta\in\operatorname{supp}(\mu_1).
\end{align*}
And by \paperref{eq:7.9}{(\ref*{eq:7.9})},

\[
\int\mathbb P(\mathcal A_\theta)\,d\mu_1(\theta)
\gtrsim\frac1{k^2}\sim\left(\log\frac1\delta\right)^{-2}
\]
\[
\Rightarrow\quad
\mu_1\left(\left\{\theta\in\operatorname{supp}(\mu_1):
\mathbb P(\mathcal A_\theta)\gtrsim\left(\log\frac1\delta\right)^{-2}\right\}\right)
\gtrsim\left(\log\frac1\delta\right)^{-2},
\]
Recall $\mathbb P$ is an $\alpha$-dimensional Frostman measure on $\mathcal A$. The intended discretization, after refinement, has the following properties.
\annot{\textbf{Proof requirement.} This refinement must retain the direction--point incidence mass while passing from the Frostman measure to a uniform discrete set. It is not enough to take an arbitrary $\delta$-net of the support.}
\begin{itemize}
\item $|\mathcal A|=\delta^{2-\alpha}$ and satisfies the non-concentration property:
\[
|\mathcal A\cap B(x,\rho)|\leq\rho^\alpha|\mathcal A|
\quad\text{for any }x\in\mathbb R^2\text{ and }\delta<\rho<1.
\]
\item $|\mathcal A_\theta|\gtrsim\delta^\varepsilon|\mathcal A|$.
\item $|\pi_\theta(\mathcal A_\theta)|\leq\delta^{1-\eta}$ for a set with $\mu_1$-measure $\gtrsim(\log(1/\delta))^{-C}$.
\end{itemize}
\annot{\textbf{Proof requirement.} The sets $\mathcal A_\theta$ depend on $\theta$. Thus the existence statement in \paperref{result:0.8}{Theorem~\ref*{result:0.8}} for the whole set does not directly give a contradiction. One needs the robust version, together with a refinement or exhaustion argument retaining their mass in the common structured pieces. After that, the bad-direction contribution is bounded by $\delta^c$ for some $c>0$, contradicting the retained lower bound $(\log(1/\delta))^{-C}$. Keeping this logarithmic lower bound is essential; replacing it by $\delta^\kappa$ does not justify comparison with $\delta^{\kappa\tau_0/10}$. This is the additional step required to complete the Hausdorff dimension proof sketch.}

\subsection{Discretized projection theorem in higher dimension}\label{sec:7.3}
Let's sketch the proof of \paperref{result:0.10}{Theorem~\ref*{result:0.10}}, the discretized projection theorem in higher dimension. We argue again by contradiction. Assume there are $\xi_1,\ldots,\xi_d,\zeta\in S^{d-1}$ such that
\[
\pi_{\xi_i}(E)<\delta^{1-\alpha/d-},\qquad
\pi_\zeta(E)<\delta^{1-\alpha/d-},
\]
\[
\&\quad E\subset A_1\times\cdots\times A_d,
\quad\text{where }A_i=\pi_{\xi_i}(E).
\]
Also, since $|E|=\delta^{d-\alpha}$, we have
\begin{itemize}
\item $|A_i|\approx\delta^{1-\alpha/d}$, $i=1,\ldots,d$.
\item $|\{x_1+\cdots+x_d:(x_1,\ldots,x_d)\in E\}|=|\pi_\zeta(E)|<\delta^{1-\alpha/d-}$.
\end{itemize}
\annot{As in the planar argument, the $d+1$ directions are chosen quantitatively transverse using the hyperplane non-concentration assumption. The change to coordinate directions and the sum direction costs $\delta^{-O_d(\tau_0)}$, absorbed into the displayed $+$ and $-$.}
At this point, we invoke a natural extension of the \paperref{result:1.6}{Balog--Szemer\'edi--Gowers Theorem} for $d\geq2$.

\begin{proposition}\label{result:7.1}
Fix $d\in\mathbb Z_+$, $d\geq2$. There exists a constant $C=C(d)$ such that the following holds. Let $A_1,\ldots,A_d$ be finite subsets of an additive group $Z$ and $G\subset A_1\times\cdots\times A_d$, such that for some $K>1$,
\begin{itemize}
\item $N<|A_i|<KN$, $i=1,\ldots,d$.
\end{itemize}

\begin{itemize}
\item $|G|>\frac1K N^d$,
\item $|\bigoplus_G A_i|=|\{x_1+\cdots+x_d:(x_1,\ldots,x_d)\in G\}|<KN$.
\end{itemize}
Then there are subsets $B_i\subset A_i$ ($1\leq i\leq d$) satisfying
\begin{itemize}
\item $|B_i|>K^{-C}N$, $i=1,\ldots,d$;
\item $|B_i+B_j|<K^CN$, $1\leq i,j\leq d$;
\item $|G\cap(B_1\times\cdots\times B_d)|>K^{-C}N^d$.
\end{itemize}
\end{proposition}
\begin{proof}
When $d=2$, this is the \paperref{result:1.6}{Balog--Szemer\'edi--Gowers theorem}. Assume $d\geq3$ and $\log K\ll\log N$.

Define
\[
H=\left\{(x_1,x_2)\in A_1\times A_2:
|G(x_1,x_2)|>\frac12K^{-3}N^{d-2}\right\},
\]
where $G(x_1,x_2)=\{(x_3,\ldots,x_d)\in A_3\times\cdots\times A_d:(x_1,x_2,x_3,\ldots,x_d)\in G\}$.

Note that
\[
|G\cap(H\times A_3\times\cdots\times A_d)|\sim|G|>\frac{N^d}{K}.
\]

\[
\&\quad |G\cap H\times A_3\times\cdots\times A_d|
=\sum_{(a_3,\ldots,a_d)\in A_3\times\cdots\times A_d}
|\{(x_1,x_2)\in H:(x_1,x_2,a_3,\ldots,a_d)\in G\}|
\]
\[
\Rightarrow\quad\exists(a_3,\ldots,a_d)\in A_3\times\cdots\times A_d\text{ s.t.}
\]
\[
H':=\{(x_1,x_2)\in H:(x_1,x_2,a_3,\ldots,a_d)\in G\}
\]
\[
\&\quad |H'|\geq\frac{|G\cap(H\times A_3\times\cdots\times A_d)|}{|A_3|\cdots|A_d|}
>\frac12K^{-d+1}N^2.
\]
Clearly,
\begin{align*}
|A_1\overset{H'}{+}A_2|
&=|\{x_1+x_2:(x_1,x_2)\in H'\}|\\
&=|\{x_1+x_2+a_3+\cdots+a_d:(x_1,x_2)\in H'\}|\\
&\leq|\bigoplus_G A_i|\\
&<KN.
\end{align*}
Apply the \paperref{result:1.6}{Balog--Szemer\'edi--Gowers theorem}, we obtain $B_1\subset A_1$, $B_2\subset A_2$ such that the following:
\begin{itemize}
\item $|B_1|>K^{-C}N$, $|B_2|>K^{-C}N$;
\item $|B_1+B_2|<K^CN$;
\item $|H'\cap(B_1\times B_2)|>K^{-C}N^2$.
\end{itemize}

Taking $G'=(B_1\times B_2\times A_3\times\cdots\times A_d)\cap G$, then
\begin{align*}
|G'|&=\sum_{(x_1,x_2)\in B_1\times B_2}|G(x_1,x_2)|\\
&\geq\sum_{(x_1,x_2)\in H'\cap(B_1\times B_2)}|G(x_1,x_2)|
>\frac12K^{-3}N^{d-2}|H'\cap(B_1\times B_2)|\\
&>\frac12K^{-C-3}N^d.
\end{align*}
Replace $A_i$ by $B_i$ for $i=1,2$, and $G$ with $G'$, we have proved proposition for $i=1,2$.

Iterate this procedure yields the desired conclusion.
\end{proof}

Now we continue to prove \paperref{result:0.10}{Theorem~\ref*{result:0.10}}. Apply \paperref{result:7.1}{Proposition~\ref*{result:7.1}} to the partial sumset $\bigoplus_E A_i$, we obtain sets $B_i\subset A_i$ ($1\leq i\leq d$) such that
\begin{itemize}
\item $N(B_i,\delta)>\delta^{-\alpha/d+}$, $i=1,\ldots,d$;
\item $N(B_i+B_j,\delta)<\delta^{-\alpha/d-}$, $1\leq i,j\leq d$;
\item $|E\cap(B_1\times\cdots\times B_d)|>\delta^{d-\alpha+}$.
\end{itemize}
Define
\[
A=\left\{x_1\in B_1:
|E(x_1)\cap(B_2\times\cdots\times B_d)|>
\delta^{d-1-\frac{d-1}{d}\alpha+}\right\}
\]
\[
\&\quad\mathcal B=E\cap(B_1\times\cdots\times B_d).
\]
Then we have

\[
|\mathcal B|\sim|E\cap(A\times B_2\times\cdots\times B_d)|>\delta^{d-\alpha+}
\]
\[
\Rightarrow\quad |A|>\delta^{1-\alpha/d+}.
\]
Again, by an appropriate choice of $\xi_1$, $A\subset\pi_{\xi_1}(E)$ satisfies the desired non-concentration property. Consider the map
\[
\varphi:\mathbb R^d\longrightarrow\mathbb R,
\quad(\xi_1,\ldots,\xi_d)\longmapsto\frac{\xi_2}{\xi_1},
\]
and let $\mu=\varphi[\mu_1]$, then
\begin{align*}
\mu(I(a,\rho))
&=\mu_1(\varphi^{-1}(I(a,\rho)))\\
&\leq\mu_1(\{(y_1,\ldots,y_d):|y_2-ay_1|\leq2\rho\})\\
&=\mu_1(V((-a,1,0,\ldots,0)^\perp,2\rho))\\
&\lesssim\rho^\kappa
\end{align*}
for any $a\in\mathbb R$, $\delta<\rho<\delta^{\tau_0}$.

From the conclusion of \paperref{result:0.7}{Theorem~\ref*{result:0.7}}, we obtain some $\varepsilon_1>0$ and $x\in\operatorname{supp}(\mu)$ such that
\[
|A+xA|>\delta^{-\varepsilon_1}|A|>\delta^{1-\alpha/d-\varepsilon_1/2}.
\]
Choose a preimage $\theta=(\theta_1,\ldots,\theta_d)\in\operatorname{supp}(\mu_1)$ of $x$ in this compact chart, so that $\theta_2/\theta_1=x$. Then

\[
|\theta_1A+\theta_2A|=|\theta_1|\,|A+xA|\gtrsim_d\delta^{1-\alpha/d-\varepsilon_1/2}
\]
\[
\Rightarrow\quad |\pi_\theta(A^d)|=|\theta_1A+\cdots+\theta_dA|
>|\theta_1A+\theta_2A|>\delta^{1-\alpha/d-\varepsilon_1/2}.
\]
By the same averaging argument as before, we have
\begin{align*}
\chi_{\pi_\theta(A^d)}
&=\chi_{A^d}\circ\pi_\theta^{-1}
\leq\frac1{|\mathcal B|}\int_{\mathcal B-A^d}
(\chi_{\mathcal B-\xi}\circ\pi_\theta^{-1})\,d\xi\\
&\leq\frac1{|\mathcal B|}\int_{\prod_{i=1}^d(B_i-A)}
\chi_{\pi_\theta(\mathcal B)-\pi_\theta(\xi)}\,d\xi,
\end{align*}
which implies
\begin{align*}
\delta^{1-\alpha/d-\varepsilon_1/2}
&<|\pi_\theta(A^d)|\leq|\mathcal B|^{-1}\cdot
\prod_{i=1}^d|B_i-A|\cdot|\pi_\theta(\mathcal B)|\\
&<\delta^{-d+\alpha-}\cdot\delta^{d-\alpha-}\cdot\delta^{1-\alpha/d-}\\
\Rightarrow\quad \delta^{-\varepsilon_1/2}&<\delta^{-o(1)},
\end{align*}
and hence a contradiction.
\annot{\textbf{Proof requirement.} This calculation gives the higher-dimensional weak projection mechanism. The full conclusion \paperref{eq:0.6}{(\ref*{eq:0.6})}, including one retained set and a large set of directions that work for every sufficiently large subset, requires the additional robust extraction argument discussed above. It is also needed in the Hausdorff-dimension passage for \paperref{result:0.11}{Theorem~\ref*{result:0.11}}.}

\subsection{Dimension of projection set in higher dimension}\label{sec:7.4}
Finally, let's briefly point out how the upper bound $d-2+\kappa$ in the \paperref{result:0.11}{Theorem~\ref*{result:0.11}} are obtained.

Set
\[
F=\{\theta\in\mathbb P^{d-1}:\dim_H\pi_\theta(E)<\eta\}
\]
and assume
\[
\dim_H F>\widetilde\kappa
\quad\Rightarrow\quad\exists\mu_1\text{ is }\widetilde\kappa\text{-dim Frostman on }\mathbb P^{d-1}.
\]
To apply \paperref{result:0.10}{Theorem~\ref*{result:0.10}}, we need to check that
\[
\mu_1(V(y^\perp,\rho))<\rho^\kappa
\quad\text{for all }\bar y\in\mathbb P^{d-1}\text{ and }\delta<\rho<\delta^{\tau_0}.
\]
Indeed,
\[
\mu_1(V(y^\perp,\rho))\sim
\sum_{\substack{\tau(x,\rho)\subset S^{d-1}\\
\text{cover }V(y^\perp,\rho)\cap S^{d-1}}}
\mu_1(\tau(x,\rho))
<\frac\rho{\rho^{d-1}}\cdot\rho^{\widetilde\kappa}
=\rho^{\widetilde\kappa-(d-2)},
\]
thus we take $\widetilde\kappa=d-2+\kappa$, as desired.

\textbf{P.S.} In higher dimension, P. Cholak, M. Cs\"ornyei, N. Lutz, P. Lutz, E. Mayordomo and D. M. Stull \cite{CholakEtAl}, have obtain sharp exceptional set estimates for orthogonal projections onto lines and hyperplanes.

\begin{theorem}\label{result:7.2}
Let $X\subset\mathbb R^n$ be Borel, then for
\[
\max\{\dim_H X-n+k,0\}<s\leq\min\{\dim_H X,k\}
\]
(for $k=1$ and $k=n-1$ respectively),
\begin{itemize}
\item
\[
\begin{aligned}
&\dim_H\{\theta\in S^{n-1}:\dim_H P_\theta(X)<s\}\\
&\qquad\leq\max\{n-2-\lfloor\dim_H X-s\rfloor,\ n-2+2s-\dim_H X\}.
\end{aligned}
\]
\item
\[
\begin{aligned}
&\dim_H\{V\in\operatorname{Gr}(n,n-1):\dim_H P_V(X)<s\}\\
&\qquad\leq\max\{2s-\dim_H X,\ \lceil s\rceil-1\}.
\end{aligned}
\]
\end{itemize}
\end{theorem}
In particular, when $n=3$, the upper bounds as following:

\begin{center}
\phantomsection\label{fig:projection-one}\includegraphics[width=\linewidth]{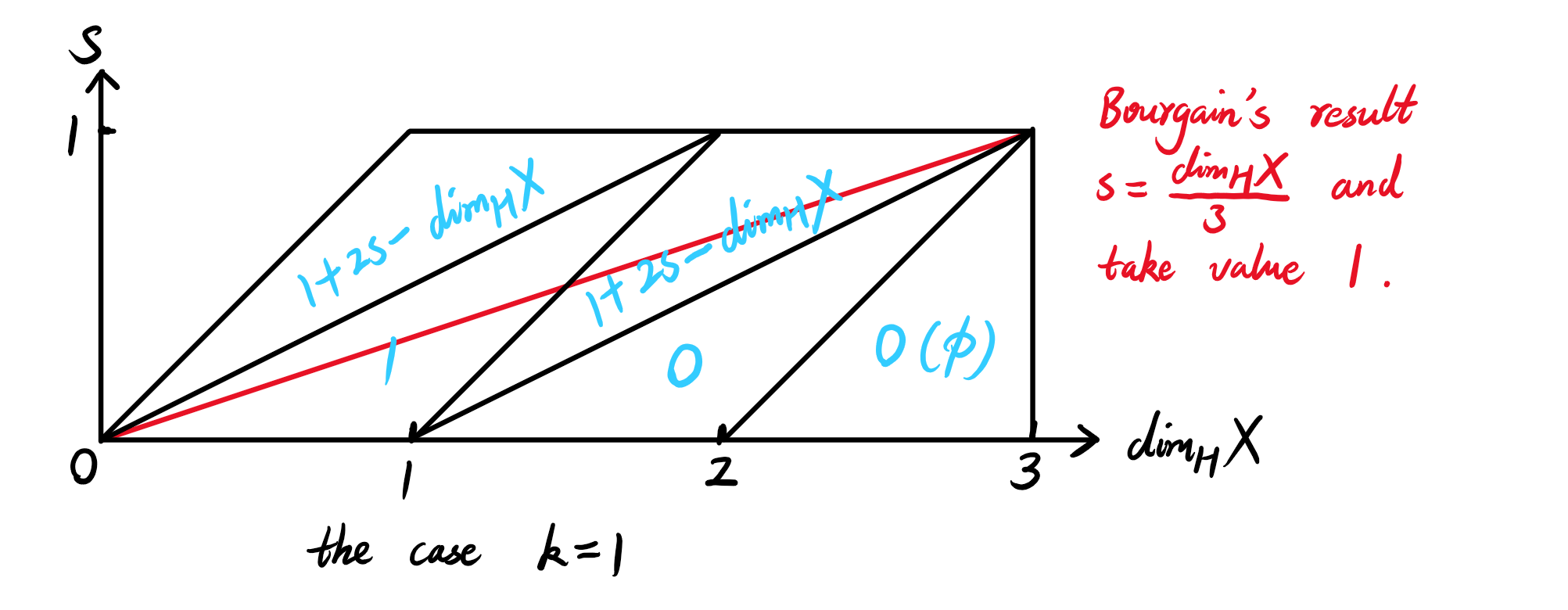}
\end{center}
\annot{The red line records Bourgain's bound $\dim_H(\text{exceptional directions})\leq1$ at $s=(\dim_H X)/3$. It does not say that the sharp bound in \paperref{result:7.2}{Theorem~\ref*{result:7.2}} is identically $1$ on this line: for example, at $\dim_H X=2$ and $s=2/3$ the sharp bound is $1/3$.}

\begin{center}
\includegraphics[width=.88\linewidth]{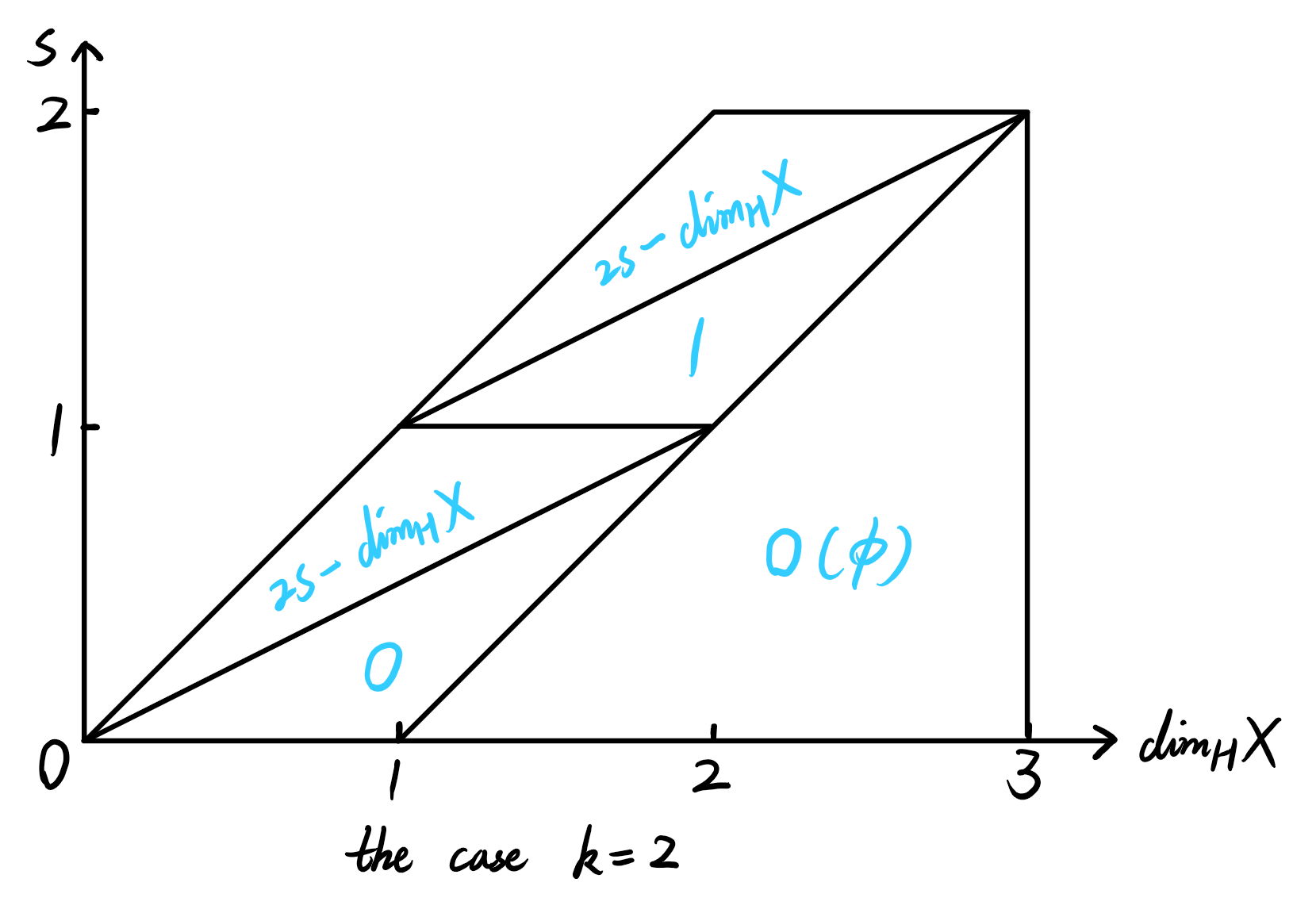}
\end{center}

\section{Reduce ring conjecture to discretized sum-product theorem}\label{sec:8}

In this section, we prove there is no Borel subring of $\mathbb R$ which has Hausdorff dimension strictly intermediate between $0$ and $1$, by using the \paperref{result:0.3}{discretized sum-product theorem}.

A tempting approach is to discretize a compact subset of the subring and
apply the \paperref{result:0.3}{discretized sum-product theorem} directly. Two difficulties prevent this
from being immediate. First, the sum and product of that compact subset
remain in the ring, but need not remain in the same set. Second,
Hausdorff dimension does not give a uniform upper bound for covering
numbers at every small scale. The closure of ring alone therefore does not
produce a finite set for which both the sumset and the product set
are small at a common scale.

The Fourier argument provides a way to construct such a set. \paperref{step:8-1}{Step 1} expresses the dimension condition as a lower bound on Fourier decay for measures supported on the ring. \paperref{step:8-2}{Step 2} uses the
preservation of this support under ring operations to force many multiplicative coincidences in a level set of a smoothed multiplicative correlation (a concept comes from signal processing). \paperref{step:8-3}{Step 3} uses the correlation structure of the same function to obtain additive coincidences. Two successive applications of the \paperref{result:1.6}{Balog--Szemer\'edi--Gowers theorem} then impose the two forms of structure on one sufficiently large, non-concentrated set.

\medskip\noindent\phantomsection\label{step:8-1}\textbf{Step 1. Reduce the Hausdorff dimension of a ring to Fourier concentration.}
Let $A\subset\mathbb R$ be a Borel ring with $\dim_H A=\alpha$, $0<\alpha<1$. Denote by $\mathcal P_c(A)$ the compactly supported probability measures whose supports are contained in $A$. For $\mu\in\mathcal P_c(A)$, take
\[
\beta(\mu)=\limsup_{N\to\infty}
\frac{\log\left(\int_{-N}^{N}|\widehat\mu(\xi)|\,d\xi\right)}{\log N},
\qquad \beta=\inf_{\mu\in\mathcal P_c(A)}\beta(\mu).
\]
Here, $\beta(\mu)$ records the average Fourier decay exponent. of the measure, rather than its pointwise Fourier decay exponent. Also, the use of a "limsup" is deliberate: later we
need suitable arbitrarily large scales, rather than a lower bound at every large scale.

\begin{claim}
$\beta=1-\alpha$.
\end{claim}

\begin{proof} (i) To prove $\beta\geq1-\alpha$, fix $\eta>0$. Since $\dim_H A<\alpha+\eta$, for every large $k_0$ there are $2^{-k-1}$-separated sets $S_k$, $k\geq k_0$, such that

\[
A\subset\bigcup_{k\geq k_0}\bigcup_{p\in S_k}I(p,2^{-k}),
\qquad |S_k|\lesssim2^{k(\alpha+\eta)}.
\]
For every $\mu\in\mathcal P_c(A)$,
\[
1\leq\sum_{k\geq k_0}\sum_{p\in S_k}\mu(I(p,2^{-k})),
\]
so for some $k\geq k_0$ the inner sum is $\gtrsim k^{-2}$. Fourier inversion and Cauchy--Schwarz give
\[
\sum_{p\in S_k}\mu(I(p,2^{-k}))
\lesssim2^{-k}
\left(\int_{|\xi|\leq C2^k}|\widehat\mu(\xi)|^2\,d\xi\right)^{1/2}
\left(\int_{|\xi|\leq C2^k}\left|\sum_{p\in S_k}e^{2\pi ip\xi}\right|^2\,d\xi\right)^{1/2}.
\]

The separation of $S_k$ implies
\[
\int_{|\xi|\leq C2^k}\left|\sum_{p\in S_k}e^{2\pi ip\xi}\right|^2\,d\xi
\lesssim2^k|S_k|.
\]
Thus
\[
\int_{|\xi|\leq C2^k}|\widehat\mu(\xi)|^2\,d\xi
\gtrsim k^{-4}2^{k(1-\alpha-\eta)}.
\]
Since $|\widehat\mu|\leq1$ and $k_0$ is arbitrary, $\beta(\mu)\geq1-\alpha-\eta$. Letting $\eta\downarrow0$ proves the lower bound.

(ii) By Frostman's lemma, there is $\mu\in\mathcal P_c(A)$ such that
\[
\mu(I(p,r))\lesssim_\eta r^{\alpha-\eta}.
\]
Write $\mu^-(E)=\mu(-E)$. Since $A$ is a ring, $\mu*\mu^-\in\mathcal P_c(A)$ and
\[
|(\mu*\mu^-)^{\wedge}(\xi)|=|\widehat\mu(\xi)|^2.
\]
we have
\[
\int_{|\xi|\leq N}|\widehat\mu(\xi)|^2\,d\xi
\lesssim N\int\mu(I(x,C/N))\,d\mu(x)
\lesssim_\eta N^{1-\alpha+\eta}.
\]
Hence,
\[
\beta\leq\beta(\mu*\mu^-)\leq1-\alpha+\eta.
\]
Let $\eta\downarrow0$, we prove the claim.
\end{proof}

\medskip\noindent\phantomsection\label{step:8-2}\textbf{Step 2. Exploit the multiplicative structure.}
We will find a large non-concentrated discretized set with small product set. Choose the preceding Frostman measure on a compact subset of $A\cap[1,2]$. For $x\in A$, let $\rho_x(a)=xa$. Then $\rho_x[\mu]\in\mathcal P_c(A)$. More generally, for $\lambda,\nu\in\mathcal P_c(A)$,
\[
\int(\rho_x[\lambda]*\lambda)\,d\nu(x)\in\mathcal P_c(A).
\]

Indeed, its support is contained in the compact set
\[
\bigcup_{x\in\operatorname{supp}(\nu)}
\bigl(x\operatorname{supp}(\lambda)+\operatorname{supp}(\lambda)\bigr)\subset A.
\]
Its Fourier transform is
\[
\int\widehat\lambda(x\xi)\widehat\lambda(\xi)\,d\nu(x).
\]
Taking $\lambda=\mu*\mu^-$ and $\nu=\mu$, and using the definition of $\beta$, gives
\begin{equation}
\limsup_{N\to\infty}
\frac{\log\left(\int_{-N}^{N}\int
|\widehat\mu(y\xi)|^2|\widehat\mu(\xi)|^2\,d\mu(y)\,d\xi\right)}{\log N}
\geq\beta.
\tag{8.1}\label{eq:8.1}
\end{equation}

Define a even bump function
\begin{align*}
\phi(x)&\approx\int \widehat{\chi}_{I(0,1)}(\frac{\xi}{N})|\widehat{\mu}(\xi)|^2e^{2\pi ix\xi}d\xi\\
 &\approx [N\chi_{I(0,\frac{1}{N})}*\mu*\mu](x).
\end{align*}
This bump function satisfies
\begin{equation}
\|\phi\|_1=1,
\qquad \|\phi\|_\infty\lesssim_\eta N^{1-\alpha+\eta}=N^{\beta+\eta}.
\tag{8.2}\label{eq:8.2}
\end{equation}
Moreover,
\begin{equation}
\int_{-cN}^{cN}\int|\widehat\mu(y\xi)|^2|\widehat\mu(\xi)|^2\,d\mu(y)\,d\xi
\lesssim\iint\phi(x)\phi(xy)\,dx\,d\mu(y).
\tag{8.3}\label{eq:8.3}
\end{equation}

\annot{\paperref{eq:8.1}{(\ref*{eq:8.1})}--\paperref{eq:8.3}{(\ref*{eq:8.3})} force the averaged correlation to $\gtrapprox N^{\beta}$, along a sequence of scales. A large correlation requires a substantial contribution from points where both factors are large. Passing to a suitable superlevel set $S_1$ turns this into overlap between $S_1$ and its dilates. Many such overlaps yield many approximate multiplicative relations between the discretized points. The first \paperref{result:1.6}{BSG} step converts these relations into a large subset with a small product set. Thus the
Fourier argument is producing the multiplicative structure that direct discretization did not provide.}

By \paperref{eq:8.1}{(\ref*{eq:8.1})}, for arbitrarily large $N$ the right-hand side is at least $N^{\beta-\eta}$. Choose
\[
T_1=N^{\beta-C_1\eta},\qquad S_1=\{x:\phi(x)>T_1\},
\]
where $C_1$ is sufficiently large. Thus
\begin{equation}
\int|S_1\cap yS_1|\,d\mu(y)\geq N^{-\beta-O(\eta)}.
\tag{8.4}\label{eq:8.4}
\end{equation}
Here $|S_1\cap y^{-1}S_1|=y^{-1}|S_1\cap yS_1|$, so replacing $y^{-1}$ by $y$ only changes fixed constants. Moreover
\[
N^{-\beta-O(\eta)}\leq|S_1|\leq T_1^{-1}=N^{-\beta+O(\eta)}.
\]
To discretize $S_1$ at the common scale $\delta$, we retain a union of intervals, still denoted $S_1$, whose centres $z_j$ satisfy $\phi(z_j)\geq cT_1$ and
\begin{equation}
S_1=\bigcup_{j=1}^{M}I(z_j,\delta),\qquad
M=N^{\alpha+O(\eta)}.
\tag{8.5}\label{eq:8.5}
\end{equation}
The centres are $\delta$-separated. For $N^{-1}\leq r\leq1$, the Frostman condition implies
\[
|S_1\cap I(p,r)|\lesssim T_1^{-1}\int_{I(p,Cr)}\phi
\lesssim_\eta N^{-\beta+O(\eta)}r^{\alpha-\eta}.
\]
For $\delta\leq r<N^{-1}$, use $|S_1\cap I(p,r)|\lesssim r$. Together these estimates give the $(\delta,\alpha)_1$ non-concentration condition.

Choose a further small $\epsilon>0$ and then $\eta\ll\alpha\epsilon$. The preceding bound shows
\[
|S_1\cap I(0,N^{-\epsilon})|\leq N^{-\beta-\alpha\epsilon+O(\eta)}.
\]
This is negligible compared with the lower bound in \paperref{eq:8.4}{(\ref*{eq:8.4})}, so we may retain only $|z_j|>N^{-\epsilon}$. Since $\phi$ is even and $y>0$, we may also take the positive half and reflect if necessary. Thus $N^{-\epsilon}<z_j<C$.

By \paperref{eq:8.4}{(\ref*{eq:8.4})}, the set of multipliers for which
\[
|S_1\cap yS_1|\geq N^{-\beta-O(\eta)}
\]
has $\mu$-measure at least $N^{-O(\eta)}$. By Frostman condition, this set contains $N^{\alpha-O(\eta)}$ multipliers $y_k$ separated by $N^{-1}$. Consequently,
\[
m_k:=\#\{(j_1,j_2):|z_{j_1}-y_kz_{j_2}|<N^{-1+O(\eta)}\}
\geq N^{\alpha-O(\eta)}.
\]
A fixed quadruple can be counted for at most $N^{O(\epsilon+\eta)}$ of these multipliers, since $z_{j_2}\geq N^{-\epsilon}$. Hence
\[
\begin{aligned}
&\#\left\{(j_1,j_2,j_3,j_4):
\left|1-\frac{z_{j_1}z_{j_3}}{z_{j_2}z_{j_4}}\right|
<N^{-1+O(\epsilon+\eta)}\right\}\\
&\qquad\gtrsim N^{-O(\epsilon+\eta)}\sum_km_k^2
\geq N^{3\alpha-O(\epsilon+\eta)}.
\end{aligned}
\]
So, the additive energy of $\{\log z_j\}$ is therefore at least $N^{3\alpha-O(\epsilon+\eta)}$. The \paperref{result:1.6}{Balog--Szemer\'edi--Gowers theorem} gives a subset $\{z_j:j\in\mathcal F_1\}$ with $|\mathcal F_1|\geq N^{\alpha-O(\epsilon+\eta)}$ and
\begin{equation}
N(\{z_{j_1}z_{j_2}:j_1,j_2\in\mathcal F_1\},\delta)
\leq N^{\alpha+O(\epsilon+\eta)}.
\tag{8.6}\label{eq:8.6}
\end{equation}

\medskip\noindent\phantomsection\label{step:8-3}\textbf{Step 3. Consider the additive property.}

We now find a subset of $\{z_j\}_{j\in\mathcal F_1}$ with small sumset. Similar as before, we define a bump function
\begin{align*}
\psi\approx\int \bigl( \widehat{\chi}_{I(0,1)}(\frac{\xi}{N}) \bigr)^{\frac{1}{2}}\widehat{\mu}(\xi)e^{2\pi ix\xi}d\xi\\
\approx [N\chi_{I(0,\frac{1}{N})}*\mu](x).
\end{align*}
Moreover, $\psi$ satisfies
\begin{equation}
\|\psi\|_1=1,\qquad \|\psi\|_\infty\lesssim_\eta N^{\beta+\eta},
\qquad \phi(z)=\int\psi(x)\psi(x-z)\,dx.
\tag{8.7}\label{eq:8.7}
\end{equation}

\annot{The useful extra information about the points retained by the first
\paperref{result:1.6}{BSG} step is that they still satisfy $\phi(z_j)\geq cT_1$. By \paperref{eq:8.7}{(\ref*{eq:8.7})},
this means that the autocorrelation of $\psi$ is large. Cauchy--Schwarz then
forces a large total of the mutual inner products of translates for $\psi$
Those mutual inner products are precisely the values $\phi(z_{j_2}-z_{j_1})$. At the common discretization scale, many pairs therefore have their differences in
relatively few intervals. Cauchy--Schwarz converts this into large
additive energy, and a second \paperref{result:1.6}{BSG} application gives a subset with a
small sumset. The source of this additive structure is the
autocorrelation identity, not the small product-set conclusion alone.}

Write $M_1=|\mathcal F_1|$. Since $\phi(z_j)\geq cT_1$, the Cauchy--Schwarz inequality gives
\[
\begin{aligned}
cT_1M_1
&\leq\sum_{j\in\mathcal F_1}\phi(z_j)
=\int\psi(x)\sum_{j\in\mathcal F_1}\psi(x-z_j)\,dx\\
&\leq\|\psi\|_2
\left(\int\left|\sum_{j\in\mathcal F_1}\psi(x-z_j)\right|^2dx\right)^{1/2}.
\end{aligned}
\]

Take $H_0=\|\psi\|_2^2=\phi(0)=\|\phi\|_\infty\leq N^{\beta+O(\eta)}$. Expanding the square, the last quantity equals
\[
H_0^{1/2}\left(\sum_{j_1,j_2\in\mathcal F_1}
\phi(z_{j_2}-z_{j_1})\right)^{1/2}.
\]
Thus, for a fixed $c_0>0$,
\[
\sum_{j_1,j_2\in\mathcal F_1}\phi(z_{j_2}-z_{j_1})
\geq c_0T_1^2M_1^2/H_0.
\]
To infer many large values at the differences, let
\[
T_2=\frac{c_0T_1^2}{2H_0},\qquad S_{\mathrm{add}}=\{\phi>T_2\}.
\]
The pairs outside this set contribute at most $T_2M_1^2$. Therefore
\[
\#\{(j_1,j_2)\in\mathcal F_1^2:z_{j_2}-z_{j_1}\in S_{\mathrm{add}}\}
\gtrsim\frac{T_1^2M_1^2}{H_0^2}
\geq N^{2\alpha-O(\epsilon+\eta)}.
\]
Moreover $|S_{\mathrm{add}}|\leq T_2^{-1}=N^{-\beta+O(\eta)}$. The additive energy consequently satisfies
\[
\begin{aligned}
&\#\{(j_1,j_2,j_3,j_4)\in\mathcal F_1^4:
|z_{j_1}-z_{j_2}+z_{j_3}-z_{j_4}|<C\delta\}\\
&\qquad\gtrsim\sum_jr_j^2
\geq\frac{(\sum_jr_j)^2}{M_{\mathrm{add}}}
\geq N^{3\alpha-O(\epsilon+\eta)}.
\end{aligned}
\]
A second application of \paperref{result:1.6}{BSG} gives $\mathcal F_2\subset\mathcal F_1$, $|\mathcal F_2|\geq N^{\alpha-O(\epsilon+\eta)}$, such that
\begin{equation}
N(\{z_{j_1}+z_{j_2}:j_1,j_2\in\mathcal F_2\},\delta)
\leq N^{\alpha+O(\epsilon+\eta)}.
\tag{8.8}\label{eq:8.8}
\end{equation}

\annot{The two \paperref{result:1.6}{BSG} selections must be nested. Two unrelated subsets, one with
a small sumset and another with a small product set, would not
contradict the \paperref{result:0.3}{sum-product theorem}. Here
$\mathcal F_2\subset\mathcal F_1$ preserves the product-set bound
from the first selection, while the second selection supplies the
sumset bound.}

Thus, the set $U=\bigcup_{j\in\mathcal F_2}I(z_j,\delta)$ remains a subset of the original non-concentrated set $S_1$. It has the required lower size, while \paperref{eq:8.6}{(\ref*{eq:8.6})} and \paperref{eq:8.8}{(\ref*{eq:8.8})} bound its product and sum covering numbers by $\delta^{-\alpha-O(\epsilon+\eta)}$. The \paperref{result:0.3}{discretized sum-product theorem} gives instead
\[
N(U+U,\delta)+N(U\cdot U,\delta)>\delta^{-\alpha-c(\alpha)}
\]
for some $c(\alpha)>0$. Choose $\epsilon$ and then $\eta$ sufficiently small relative to this gain. We obtain a contradiction.

Therefore $\alpha=0$ or $1$, completing the proof.

\section{Fourier decay of multiplicative convolutions of measures}\label{sec:9}

In this section, we prove \paperref{result:0.12}{Theorem~\ref*{result:0.12}}. The standard proof of \paperref{result:0.13}{Theorem~\ref*{result:0.13}} uses Fourier analysis and potential theory, and is omitted.

Let us first explain the route from sum-product growth to Fourier decay. Multiplying independent variables averages the Fourier transform over dilated frequencies, so the required cancellation is not immediate. We first work with integrated high moments, which measure how much of the Fourier transform remains large. The \paperref{result:9.1}{finite-scale convolution estimate below} forces these moments to decrease as long as their exponents stay above a prescribed small threshold. Once this threshold is reached, one more independent multiplier, together with non-concentration, converts the integral estimate into pointwise decay. Finally, additive convolution raises the Fourier transform to a power and makes the resulting decay strong enough to control the mass of intervals. These are the separate roles of multiplication, high moments, and addition in the proof.

\subsection{A finite-scale convolution estimate.}\label{sec:9.1}

First, we establish a slight refinement of the additive/multiplicative measure convolution theorem obtained in \paperref{sec:8}{section~8}. The proof is the same, except that \paperref{result:0.6}{Theorem~\ref*{result:0.6}} is used to instead of \paperref{result:0.3}{Theorem~\ref*{result:0.3}} to derive the contradiction.

\begin{lemma}\label{result:9.1}
For $0<\sigma<1$ and $\kappa>0$, there is $\varepsilon>0$ such that the following holds.

Suppose $\mu\in\mathcal P([-1,1])$ and $N\in\mathbb{Z}_+$ a large integer. Assume
\[
E_N(\mu):=\int_{-N}^{N}|\widehat\mu(\xi)|^2\,d\xi=N^\sigma,
\]
\centerline{\textcolor{cyan}{(play the same role with $\mu*\mu\in\mathcal{P}_c(A)$ as before)}}
and for all $N^{-1}<\rho<N^{-\varepsilon}$,
\[
\sup_a\mu(I(a,\rho))\leq\rho^\kappa
\quad(N^{-1}<\rho<N^{-\varepsilon}).
\]
\centerline{\textcolor{cyan}{(the reason why we use \paperref{result:0.6}{Theorem~\ref*{result:0.6}} rather than \paperref{result:0.3}{Theorem~\ref*{result:0.3}})}}
Then
\begin{equation}
\int_{-N}^{N}\int|\widehat\mu(\xi)|^2|\widehat\mu(x\xi)|^2
\,d(\mu*\mu)(x)\,d\xi<N^{\sigma-\varepsilon}.
\tag{9.1}\label{eq:9.1}
\end{equation}
\centerline{\textcolor{cyan}{(compared to \paperref{eq:8.1}{(\ref*{eq:8.1})}, we replace $\nu=\mu$ by $\mu*\mu$)}}
The gain may be chosen uniformly for $\sigma$ in a fixed compact subinterval of $(0,1)$. The same assertion holds for measures supported in any fixed bounded interval.
\end{lemma}

\subsection[Version about Fourier decay of Theorem 0.12]{Version about Fourier decay of \paperref{result:0.12}{Theorem 0.12}}\label{sec:9.2}
\begin{theorem}\label{result:9.2}
For every $\kappa>0$, there are $\varepsilon_1,\varepsilon_2>0$ and $k\in\mathbb Z_+$ such that the following holds.

Let $\mu\in\mathcal P([\frac12,1])$ and $N$ be a large integer. Assume for all $N^{-1}<\rho<N^{-\varepsilon_1}$ that
\[
\sup_a\mu(I(a,\rho))\leq\rho^\kappa.
\]
Then for all $\xi\in\mathbb{R}$, $|\xi|\sim N$,
\[
\left|\int e^{2\pi i\xi x_1\cdots x_k}\,d\mu(x_1)\cdots d\mu(x_k)\right|
<N^{-\varepsilon_2}.
\]
\end{theorem}
\begin{proof}
Let $\mu_k$ be the image of $\mu^{\otimes k}$ under $(x_1,\ldots,x_k)\mapsto x_1\cdots x_k$.

With $e(z)=e^{2\pi iz}$,
\[
\widehat\mu_k(\xi)=\int e(\xi x_1\cdots x_k)\,d\mu(x_1)\cdots d\mu(x_k),
\qquad \widehat\mu_{k+l}(\xi)=\int\widehat\mu_k(y\xi)\,d\mu_l(y).
\]
Take $\Omega_{k,\tau}=\{\xi:|\widehat\mu_k(\xi)|>N^{-\tau}\}$. If $\xi\in\Omega_{2k,\tau}$, Jensen's inequality gives, for every integer $r\geq1$,
\[
N^{-2r\tau}<\int|\widehat\mu_k(y\xi)|^{2r}\,d\mu_k(y).
\]
Expanding the even power,
\[
|\widehat\mu_k(y\xi)|^{2r}
=\int e\bigl(y\xi(x_1-x_2+\cdots+x_{2r-1}-x_{2r})\bigr)
\,d\mu_k(x_1)\cdots d\mu_k(x_{2r}).
\]

Let $\lambda_k=\mu_k*\mu_k^-$ and write $\lambda^{*r}$ for additive convolution powers. Fubini gives
\[
N^{-2r\tau}<\int\widehat\mu_k(\xi z)\,d\lambda_k^{*r}(z).
\]
Replacing $r$ by $2r$ and applying H\"older,
\[
N^{-4r\tau}
<\int\widehat\mu_k(\xi z)\,d\lambda_k^{*(2r)}(z)
\leq\left(\int|\widehat\mu_k(\xi z)|^{4r}\,d\lambda_k^{*(2r)}(z)\right)^{1/(4r)}.
\]
Since $|\widehat\mu_{2k}|\leq1$ and $4r^2\geq4r$, we may decrease the integrand by replacing the latter exponent $4r$ with $4r^2$. Consequently
\begin{equation}
\begin{aligned}
&\int_{-N}^{N}\int|\widehat\mu_{2k}(\xi)|^{4r}
|\widehat\mu_k(z\xi)|^{4r}\,d\lambda_k^{*(2r)}(z)\,d\xi\\
&\qquad>N^{-16r^2\tau}\int_{\Omega_{2k,\tau}\cap[-N,N]}
|\widehat\mu_{2k}(\xi)|^{4r^2}\,d\xi\\
&\qquad>N^{-20r^2\tau}|\Omega_{2k,\tau}\cap[-N,N]|.
\end{aligned}\tag{9.2}\label{eq:9.2}
\end{equation}
Set $E_{k,r}(N)=\int_{-N}^{N}|\widehat\mu_k(\xi)|^r\,d\xi=N^{\sigma_{k,r}}$.
\annot{Since $|\widehat\mu_k|\leq1$, increasing $r$ suppresses the frequencies at which the Fourier transform is already small. A high moment therefore captures the part where it remains large. The sets $\Omega_{k,\tau}$ make this distinction explicit: on their complement, the integrand in a moment of order $r$ is at most $N^{-r\tau}$, whereas on the large value set it is at most $1$. Consequently, a bound for the size of the large value set can be converted into a bound for a sufficiently high moment. The iteration will use this freedom to increase both the number of multiplicative factors and the moment order.}

These values decrease in $r$, and
\[
E_{k+1,r}(N)\leq\int\int_{-N}^{N}|\widehat\mu_k(y\xi)|^r\,d\xi\,d\mu(y)
\leq2E_{k,r}(N).
\]
Their exponents are therefore essentially decreasing in $k$ as well. Furthermore,
\[
E_{1,2}(N)\lesssim N\sup_a\mu(I(a,C/N))\lesssim N^{1-\kappa}.
\]
Hence all the values used below have exponent at most $1-\kappa+o(1)$.

Fix a target $\sigma'>0$. Next, we start at $\sigma_{1,2}$ to decrease recursively $\sigma_{k,r}$, until $\sigma_{k,r}$ becomes $\sigma'$ sufficiently small. We use the uniform version of \paperref{result:9.1}{Lemma~\ref*{result:9.1}} on a compact exponent interval, with non-concentration exponent $\kappa/2$; write its gain as $\varepsilon_*>0$. At each stage, apply \paperref{result:9.1}{Lemma~\ref*{result:9.1}} with $\mu$ replaced with
\[
\nu=\frac12\bigl(\lambda_k^{*r}+\lambda_{2k}^{*r}\bigr),
\qquad \widehat\nu=\frac12\bigl(|\widehat\mu_k|^{2r}+|\widehat\mu_{2k}|^{2r}\bigr)\geq0.
\]
Then
\[
E_N(\nu)\lesssim_k E_{k,2r}(N),
\qquad E_{2k,4r}(N)\leq4E_N(\nu).
\]
If $E_N(\nu)<N^{\sigma'/2}$, the second inequality already gives $E_{2k,4r}(N)<N^{\sigma'}$ for large $N$, and the procedure terminates.

Otherwise the exponent of $E_N(\nu)$ lies in a compact interval contained in $(0,1)$, with lower endpoint $\sigma'/2$ and upper endpoint $1-\kappa/2$. Clearly, for $N^{-1}<\rho<N^{-\varepsilon_1}$,
\[
\begin{aligned}
\sup_a\nu(I(a,\rho))\leq\sup_a\mu_k(I(a,\rho))+\sup_a\mu_{2k}(I(a,\rho))\lesssim\rho^{\kappa/2}.
\end{aligned}
\]
Hence, by \paperref{result:9.1}{Lemma~\ref*{result:9.1}},
\[
\int_{-N}^{N}\int|\widehat\nu(\xi)|^2|\widehat\nu(z\xi)|^2
\,d(\nu*\nu)(z)\,d\xi
<N^{\sigma_{k,2r}-\varepsilon_*}.
\]
Now $|\widehat\nu(\xi)|^2\geq\frac14|\widehat\mu_{2k}(\xi)|^{4r}$,
$|\widehat\nu(z\xi)|^2\geq\frac14|\widehat\mu_k(z\xi)|^{4r}$, and
$\nu*\nu\geq\frac14\lambda_k^{*(2r)}$ as measures. Combining these three positive dominations with \paperref{eq:9.2}{(\ref*{eq:9.2})} yields
\[
N^{-20r^2\tau}|\Omega_{2k,\tau}\cap[-N,N]|
\lesssim N^{\sigma_{k,2r}-\varepsilon_*}.
\]
Therefore
\[
E_{2k,\bar r}(N)
\leq|\Omega_{2k,\tau}\cap[-N,N]|+2N^{1-\tau\bar r}
\lesssim N^{\sigma_{k,2r}-\varepsilon_*+20r^2\tau}+N^{1-\tau\bar r}.
\]

Take $\tau=\varepsilon_*/(100r^2)$ and choose an even integer $\bar r\geq2/\tau$, we obtain 
\[
\sigma_{2k,\bar r}\leq\sigma_{k,2r}-\tfrac12\varepsilon_*.
\]
Continuing the process show that there exist $k$ and $r$ depending on $\kappa$ and $\sigma'$ such that
\[
E_{k,r}(N)<N^{\sigma'}.
\]

We now turn the small moment into decay at a specific frequency. Fix $|\xi_0|\sim N$. It follows that
\begin{align*}
|\widehat\mu_{k+1}(\xi_0)|^r &\leq\int|\widehat\mu_k(y\xi_0)|^rd\mu(y)\\
 &\lesssim \int_{\frac{1}{2}}^1 |\widehat{\mu}_k(y\mu_0)|^r\bigl(\mu*(N\chi_{I(0,\frac{1}{N})})\bigr)(y)dy\\
 &\lesssim \int_0^N |\widehat{\mu}_k(\xi)|^r\bigl(\mu*(N\chi_{I(0,\frac{1}{N})})\bigr)(\frac{\xi}{|\xi_0|})\frac{d\xi}{|\xi_0|}\\
 &\le \Vert \mu*(N\chi_{I(0,\frac{1}{N})}) \Vert_{\infty}\cdot|\xi_0|^{-1}\cdot\int_0^N|\widehat{\mu}_k(\xi)|^rd\xi\\
 &<N\cdot N^{-\kappa}\cdot N^{-1}\cdot N^{\sigma'}\\
 &=N^{\sigma'-\kappa}.
\end{align*}
Taking $\sigma'=\kappa/2$, and choosing $\varepsilon_2=\kappa/2r$, proves the result.
\end{proof}

\subsection[From Fourier decay to Theorem 0.12.]{From Fourier decay to \paperref{result:0.12}{Theorem 0.12}.}\label{sec:9.3}
For all sufficiently large $N\in\mathbb{Z}_+$ with $N<1/\delta$, clearly the measure $\mu$ satisfies the condition in \paperref{result:9.2}{Theorem~\ref*{result:9.2}}. Hence, if $k$ is sufficiently large, we have
\[
|\widehat\mu_{k}(\xi)|\lesssim|\xi|^{-\varepsilon_2}\qquad(|\xi|<1/\delta).
\]
Note that
\[
\underbrace{\mu*\cdots*\mu}_{k-\text{fold}}=\psi[\underbrace{\mu\times\cdots\times\mu}_{k-\text{fold}}],
\]
where
\[
\phi:\ \mathbb{R}^k\to\mathbb{R},\quad(x_1,\ldots,x_k)\mapsto x_1+\cdots+x_k.
\]
The image measure $\nu$ in \paperref{result:0.12}{Theorem~\ref*{result:0.12}} is $k$-fold convolutions of $\mu_k$, i.e., $\nu=\mu_k*\cdots*\mu_k$. Let $k>2/\varepsilon_2$ yield
\[
|\widehat{\nu}(\xi)|=|\widehat{\mu}_k(\xi)|^k\lesssim|\xi|^{-k\varepsilon_2}<|\xi|^{-2}
\]
for $|\xi|<1/\delta$. Now, for $\delta<\rho<1$,
\begin{align*}
\nu(I(a,\rho)) &=\int \chi_{I(0,\rho)}(x-a)d\nu(x)\\
&\sim\int\rho\widehat{\chi}_{I(0,1)}(\rho\xi)e^{-2\pi ia\xi}\widehat{\nu}(\xi)d\xi\\
&\lesssim\rho\int_{|\xi|<1/\rho}|\widehat\nu(\xi)|d\xi\\
&\lesssim\rho\int_{1<|\xi|<1/\rho}|\xi|^{-2}d\xi\\
&<\rho
\end{align*}
This concludes the proof of \paperref{result:0.12}{Theorem~\ref*{result:0.12}}.
\clearpage
\phantomsection

\end{document}